\documentclass[opre,nonblindrev]{informs3} 

\OneAndAHalfSpacedXI 

\usepackage{natbib}
 \bibpunct[, ]{(}{)}{,}{a}{}{,}%
 \def\bibfont{\footnotesize}%
 \def\BIBand{and}%

\usepackage{float}
\usepackage{graphicx}
\usepackage[hidelinks]{hyperref}
\usepackage{blindtext}
\usepackage{standalone}
\usepackage[misc]{ifsym}
\usepackage{amssymb}
\usepackage{bm}
\usepackage{bbm}
\usepackage{mathtools}
\usepackage{pgfplots}
\usepackage{stmaryrd}
\usepackage{amsmath}
\usepackage{pifont}
\newcommand{\cmark}{\text{\ding{51}}}
\newcommand{\xmark}{\text{\ding{55}}}

\usepackage{mathrsfs}
\usepackage{tikz}

\usetikzlibrary{positioning}

\usetikzlibrary{arrows,automata}
\usetikzlibrary{decorations.text}

\usepackage{booktabs}
\usepackage{algorithm}
\usepackage{algorithmic}
\DeclarePairedDelimiterX{\inp}[2]{\langle}{\rangle}{#1, #2}

\usepackage{esvect}

\usepackage{caption}
\usepackage{subcaption}
\usepackage{array}

\usepackage[compact]{titlesec}
\titlespacing{\section}{0pt}{2ex}{1ex}
\titlespacing{\subsection}{0pt}{1ex}{0ex}
\titlespacing{\subsubsection}{0pt}{0.5ex}{0ex}

\usepackage[utf8]{inputenc}
\usepackage{placeins}

\TheoremsNumberedThrough
\ECRepeatTheorems
\EquationsNumberedThrough

\usepackage{mdframed}
\mdfsetup{skipabove=22pt,skipbelow=22pt}
\newmdenv[
  topline=true,
  bottomline=true,
  rightline=true,
  leftline=true,
  innertopmargin=11pt,
  linewidth=.75pt
]{exbox}

\begin{document}
\RUNAUTHOR{Bertsimas and Cory-Wright}

\RUNTITLE{A Scalable Algorithm For Sparse Portfolio Selection}

\TITLE{A Scalable Algorithm for Sparse Portfolio Selection}

\ARTICLEAUTHORS{%
\AUTHOR{Dimitris Bertsimas}
\AFF{Sloan School of Management, Massachusetts Institute of Technology, Cambridge, MA, USA.\\ ORCID: \href{https://orcid.org/0000-0002-1985-1003}{$0000$-$0002$-$1985$-$1003$}\\ \EMAIL{dbertsim@mit.edu}} 
\AUTHOR{Ryan Cory-Wright}
\AFF{Operations Research Center, Massachusetts Institute of Technology, Cambridge, MA, USA, \\ORCID: \href{https://orcid.org/0000-0002-4485-0619}{$0000$-$0002$-$4485$-$0619$}\\\EMAIL{ryancw@mit.edu}}
} 

\ABSTRACT{%
The sparse portfolio selection problem is one of the most famous and frequently-studied problems in the optimization and financial economics literatures. In a universe of risky assets, the goal is to construct a portfolio with maximal expected return and minimum variance, subject to an upper bound on the number of positions, linear inequalities and minimum investment constraints. Existing certifiably optimal approaches to this problem {\color{black}have not been shown to} converge within a practical amount of time at real-world problem sizes with more than $400$ securities. In this paper, we propose a more scalable approach. By imposing a ridge regularization term, we reformulate the problem as a convex binary optimization problem, which is solvable via an efficient outer-approximation procedure. We propose various techniques for improving the performance of the procedure, including a heuristic which supplies high-quality warm-starts, {\color{black} and }a {\color{black} second heuristic for generating additional cuts which strengthens the root relaxation}. We also study the problem's {\color{black}continuous} relaxation, establish that it is second-order-cone representable, and supply a sufficient condition for its tightness. In numerical experiments, we establish that {\color{black}a conjunction of the imposition of ridge regularization and the use of }the outer-approximation procedure gives rise to dramatic speedups for sparse portfolio selection problems.
}%

\KEYWORDS{Sparse Portfolio Selection, Binary Convex Optimization, Outer Approximation.}
\HISTORY{This paper was first submitted in November $2020$. A revision was submitted in March $2021$.}
\maketitle
 \SUBJECTCLASSname{programming: integer; non-linear: quadratic; finance: portfolio}

\AREAOFREVIEWname{Design and Analysis of Algorithms-Discrete} 
\section{Introduction}
Since the Nobel-prize winning work of \citet{markowitz1952portfolio}, the problem of selecting an optimal portfolio of securities has received an enormous amount of attention from practitioners and academics alike. In a universe containing $n$ distinct securities with expected marginal returns $\bm{\mu} \in \mathbb{R}^n$ and {a }variance-covariance matrix of the returns {\color{black}${\bm{\Sigma} \in S^n_+}$}, the {\color{black} scalarized} Markowitz model selects a portfolio which provides the highest expected return for a given amount of variance, by solving:
\begin{equation}\label{prob:ps}
\begin{aligned}
    \min_{\bm{x} \in \mathbb{R}^n_+}\ \frac{\sigma}{2} \bm{x}^\top \bm{\Sigma} \bm{x}-\bm{\mu}^\top \bm{x} \quad \text{s.t.} \quad \bm{e}^\top \bm{x}=1,
\end{aligned}
\end{equation}
{where $\sigma \geq 0$ is a parameter that controls the trade-off between the portfolios risk and return,} {\color{black}and $\bm{e} \in \mathbb{R}^n$ denotes the vector of all ones.}

To improve its realism, many authors have proposed augmenting Problem \eqref{prob:ps} with minimum investment, maximum investment, and cardinality constraints among others \citep[see, e.g., ][]{jacob1974limited, perold1984large, chang2000heuristics}. Unfortunately, these constraints are disparate and sometimes imply each other, which makes defining a canonical portfolio selection model challenging. We refer to \citet{mencarelli2019complex} for a survey of real-life portfolio selection constraints.

\citet{bienstock1996computational} \citep[see also][]{bertsimas1999portfolio} defined a realistic portfolio selection model by augmenting Problem \eqref{prob:ps} with two sets of inequalities. The first {\color{black} set is a generic system of linear inequalites $\bm{l} \leq \bm{A}\bm{x} \leq \bm{u}$ which, {\color{black}through an appropriate choice of data $\bm{l} \in \mathbb{R}^m, \bm{u} \in \mathbb{R}^m, \bm{A} \in \mathbb{R}^{m \times n}$,} ensures that various real-world constraints such as allocating an appropriate amount of capital to each market sector hold.} The second inequality limits the number of non-zero positions held to {\color{black} $k \ll n$}, by requiring that the portfolio is $k$-sparse, i.e., $\Vert \bm{x} \Vert_0 \leq k$. The sparsity constraint is important because (a) managers incur monitoring costs for each non-zero position, and (b) investors believe that portfolio managers who do not control the number of positions held perform index-tracking while charging active management fees \citep[see][for an implementation of portfolio selection with sparsity constraints at a real-world asset management company]{bertsimas1999portfolio}{\color{black}; see also \citep[][Chap. 11.1]{bertsimas2005optimization} for further discussion on the motivation for setting $k \ll n$}. Imposing the real-world constraints yields the following portfolio selection model:
\begin{align}
\label{mainproblemnonreg}
    \min_{\bm{x} \in \mathbb{R}^n_+}\ \frac{\sigma}{2} \bm{x}^\top \bm{\Sigma} \bm{x}-\bm{\mu}^\top \bm{x} \quad \text{s.t.} \quad  \bm{l} \leq \bm{A}\bm{x} \leq \bm{u},\ \bm{e}^\top \bm{x}=1, \ \Vert \bm{x} \Vert_0 \leq k.
\end{align}
{\color{black}Note that Problem \eqref{mainproblemnonreg} is NP-hard—even in the absence of linear inequalities \citep{gao2013optimal}.}

By introducing binary variables $z_i$ which {\color{black}model whether $x_i$ takes non-zero values by requiring that $x_i=0$ if $z_i=0$, we rewrite the above problem as a {\color{black}convex} mixed-integer quadratic problem:}
\begin{align}
\label{mainproblemnonreg_2}
    \min_{\bm{z} \in \{0, 1\}^n: \bm{e}^\top \bm{z} \leq k, \ \bm{x} \in \mathbb{R}^n_+}\ \frac{\sigma}{2} \bm{x}^\top \bm{\Sigma} \bm{x}-\bm{\mu}^\top \bm{x} \quad \text{s.t.} \quad  \bm{l} \leq \bm{A}\bm{x} \leq \bm{u},\ \bm{e}^\top \bm{x}=1, \ x_i=0 \ \text{if} \ z_i=0 \quad \forall i \in [n].
\end{align}

In the past $20$ years, a number of authors have proposed approaches for solving Problem \eqref{mainproblemnonreg} to certifiable optimality. However, no method {\color{black}has been shown to }scale to real-world problem sizes\footnote{{\color{black}By ``real-world'', we refer to problem instances where the portfolio manager optimizes over an index at least as large as the S$\&$P $500$.} To our knowledge, the Wilshire $5000$ index, which contains around $3,200$ frequently traded securities, is the largest index by number of securities. As portfolio optimization problems generally involve optimizing over securities within an index, we have written $3,200$ here as an upper bound, although one could conceivably also optimize over securities from {\color{black}the union of} multiple stock indices.} where $500 \leq n \leq 3,200$ {\color{black} and $k \ll n$}. This lack of scalability presents a challenge for practitioners and academics alike, because a scalable algorithm for Problem \eqref{mainproblemnonreg} has numerous financial applications, while algorithms which do not scale to this problem size are less practically useful.
\subsection{Problem Formulation and Main Contributions}
In this paper, we provide two main contributions. Our first contribution is augmenting Problem \eqref{mainproblemnonreg} with a ridge regularization term, namely $1/(2\gamma) \cdot \Vert \bm{x}\Vert_2^2$—where $\gamma>0$ is fixed—to yield:
\begin{align}
\label{mainproblem}
    \min_{\bm{x} \in \mathbb{R}^n_+}\ \frac{\sigma}{2} \bm{x}^\top \bm{\Sigma} \bm{x}+\frac{1}{2\gamma}\Vert \bm{x} \Vert_2^2-\bm{\mu}^\top \bm{x} \quad \text{s.t.} \quad  \bm{l} \leq \bm{A}\bm{x} \leq \bm{u},\ \bm{e}^\top \bm{x}=1, \ \Vert \bm{x} \Vert_0 \leq k.
\end{align}

{\color{black}
This problem is more practically tractable\footnote{\color{black}We remark that while Problem \eqref{mainproblem} is more practically tractable than Problem \eqref{mainproblemnonreg}, both problems are NP-hard. Indeed, while the NP-hardness of \eqref{mainproblem} has not been explicitly written down, it follows directly from \cite[Section E.C.1]{gao2013optimal}, because the covariance matrix used in their proof of NP-hardness is positive definite and can therefore be split into a positive semidefinite matrix plus a diagonal regularization matrix as described in Section \ref{sec:cuttingolanemethodcont}.} than Problem \eqref{mainproblemnonreg}, for two reasons. First, {\color{black}as we formally establish in Section \ref{ssec:socp}, the duality gap between Problem \eqref{mainproblem} and its second-order cone relaxation decreases as we decrease $\gamma$ and becomes $0$ at some finite $\gamma>0$. Second, as we numerically establish in Section \ref{sec:compexperiments}, the algorithms developed here converge more rapidly when $\gamma$ is smaller.

In addition to being more practically tractable, Problem \eqref{mainproblem} is a computationally useful surrogate for Problem \eqref{mainproblemnonreg}. Indeed, as we formally establish in Section \ref{ssec:sensitivityanalysis}, any optimal solution to Problem \eqref{mainproblem} is a $1/(2\gamma)$-optimal solution to Problem \eqref{mainproblemnonreg}. Moreover, one can find a solution to Problem \eqref{mainproblemnonreg} which is—often substantially—better than this, by (a) solving Problem \eqref{mainproblem} and (b) solving a simple quadratic optimization problem over the set of securities with the same support as Problem \eqref{mainproblem}'s solution and an unregularized objective. Indeed, {\color{black}since there are finitely many $k$-sparse binary support vectors, this strategy recovers an optimal solution to \eqref{mainproblemnonreg} for any sufficiently large $\gamma$. }}}

Our second main contribution is a scalable outer-approximation algorithm for Problem \eqref{mainproblem}. {\color{black}By utilizing Problem \eqref{mainproblem}'s regularization term, we question the modeling paradigm of writing the logical constraint ``$x_i=0$ if $z_i=0$'' as $x_i \leq z_i$ in Problem \eqref{mainproblem}, by substituting the
equivalent but non-convex term $x_i z_i$ for $x_i$ and invoking strong duality {\color{black}and that $z_i^2=z_i$} to {\color{black}obtain a convex mixed-integer quadratic reformulation of the problem.}}
This allows us to propose a new outer-approximation algorithm {\color{black}in the spirit of the methods of \cite{duran1986outer,fletcher1994solving} but applied to a perspective reformulation \citep{frangioni2006perspective, gunluk2012perspective} of the problem} which solves large-scale sparse portfolio selection problems with up to $3,200$ securities to certifiable optimality {\color{black}in $100$s or $1000$s of seconds}.

{\color{black}
\subsubsection{Connection Between Regularization and Robustness: \quad }
While we have introduced ridge regularization as a device which improves the problems tractability while only marginally affecting the optimal objective value, one can actually interpret the regularizer as a robustification technique which improves the overall quality of the selected portfolio. Indeed, \citet{ledoit2004honey} \citep[see also][]{carrasco2011optimal} have demonstrated that, in mean-variance portfolio selection problems, the largest eigenvalues of the sample covariance matrix $\bm{\Sigma}$ are systematically biased upwards and the smallest eigenvalues of $\bm{\Sigma}$ are systematically biased downwards. As a result, imposing a ridge regularization term (with a properly cross-validated $\gamma$) leads to portfolios which perform better out-of-sample. In a similar vein, \citet{demiguel2009optimal} has shown that the strategy of allocating an identical amount of capital to each security outperforms $13$ other popular investment strategies. Since a ridge regularization term encourages investing a more equal amount in each security, \citet{demiguel2009optimal}'s work suggests that a ridge regularization term is beneficial.
}
{\color{black}
{\color{black}

}}
\subsection{The Scalability of State-of-the-Art Approaches}\label{sec:thescalabilityofsoa}
{\color{black}We now review the scalability of existing state-of-the-art approaches, as stated by their authors, in order to justify our claim that no existing method has been shown to scale to problem sizes where ${500 \leq n \leq 3,200}$ and ${20 \leq k \leq 50}$. In this direction, Table \ref{tab:instancesizesummary} depicts the largest problem solved by each approach, as reported by its authors, and Table \ref{tab:instancedata2} depicts the constraints imposed.

{We remark that, as indicated in Table \ref{tab:instancedata2}, each of the existing approaches was benchmarked on a different data set using a different solver and a different processor, which makes Table \ref{tab:instancesizesummary} an imperfect comparison with our method. Nonetheless, Table \ref{tab:instancesizesummary} appears to be about as accurate a comparison as we can reasonably hope to achieve within the scope of this work. Indeed, a more accurate comparison would require re-implementing these methods from scratch {\color{black} and benchmarking them using the same machine/solver}, which, as {\color{black}we are not aware of publicly available source code for any of these methods other than \citet{vielma2017extended}\footnote{\color{black}Note that the techniques in \cite{vielma2017advances}'s approach were incorporated within CPLEX as of version $12.6.2$ \citep{vielma2017advances}. Accordingly, the numerical comparison against CPLEX $12.8$'s MISOCP solver in Section \ref{sec:compexperiments} can be viewed as a comparison against \citet{vielma2017extended}'s approach.}}, would require a separate review paper.

}
}

\begin{table}[h]
\centering\footnotesize
\caption{Largest portfolio instance solved during benchmarking, by approach. ``$k_{\max}$'' denotes the largest cardinality-constraint right-hand-side imposed when benchmarking. ``n/a'' indicates that a cardinality constraint was not imposed. {\color{black} ``SDO'' denotes a method which involves solving an auxiliary semidefinite optimization problem.}}
\begin{tabular}{@{}r l r r@{}} \toprule
Reference &  Solution method & Largest instance solved & $k_{\max}$\\
& & (no. securities) & \\ \midrule
\citet{vielma2008lifted} & Nonlinear Branch-and-Bound & $100$ & $10$\\
\cite{bonami2009exact} & Nonlinear Branch-and-Bound & $200$ & n/a\\
\cite{frangioni2009computational} & Branch-and-Cut+SDO & $400$ & n/a\\
\cite{gao2013optimal} & Nonlinear Branch-and-Bound & $300$ & $20$\\
\cite{cui2013convex} & Nonlinear Branch-and-Bound & $300$ & $10$\\
\cite{zheng2014improving} & Branch-and-Cut+SDO & $400$ & $12$\\
\citet{frangioni2016approximated} & Branch-and-Cut+SDO  & $400$ & $10$\\
\citet{vielma2017extended} & Branch-and-Bound & $200$ & $10$\\
\citet{frangioni2017improving} & Branch-and-Cut+SDO & $400$ & $10$\\
\bottomrule
\end{tabular}
\label{tab:instancesizesummary}
\end{table}

\begin{table}[h]
    \footnotesize
    \caption{Constraints imposed and solver used, by reference; see also \citep[Table 1]{mencarelli2019complex}. We use the following notation to refer to the constraints imposed: \textbf{C}: A Cardinality constraint $\Vert \bm{x}\Vert_0 \leq k$; \textbf{MR}: A Minimum Return constraint $\bm{\mu}^\top \bm{x}\geq \bar{r}$; \textbf{SC}: A Semi-Continuous, or minimum investment, constraint $x_i \in \{0\} \cup [l_i, u_i] \ \forall i \in [n]$; \textbf{SOC}: A Second-Order-Cone approximation of a chance constraint: $\bm{\mu}^\top \bm{x}+F^{-1}_{\bm{x}}(1-p)\sqrt{\bm{x}^\top \bm{\Sigma}\bm{x}} \geq R$; \textbf{LS}: A Lot-sizing constraint $x_i = M \rho_i: \rho_i \in \mathbb{Z}$.}
    \centering
    \begin{tabular}{@{} r l l r r r r r l@{}}\toprule
    Reference & Solver & & C & MR  & SC & SOC& LS & Data Source\\\midrule
    \citet{vielma2008lifted} & \verb|CPLEX| $10.0$ & & \cmark & \xmark & \xmark & \cmark & \xmark & $20$ instances generated\\
    & & & & & & & & using S$\&$P $500$ daily returns\\
    \citet{bonami2009exact} & \verb|CPLEX| $10.1$, & & {\color{black}\xmark} & \xmark & \cmark & \cmark & \cmark & $36$ instances generated\\
    & \verb|Bonmin| & & & & & & & using S$\&$P $500$ daily returns\\
    \citet{frangioni2009computational} & \verb|CPLEX| $11$ & & \xmark & \cmark & \cmark & \xmark & \xmark & \citet{frangioni2006perspective}\\
    \citet{cui2013convex} & \verb|CPLEX| $12.1$ & & \cmark & \cmark & \cmark & \xmark & \xmark & $20$ self-generated instances\\
    \citet{gao2013optimal} & \verb|CPLEX| $12.3$, & & \cmark & \cmark & \xmark & \xmark & \xmark & $58$ instances generated\\
    & \verb|MOSEK| & & & & & & & using S$\&$P $500$ daily returns\\
    \citet{zheng2014improving} & \verb|CPLEX| $12.4$ & & \cmark & \cmark & \cmark & \xmark & \xmark & \citet{frangioni2006perspective}\\
    & & & & & & & & OR-library \citep{beasley1990or}\\
    \citet{frangioni2016approximated} & \verb|CPLEX| $12.6$ & & \cmark & \cmark & \cmark & \xmark & \xmark & \citet{frangioni2006perspective}\\
    \citet{vielma2017extended} & \verb|CPLEX| $12.6$, & & \cmark & \xmark & \xmark & \cmark & \xmark & $200$ instances generated by
    \\
   & \verb|Gurobi| $5.6$ & & & & & & &\citet{vielma2008lifted}\\
    \citet{frangioni2017improving} & \verb|CPLEX| $12.7$ & & \cmark & \cmark & \cmark & \xmark & \xmark & \citet{frangioni2006perspective}\\
    This paper & \verb|CPLEX| $12.8$ & & \cmark & \cmark & \cmark & \xmark & \xmark & \citet{frangioni2006perspective}\\
    & \verb|MOSEK| & & & & & & & OR-library \cite{beasley1990or}\\
    \bottomrule
    \end{tabular}
        \label{tab:instancedata2}
\end{table}

\subsection{Background and Literature Review}\label{ssec:background}
Our work touches on three different strands of the mixed-integer non-linear optimization literature, each of which propose certifiably optimal methods for solving Problem \eqref{mainproblemnonreg}: (a) branch-and-bound methods which solve a sequence of relaxations, (b) decomposition methods which separate the discrete and continuous variables in Problem \eqref{mainproblemnonreg}, and (c) perspective reformulation methods which obtain tight relaxations by linking the discrete and the continuous in a non-linear manner.

\subsubsection*{Branch-and-bound algorithms: \quad} A variety of branch-and-bound algorithms have been proposed for solving Mixed-Integer Nonlinear Optimization problems (MINLOs) to certifiable optimality, since the work of \citet{glover1975improved}, who proposed linearizing logical constraints ``$x=0$ if $z=0$'' by rewriting them as $-M z \leq x \leq M z$ for some $M>0$. This approach is known as the big-$M$ method.

The first branch-and-bound algorithm for solving Problem \eqref{mainproblemnonreg} to certifiable optimality was proposed by \citet{bienstock1996computational}. This algorithm {does not make use of binary variables. Instead, it} reformulates the sparsity constraint {implicitly, by recursively branching on subsets of the universe of buyable securities and obtaining relaxations by imposing constraints of the form $\sum_i \frac{x_i}{M_i} \leq K$, where $M_i$ is an upper bound on $x_i$.} Similar branch-and-bound schemes {(which make use of binary variables)} are studied in \citet{bertsimas2009algorithm, bonami2009exact}, who solve instances of Problem \eqref{mainproblemnonreg} with up to $50$ (resp. $200$) securities to certifiable optimality. Unfortunately, these methods do not scale well, because reformulating a sparsity constraint via the big-M method often yields weak relaxations in practice\footnote{Indeed, if all securities are i.i.d. then investing $\frac{1}{k}$ in $k$ randomly selected securities constitutes an optimal solution to Problem \eqref{mainproblemnonreg}, but, as proven in \cite{bienstock2010eigenvalue}, branch-and-bound must expand $2^{\frac{n}{10}}$ nodes to improve upon a naive sparsity-constraint free bound by $10\%$, and expand all $2^n$ nodes to certify optimality.}

Motivated by the need to obtain tighter relaxations, more sophisticated branch-and-bound schemes have since been proposed, which obtain higher-quality bounds by lifting the problem to a higher-dimensional space. The first lifted approach was proposed by \citet{vielma2008lifted}, who successfully solved instances of Problem \eqref{mainproblemnonreg} with up to $100$ securities to certifiable optimality, by taking efficient polyhedral relaxations of second order cone constraints. This approach has since been improved by \citet{gao2013optimal,cui2013convex}, who derive non-linear branch-and-bound schemes which use even tighter second order cone and semi-definite relaxations to solve problems with up to $300$ securities to certifiable optimality.


\subsubsection*{Decomposition algorithms: \quad}
A well-known method for solving MINLOs such as Problem \eqref{mainproblemnonreg} is called outer approximation (OA), which was first proposed by \citet{duran1986outer} (building on the work of \citet{kelley1960cutting, benders1962partitioning, geoffrion1972generalized}), who prove its finite termination. OA separates a difficult MINLO into a finite sequence of \textit{master} mixed-integer linear problems and non-linear \textit{subproblems} (NLOs). This is often a good strategy, because linear integer and continuous conic solvers are usually much more powerful than MINLO solvers.

Unfortunately, OA has not yet been successfully applied to Problem \eqref{mainproblemnonreg}, because it requires informative {\color{black}sub}gradient inequalities from each subproblem to attain a fast rate of convergence. Among others, \citet{borchers1997computational, fletcher1998numerical} have compared OA to branch-and-bound, and found that branch-and-bound outperforms OA for Problem \eqref{mainproblemnonreg}. 

In the present paper, by invoking strong duality, we derive a new {\color{black}sub}gradient inequality, redesign OA using this inequality, and solve Problem \eqref{mainproblem} to certifiable optimality via OA. The numerical success of our decomposition scheme can be explained by {\color{black}two} ingredients: (a) the strength of the {\color{black}sub}gradient inequality, and (b) the tightness of our non-linear reformulation of a sparsity constraint, as further investigated in a more general setting in \cite{bertsimas2019unified}.

\subsubsection*{Perspective reformulation algorithms: \quad} An important aspect of solving Problem \eqref{mainproblemnonreg} is understanding its objective's convex envelope, since approaches which exploit the envelope perform better than approaches which use looser approximations of the objective \citep{klotz2013practical}. An important step in this direction was taken by \cite{frangioni2006perspective}, who built on the work of \citet{ceria1999convex} to derive Problem \eqref{mainproblemnonreg}'s convex envelope under an assumption that $\bm{\Sigma}$ is diagonal, and reformulated the envelope as a semi-infinite piecewise linear function. By splitting a generic covariance matrix into a diagonal matrix plus a positive semidefinite matrix, they subsequently derived a class of perspective cuts which provide bound gaps of $<1\%$ for instances of Problem \eqref{mainproblemnonreg} with up to $200$ securities. This approach was subsequently refined by \cite{frangioni2007sdp, frangioni2009computational}, who solved auxiliary semidefinite problems to extract larger diagonal matrices, and thereby solve instances of Problem \eqref{mainproblemnonreg} with up to $400$ securities.

The perspective reformulation approach has also been extended by other authors. An important work in the area is \citet{akturk2009strong}, who, building on the work of \citet[pp. 88, item 5]{ben2001lectures}, prove that if $\bm{\Sigma}$ is positive definite, i.e., $\bm{\Sigma} \succ \bm{0}$, then after extracting a diagonal matrix $\bm{D} \succ \bm{0}$ such that $\sigma\bm{\Sigma}-\bm{D} \succeq \bm{0}$, Problem \eqref{mainproblemnonreg} is equivalent to the following mixed-integer second order cone optimization problem (MISOCO):
\begin{equation}\label{mainproblemmisocp_orig}
    \begin{aligned}
        \min_{\bm{z} \in \mathcal{Z}_k^n, \ \bm{x} \in \mathbb{R}_+^n, \ \bm{\theta} \in \mathbb{R}_+^n}\quad & \frac{\sigma}{2} \bm{x}^\top \bm{\Sigma} \bm{x}+\frac{1}{2}\sum_{i=1}^n D_{i,i} \theta_i-\bm{\mu}^\top \bm{x} \\
    \text{s.t.} \quad &  \bm{l} \leq \bm{A}\bm{x} \leq \bm{u},\ \bm{e}^\top \bm{x}=1, \ x_i^2 \leq \theta_i z_i \quad \forall i \in [n].
    \end{aligned}
\end{equation}

In light of the above MISOCO, a natural question to ask is \textit{what is the best matrix $\bm{D}$ to use?} This question was partially\footnote{Note that weaker continuous relaxations may in fact perform better after branching, {\color{black}as discussed by \cite{dong2016relaxing}.}} answered by \citet{zheng2014improving}, who demonstrated that the matrix $\bm{D}$ which yields the tightest continuous relaxation is computable via semidefinite optimization, and invoked this observation to solve problems with up to $400$ securities to optimality \citep[see also][who derive a similar perspective reformulation of sparse regression problems]{dong2015regularization}. We refer the reader to \citet{gunluk2012perspective} for a survey of perspective reformulation approaches.

\subsubsection*{Connection to our approach: \quad} An unchallenged assumption in \textit{all} perspective reformulation approaches is that Problem \eqref{mainproblemnonreg} \textit{must not be modified}. Under this assumption, perspective reformulation approaches separate $\bm{\Sigma}$ into a diagonal matrix $\bm{D} \succeq \bm{0}$ plus a positive semidefinite matrix $\bm{H}$, such that $\bm{D}$ is as diagonally dominant as possible. Recently, this approach was challenged by \citet{BvP_AS_2016}. Following a standard statistical learning theory paradigm, they imposed a ridge regularizer and set $\bm{D}$ equal to $1/\gamma \cdot \mathbb{I}$, where $\mathbb{I}$ denotes an identity matrix of appropriate dimension. Subsequently, they derived a cutting-plane method which exploits the regularizer to solve large-scale sparse regression problems to certifiable optimality. In the present paper, we join \citet{BvP_AS_2016} in imposing a ridge regularizer, and derive a cutting-plane method which solves convex MIQOs \textit{with constraints}. {\color{black} We also unify their approach with the perspective reformulation approach, in two steps. First, we note that \citet{BvP_AS_2016}'s algorithm can be improved by setting $\bm{D}$ equal to $1/\gamma\cdot \mathbb{I}$ \textit{plus} a perspective reformulation's diagonal matrix, and this is particularly effective when $\bm{\Sigma}$ is diagonally dominant (see Section \ref{sec:cuttingolanemethodcont}, \ref{sec:buyin}){. Second, we observe that the cutting-plane approach also helps solve the unregularized problem, indeed, as mentioned previously it successfully supplies a $1/(2\gamma)$-optimal solution to Problem \eqref{mainproblemnonreg}.}}

\subsection{Structure}\label{sec:structure}
The rest of this paper is laid out as follows:
\begin{itemize}
    \item {In Section \ref{sec:dualperspective}, we
    lay the groundwork for our approach, by observing an equivalence between regression and portfolio selection, and rewriting Problem \eqref{mainproblem} as a constrained regression problem.}
    \item In Section \ref{sec:oamethod}, we propose an efficient numerical strategy for solving Problem \eqref{mainproblem}. By
    observing that Problem \eqref{mainproblem}'s inner dual problem supplies subgradients with respect to the positions held, we design an outer-approximation procedure which solves Problem \eqref{mainproblem} to provable optimality. We also discuss practical aspects of the procedure, including a computationally efficient subproblem strategy and a prepossessing technique for decreasing the bound gap at the root node. {\color{black} In addition, we {\color{black}study the problems sensitivity to $\gamma$}, and {\color{black}establish theoretically} that the support of an optimal portfolio (although not the amount allocated to each security) is stable under small changes in $\gamma$.}
    \item In Section \ref{sec:improvedcuttingplanemethod}, we propose techniques for obtaining certifiably near-optimal solutions quickly. First, we introduce a heuristic which supplies high-quality warm-starts. Second, we observe that Problem \eqref{mainproblem}'s continuous relaxation supplies a near-exact second-order cone representable lower bound, and exploit this observation by deriving a sufficient condition for the bound to be exact.
    \item In Section \ref{sec:compexperiments}, we apply the cutting-plane method to the problems described in \cite{chang2000heuristics}, \cite{frangioni2006perspective}, and three larger scale data sets: the S$\&$P $500$, Russell $1000$, and Wilshire $5000$. We also explore Problem \eqref{mainproblem}'s sensitivity to its hyperparameters, and establish empirically that
    optimal support indices tend to be stable for reasonable hyperparameter choices.
\end{itemize}
\subsubsection*{Notation: \quad }
We let nonbold face characters denote scalars, lowercase bold faced characters such as $\bm{x} \in \mathbb{R}^n$ denote vectors, uppercase bold faced characters such as $\bm{X} \in \mathbb{R}^{n \times r}$ denote matrices, and calligraphic uppercase characters such as $\mathcal{X}$ denote sets. We let $\mathbf{e}$ denote a vector of all $1$'s, $\bm{0}$ denote a vector of all $0$'s, and $\mathbb{I}$ denote the identity matrix, with dimension implied by the context. If $\bm{x}$ is a {\color{black}$n$}-dimensional vector then $\mathrm{Diag}(\bm{x})$ denotes the $n \times n$ diagonal matrix whose diagonal entries are given by $\bm{x}$. We let $[n]$ denote the set of running indices $\{1, ..., n\}$, $\bm{x}\circ\bm{y}$ denote the elementwise, or Hadamard, product between two vectors $\bm{x}$ and $\bm{y}$, and $\mathbb{R}_+^n$ denote the $n$-dimensional non-negative orthant. We let $\mathrm{relint}(\mathcal{X})$ denote the relative interior of a convex set $\mathcal{X}$, i.e., the set of points {\color{black}in} the interior of the affine hull of $\mathcal{X}$ \citep[see][Section 2.1.3]{boyd2004convex}. Finally, we let $\mathcal{Z}_k^n$ denote the set of $k$-sparse binary vectors, i.e, $\mathcal{Z}_k^n:=\{\bm{z} \in \{0, 1\}^n: \ \bm{e}^\top \bm{z} \leq k\}.$

\section{Equivalence Between Portfolio Selection and Constrained Regression}{
\label{sec:dualperspective}
We now lay the groundwork for our outer-approximation procedure, by rewriting Problem \eqref{mainproblem} as a constrained sparse regression problem. 
To achieve this, we take a Cholesky decomposition of $\bm{\Sigma}$ and complete the square. This is justified, because $\bm{\Sigma}$ is positive semidefinite and rank-$r$, meaning there exists an $\bm{X} \in \mathbb{R}^{r \times n}:$ $\bm{\Sigma}=\bm{X}^\top \bm{X}$. Therefore, by scaling $\bm{\Sigma} \leftarrow \sigma\bm{\Sigma}$ and letting: \begin{align}
    \bm{y}:=&\left(\bm{X} \bm{X}^\top\right)^{-1}\bm{X}\bm{\mu},\\ \bm{d}:=&\left(\bm{X}^\top \left(\bm{X} \bm{X}^\top\right)^{-1}\bm{X}-\mathbb{I}\right)\bm{\mu},
\end{align}
be the projection{\color{black}s} of the return vector $\bm{\mu}$ onto the span and nullspace of $\bm{X}$ {\color{black}respectively}, completing the square yields the following equivalent problem, where we add the constant $\frac{1}{2}\bm{y}^\top \bm{y}$ without loss of generality:
\begin{equation}
\label{eqn:constrainedregression}
\begin{aligned}
\min_{\bm{x} \in \mathbb{R}^n_+} \quad & \frac{1}{2\gamma}\left\Vert \bm{x} \right\Vert_2^2+\frac{1}{2} \left\Vert \bm{X}\bm{x}-\bm{y} \right\Vert_2^2+\bm{d}^\top \bm{x} \quad \text{s.t.} \quad  \bm{l} \leq \bm{A}\bm{x} \leq \bm{u}, \
 \bm{e}^\top \bm{x}=1,\ \Vert \bm{x} \Vert_0 \leq k.
\end{aligned}
\end{equation}
{That is, sparse portfolio selection and sparse regression with constraints are equivalent. }

{We remark that while this connection has not been explicitly noted in the literature, Cholesky decompositions of covariance matrices have successfully been applied in perspective reformulation approaches \citep[see][pp 233]{frangioni2006perspective},} {\color{black} and more recently for robust portfolio selection problems with transaction costs \citep[see][]{olivares2018robust}.}
}

\section{A Cutting-Plane Method}\label{sec:oamethod}
In this section, we present an efficient outer-approximation method for solving Problem \eqref{mainproblem}{\color{black}, via its reformulation \eqref{eqn:constrainedregression}}. The key step in deriving this method is enforcing the logical constraint \textit{$x_i=0$ if $z_i=0$} in a tractable fashion. We achieve this by replacing $x_i$ with $z_i x_i$, and rewriting \eqref{mainproblem} as:
\begin{equation}\label{bilevelmainproblem}
\begin{aligned}
    \min_{\bm{z} \in \mathcal{Z}_k^n}\Big[f(\bm{z})\Big],\\
\end{aligned}
\end{equation}
where
\begin{align}\label{bilevelsubproblem}
    f(\bm{z}):=\min_{\bm{x} \in \mathbb{R}^{n}}\quad &\frac{1}{2\gamma}\bm{x}^\top \bm{x}+\frac{1}{2}\left\Vert \bm{X}\bm{Z} \bm{x} -\bm{y} \right\Vert_2^2+\bm{d}^\top \bm{Z}\bm{x} & \text{s.t.} \quad & \bm{l} \leq \bm{A}\bm{Z}\bm{x} \leq \bm{u}, \ \bm{e}^\top\bm{Z} \bm{x}=1,\ \bm{Z}\bm{x} \geq \bm{0},
\end{align} $\bm{Z}=\mathrm{Diag}(\bm{z})$ is a diagonal matrix such that $Z_{i,i}=z_i$, {\color{black}and we do not associate a $\bm{Z}$ term with $\frac{1}{2\gamma}\bm{x}^\top \bm{x}$ in order to enforce the logical constraint $x_i=0$ if $z_i=0$. Observe that the subproblem generated by $f(\bm{z})$ attains its minimum by setting $x_i=0$ whenever $z_i=0$, which verifies that we have successfully enforced the logical constraint.}

{\color{black}
\begin{remark}
The formulation proposed in \eqref{bilevelmainproblem}-\eqref{bilevelsubproblem} is reminiscent of a classical complementarity formulation of a cardinality constraint \citep[see, e.g.,][]{burdakov2016mathematical}, where the logical constraint is enforced directly by imposing  $x_i(1-z_i)=0 \ \forall i \in [n]$. Perhaps the main point of difference between \eqref{bilevelmainproblem} and a complementarity formulation is that this paper effectively replaces the logical constraint with a penalty term in the objective and uses the strongly convex regularizer to implicitly enforce the logical constraints, while the complementarity formulation explicitly imposes the logical constraints, which does not, by itself, give rise to new algorithms for solving {\color{black}\eqref{mainproblem}}.
\end{remark}
}

Problem \eqref{bilevelmainproblem}'s formulation might appear to be intractable, because it appears to {\color{black}define $f(\bm{z})$ as a function which is non-convex in $\bm{z}$.} However, $f$ is actually convex in $\bm{z}$. Indeed, in {\color{black}Theorem \ref{minimax}, which can be found in Section \ref{sec:cuttingplanealgorithm},} we invoke duality to demonstrate that $f(\bm{z})$ can be rewritten as the supremum of functions which are linear in $\bm{z}$, thus proving that $f$ is actually convex in $\bm{z}$.

As $f(\bm{z})$ is convex in $\bm{z}$, a natural strategy for solving \eqref{bilevelmainproblem} is to iteratively minimize and refine a piecewise linear underestimator of $f(\bm{z})$. This strategy is called outer-approximation (OA), and was originally proposed by \citet{duran1986outer}. OA works as follows: by assuming that at each iteration $t$ we have access to $f(\bm{z}_i)$ and a set of subgradients $\bm{g}_{\bm{z}_i}$ at the points $\bm{z}_i: i \in [t]$, we construct the following underestimator:
\begin{align*}
    f_t(\bm{z})=\max_{1 \leq i \leq t} \left\{f(\bm{z}_i)+\bm{g}_{\bm{z}_i}^\top (\bm{z}-\bm{z}_i)\right\}.
\end{align*}
By iteratively minimizing $f_t(\bm{z})$ over $\mathcal{Z}_n^k$ to obtain $\bm{z}_t$, and evaluating $f(\cdot)$ and its subgradient at $\bm{z}_t$, we obtain a non-decreasing sequence of underestimators $f_t(\bm{z}_t)$
which converge to the optimal value of $f(\bm{z})$ within a finite number of iterations, since $\mathcal{Z}_n^k$ is a finite set and OA never visits a point twice. Additionally, we can avoid solving a different {\color{black}mixed-integer linear optimization problem} (MILO) at each OA iteration by integrating the entire algorithm within a single branch-and-bound tree, as first proposed by \cite{padberg1991branchcut, quesada1992lp}, using {\color{black} dynamic cut generation, which could be implemented using a cut pool and \verb|lazy constraint callbacks|
; see \citet{fischetti2016redesigning} for an implementation of this approach for facility location problems.} {\color{black}Indeed, \verb|lazy constraint callbacks| are now standard components of modern MILO solvers such as \verb|Gurobi|, \verb|CPLEX| and \verb|GLPK| which allow users to easily implement single-tree OA schemes—we make use of these callbacks when we develop our single-tree OA scheme in Section \ref{sec:cuttingolanemethodcont}.}

As we mentioned in Section \ref{ssec:background}, OA is a widely known procedure. However, it has not yet been successfully applied to {\color{black}sparse portfolio selections problems}, largely because of a lack of efficient separation oracles which provide both zeroth and first order information {\color{black}concerning the value of $f$ at the point $\bm{z}$}. Therefore, we now outline such a procedure {\color{black}which addresses Problem \eqref{mainproblem}}.
\subsection{Efficient Subgradient Evaluations}
\label{sec:cuttingplanealgorithm}
We now rewrite Problem \eqref{bilevelmainproblem} as a saddle-point problem, in the following theorem:
\begin{theorem}
\label{minimax}
Suppose that Problem \eqref{bilevelmainproblem} is feasible. Then, it is equivalent to the following problem:
\begin{equation}\label{saddlepointreformulation}
\begin{aligned}
    \min_{\bm{z} \in \mathcal{Z}_k^n} \ \max_{\substack{\bm{\alpha} \in \mathbb{R}^r, \ \bm{w} \in \mathbb{R}^n, \\ \bm{\beta}_l,\ \bm{\beta}_u \in \mathbb{R}^m_+,\ \lambda \in \mathbb{R}}} \quad & -\frac{1}{2} \bm{\alpha}^\top \bm{\alpha} - \frac{\gamma}{2} \sum_i z_i w_i^2 +\bm{y}^\top \bm{\alpha} +\bm{\beta}_l^\top \bm{l}-\bm{\beta}_u^\top \bm{u}+\lambda\\
    {\mbox{\rm s.t.}} \quad & \bm{w}\geq \bm{X}^\top \bm{\alpha}+\bm{A}^\top\left(\bm{\beta}_l-\bm{\beta}_u\right)+\lambda \bm{e}-\bm{d}.
\end{aligned}
\end{equation}
\end{theorem}

\begin{remark}
Theorem \ref{minimax} proves that $f(\bm{z})$ is convex in $\bm{z}$, by rewriting $f(\bm{z})$ as the pointwise maximum of functions which are linear in $\bm{z}$ \citep[][Section 3.2.3]{boyd2004convex}. This justifies our application of OA, which converges under a convexity assumption over finite sets such as $\mathcal{Z}^n_k$. 
\end{remark}


\begin{remark}
In the proof of Theorem \ref{minimax}, we obtain the following relationship between the optimal primal and dual variables:
\begin{gather}\label{primaldualkkt}
    \bm{x}^*=\gamma\mathrm{Diag}(\bm{z}^\star)\bm{w}^\star.
\end{gather}
Notably, this proof carries through if we replace $\bm{z} \in \mathcal{Z}^k_n$ with $\mathrm{Conv}(\mathcal{Z}_k^n)$. Therefore, Equation \ref{primaldualkkt} also holds on $\mathrm{int}(\mathcal{Z}_k^n)$.
\end{remark}

{\color{black}
\proof{Proof of Theorem \ref{minimax}}
We prove this result by invoking strong duality. First, observe that, for each fixed $\bm{z} \in \mathcal{Z}_n^k$, Problem \eqref{bilevelsubproblem}'s dual maximization problem \eqref{saddlepointreformulation} is feasible, because $\bm{w}$ can be increased without bound. {\color{black}Moreover, each inner primal minimization problem is a convex quadratic with linear constraints, and thus satisfies the Abadie constraint qualification whenever it is feasible \citep{abadie1967kuhn}.} Therefore, for each fixed $\bm{z} \in \mathcal{Z}_n^k$, either the inner optimization problem (with respect to $\bm{x}$) is infeasible or strong duality holds.

Let us now introduce an auxiliary vector of variables $\bm{r}$ such that $\bm{r}=\bm{y}-\bm{X}\bm{Z}\bm{x}$. This allows us to rewrite Problem \eqref{bilevelmainproblem} as {\color{black} the following problem, where each set of constraints is associated with a vector of dual variables in square brackets}
\begin{equation}
\label{primalwithresidual2}
\begin{aligned}
\min_{\bm{z} \in \mathcal{Z}_n^k, \ \bm{x} \in \mathbb{R}^n,\ \bm{r} \in \mathbb{R}^r} \quad & \frac{1}{2\gamma}\Vert \bm{x} \Vert_2^2+\frac{1}{2} \Vert \bm{r} \Vert_2^2+\bm{d}^\top \bm{Z} \bm{x} \\
\text{s.t.} \quad & \bm{y}-\bm{X}\bm{Z}\bm{x}=\bm{r}, \ [\bm{\alpha}], \ \bm{A}\bm{Z}\bm{x} \geq \bm{l}, \ [\bm{\beta}_l], \ \bm{A}\bm{Z}\bm{x} \leq \bm{u},  \ [\bm{\beta}_u], \\
 & \bm{e}^\top \bm{Z}\bm{x}=1, \ [\lambda],\ \bm{Z} \bm{x}\geq \bm{0}, \ [\bm{\pi}],\\
\end{aligned}
\end{equation}
{\color{black} where $\bm{\alpha} \in \mathbb{R}^r, \ \bm{w} \in \mathbb{R}^{n},\ \bm{\beta}_l,\ \bm{\beta}_u \in \mathbb{R}^m_+,\ \lambda \in \mathbb{R}, \bm{\pi} \in \mathbb{R}^n_+$.}

This problem has the following Lagrangian:
\begin{align*}
    \mathcal{L}\quad =& \frac{1}{2\gamma} \bm{x}^\top \bm{x} + \frac{1}{2} \bm{r}^\top \bm{r}+\bm{d}^\top \bm{Z}\bm{x}+\bm{\alpha}^\top \left(\bm{y}-\bm{X} \bm{Z} \bm{x}-\bm{r}\right)-\bm{\pi}^\top \bm{Z} \bm{x}\\
    &-\lambda(\bm{e}^\top \bm{Z}\bm{x}-1)
    -\bm{\beta}_l^\top(\bm{A}\bm{Z}\bm{x}-\bm{l})+\bm{\beta}_u^\top(\bm{A}\bm{Z}\bm{x}-\bm{u}).
\end{align*}

For a fixed $\bm{z}$, minimizing this Lagrangian is equivalent to solving the following KKT conditions:
\begin{align*}
    \nabla_{\bm{x}} \mathcal{L}=\bm{0}& \implies \frac{1}{\gamma}\bm{x}+\bm{Z}\left(\bm{d}-\bm{X}^\top \bm{\alpha}-\bm{\pi}-\lambda\bm{e}- \bm{A}^\top(\bm{\beta}_l-\bm{\beta}_u)\right)=\bm{0}, \\
    &\implies \bm{x}=\gamma \bm{Z} \left(\bm{X}^\top \bm{\alpha}+ \bm{\pi}+\lambda\bm{e}+ \bm{A}^\top(\bm{\beta}_l-\bm{\beta}_u)-\bm{d}\right),\\
    \nabla_{\bm{r}} \mathcal{L}=\bm{0}&  \implies \bm{r}-\bm{\alpha}=0 \implies \bm{r}=\bm{\alpha}.
\end{align*}
Substituting the above expressions for $\bm{x}$, $\bm{r}$ into $\mathcal{L}$ then defines the Lagrangian dual, where we eliminate $\bm{\pi}$ {\color{black}(by replacing it with a non-negativity constraint)} and introduce $\bm{w}$ such that $\bm{x}:=\gamma \bm{Z}\bm{w}$ for brevity. The Lagrangian dual reveals that for any $\bm{z}$ such that Problem \eqref{bilevelsubproblem} is feasible:
\begin{align*}
    f(\bm{z})=\max_{\substack{\bm{\alpha} \in \mathbb{R}^r, \ \bm{w} \in \mathbb{R}^{n},\\ \bm{\beta}_l,\ \bm{\beta}_u \in \mathbb{R}^m_+,\ \lambda \in \mathbb{R}}} \quad & -\frac{1}{2} \bm{\alpha}^\top \bm{\alpha} - \frac{\gamma}{2} \bm{w}^\top \bm{Z}^2 \bm{w} +\bm{y}^\top \bm{\alpha} +\bm{\beta}_l^\top \bm{l}-\bm{\beta}_u^\top \bm{u}+\lambda\\
    \text{s.t.} \quad & \bm{w}\geq \bm{X}^\top \bm{\alpha}+\bm{A}^\top\left(\bm{\beta}_l-\bm{\beta}_u\right)+\lambda\bm{e}-\bm{d}.
\end{align*}

Moreover, at binary points $\bm{z}$, $z_i^2=z_i$ and therefore the above problem is equivalent to solving:
\begin{align*}
    f(\bm{z})=\max_{\substack{\bm{\alpha} \in \mathbb{R}^r, \ \bm{w} \in \mathbb{R}^n, \\ \bm{\beta}_l,\ \bm{\beta}_u \in \mathbb{R}^m_+,\ \lambda \in \mathbb{R}}} \quad & -\frac{1}{2} \bm{\alpha}^\top \bm{\alpha} - \frac{\gamma}{2} \sum_i z_i w_i^2 +\bm{y}^\top \bm{\alpha} +\bm{\beta}_l^\top \bm{l}-\bm{\beta}_u^\top \bm{u}+\lambda\\
    \text{s.t.} \quad & \bm{w} \geq \bm{X}^\top \bm{\alpha}+\bm{A}^\top\left(\bm{\beta}_l-\bm{\beta}_u\right)+\lambda \bm{e}-\bm{d}.
\end{align*}

Minimizing $f(\bm{z})$ over $\mathcal{Z}_k^n$ then yields the result, where we ignore choices of $\bm{z}$ which yield infeasible primal subproblems without loss of generality, as their dual problems are feasible and therefore unbounded by weak duality, and a choice of $\bm{z}$ such that $f(\bm{z})=+\infty$ is certainly suboptimal.  \Halmos
\endproof}

Theorem \ref{minimax} supplies objective function evaluations $f(\bm{z}_t)$ and subgradients $\bm{g}_t$ after solving a single convex quadratic optimization problem. We formalize this observation in the following corollary:

\begin{corollary}\label{subgradientvalues}
Let $\bm{w}^\star(\bm{z})$ be an optimal choice of $\bm{w}$ for a particular subset of securities $\bm{z}$. {\color{black} Then, a valid subgradient $\bm{g}_{\bm{z}} \in \partial f(\bm{z})$ has components given by the following expression for each $i \in [n]$:}
\begin{align}\label{subgradexpression}
    g_{\bm{z},i}=-\frac{\gamma}{2}w_i^{\star}(\bm{z})^2.
\end{align}
\end{corollary}

In latter sections of this paper, we design aspects of our numerical strategy by assuming that $f(\bm{z})$ is Lipschitz continuous in $\bm{z}$. It turns out that this assumption is valid whenever we can bound\footnote{Such a bound always exists, because $\mathcal{Z}^n_k$ has a finite number of elements. However, it will, in general, depend upon the problem data. } $\vert\bm{w}_i^\star(\bm{z})\vert$ for each $\bm{z}$, as we now establish in the following corollary ({\color{black}proof due to the well-known subgradient inequality \citep[see][]{boyd2004convex},} deferred to Appendix \ref{sec:ommitedproofs}):
\begin{corollary}\label{corr:lcont}
Let $\bm{w}^\star(\bm{z})$ be an optimal choice of $\bm{w}$ for a given subset of securities $\bm{z} \in \mathcal{Z}^n_k$. Then: \begin{align}
    f(\bm{z})-f(\hat{\bm{z}}) \leq \frac{\gamma}{2} \sum_i (\hat{z}_i-z_i)w_i^\star(\bm{z})^2.
\end{align}
\end{corollary}

\subsection{A Cutting-Plane Method-Continued}\label{sec:cuttingolanemethodcont}
Corollary \ref{subgradientvalues} shows that evaluating $f(\hat{\bm{z}})$ yields a first-order underestimator of $f(\bm{z})$, namely
\begin{align}
f(\bm{z}) \geq f(\hat{\bm{z}})+\bm{g}_{\hat{\bm{z}}}^\top(\bm{z}-\hat{\bm{z}})
\end{align} at no additional cost. Consequently, a numerically efficient strategy for minimizing $f(\bm{z})$ is the previously discussed OA method. We formalize this procedure in Algorithm \ref{alg:cuttingPlaneMethod}. Note that we add the OA cuts via {\color{black}dynamic cut generation} to maintain a single tree of partial solutions throughout the process, and avoid the cost otherwise incurred in rebuilding the tree whenever a cut is added.

\begin{algorithm*}
\caption{An outer-approximation method for Problem \eqref{mainproblem}}
\label{alg:cuttingPlaneMethod}
\begin{algorithmic}
\REQUIRE Initial solution $\bm{z}_1$
\STATE $t \leftarrow 1 $
\REPEAT
\STATE Compute $\bm{z}_{t+1}, \theta_{t+1}$ solution of
{\vspace{-5mm}
\begin{align*}
\min_{\bm{z} \in \mathcal{Z}_k^n, \theta} \: \theta \quad \mbox{ s.t. } \quad \theta \geq f(\bm{z}_i) + \bm{g}_{\bm{z}_i}^\top (\bm{z}-\bm{z}_i) \quad \forall i \in [t],
\end{align*}}\vspace{-5mm}
\STATE Compute $f(\bm{z}_{t+1})$ and $\bm{g}_{\bm{z}_{t+1}} \in \partial f(\bm{z}_{t+1})$
\STATE $t \leftarrow t+1 $
\UNTIL{$ f(\bm{z}_t)-\theta_t \leq \varepsilon$}
\RETURN $\bm{z}_t$
\end{algorithmic}
\end{algorithm*}
{\color{black}
As noted by \citet[Theorem 3]{duran1986outer} and \citet[Theorem 2]{fletcher1994solving}, Algorithm \ref{alg:cuttingPlaneMethod} terminates with an optimal solution to Problem \eqref{mainproblem} in a finite number of iterations, because $f(\bm{z})$ is convex in $\bm{z}$, $\mathcal{Z}_k^n$ is a finite set and Algorithm \ref{alg:cuttingPlaneMethod} never selects a binary vector twice. Interestingly, while OA methods which leverage strong duality typically require a constraint qualification assumption to guarantee correctness, the proof of Theorem \ref{minimax} reveals that this is not necessary to ensure the correctness of Algorithm \ref{alg:cuttingPlaneMethod}. Rather, since each subproblem generated by $f(\bm{z})$ is a convex quadratic with linear constraints, the Abadie constraint qualification \citep{abadie1967kuhn} holds automatically, and thus Algorithm \ref{alg:cuttingPlaneMethod} converges within ${n \choose k}\leq \left(\frac{e n}{k}\right)^k$ iterations in the worst case, and typically far fewer iterations in practice.
}

As Algorithm \ref{alg:cuttingPlaneMethod}'s rate of convergence depends heavily upon its implementation, we now discuss some practical aspects of the method:
\begin{itemize}
    \item \textbf{A computationally efficient subproblem strategy: } For computational efficiency purposes, we would like to solve subproblems which only involve active indices, i.e., indices where $z_i=1$, since $k \ll n$.
    At a first glance, this does not appear to be possible, because we must supply an optimal choice of $w_i$ for all $n$ indices in order to obtain valid subgradients. Fortunately, we can in fact supply a full OA cut after solving a subproblem in the active indices, by exploiting the structure of the saddle-point reformulation. Specifically, we optimize over the $k$ indices where $z_i=1$ and set $w_i=\max\left(\bm{X}_i^\top \bm{\alpha}^\star +\bm{A}_i^\top(\bm{\beta}_l^\star-\bm{\beta}_u^\star) +\lambda^\star -d_i,0\right)$ for the remaining $n-k$ $w_i$'s. This procedure yields an optimal choice of $w_i$ for each index $i$, because it is a feasible choice and the remaining $w_i$'s have a weight of $0$ in the objective function.

    Observe that this procedure yields the strongest possible cut which can be generated at a given $\bm{z}$ \textit{for this set of dual variables}, because (a) the procedure yields the minimum feasible absolute magnitude of $w_i$ whenever $z_i=0$, and (b) there is a unique optimal choice of $w_i$ for the remaining indices, since Problem \eqref{saddlepointreformulation} is strongly concave in $w_i$ when $z_i>0$. In fact, if $w_i=0$ is feasible for some security $i$ such that $z_i=0$, then we cannot improve upon the current iterate $\bm{z}$ by setting $z_i=1$, as our lower approximation gives $$f(\bm{z}+\bm{e}_i) \geq  f(\bm{z})+\bm{g}_{\bm{z}}^\top (\bm{z}+\bm{e}_i-\bm{z})=f(\bm{z}) +\bm{g}_{\bm{z}}^\top \bm{e}_i\geq f(\bm{z}),$$ {where the last inequality holds because $g_{z,i}=-\gamma/2 \cdot w_i^2=0$.}

    \item \textbf{Cut generation at the root node: } Another important aspect of decomposition methods is supplying as much information\footnote{{\color{black} For completeness, we should note that while supplying this information often substantially accelerates decomposition schemes, it can sometimes do more harm than good, particularly if the cuts applied are nearly parallel to one another. In the context of sparse portfolio selection, our numerical experiments in Section \ref{sec:compexperiments} demonstrate the efficacy of supplying information at the root node for all but the simplest sparse porfolio selection problems.}} as possible to the solver before commencing branching, as advocated by \citet[Section 4.2]{fischetti2016redesigning}. One effective way to achieve this is to relax the integrality constraint $\bm{z} \in \{0, 1\}^n$ to $\bm{z} \in [0, 1]^n$ in Problem \eqref{saddlepointreformulation}, run a cutting-plane method on this continuous relaxation and apply the resulting cuts at the root node before solving the binary problem. Traditionally, this relaxation is solved using \citet{kelley1960cutting}'s cutting-plane method.
    However, as \citet{kelley1960cutting}'s method often converges slowly in practice, we instead solve the relaxation using an \verb|in-out| bundle method \citep{ben2007acceleration, fischetti2016redesigning}. We supply pseudocode for our implementation of the \verb|in-out| method in Appendix \ref{sec:auxpsueodocode}.

    In order to further accelerate OA, a variant of the root node processing technique which is often effective is to also run the \verb|in-out| method at some additional nodes in the tree, as proposed in \cite[Section 4.3]{fischetti2016redesigning}. This can be implemented {\color{black} using the cut separation method derived in the previous section (and implemented }via a \verb|user cut callback|), by using the current LO solution $\bm{z}^\star$ as a stabilization point for the \verb|in-out| method.

    To avoid generating too many cuts at continuous points, we impose a limit of $200$ cuts at the root node, $20$ cuts at all other nodes, and do not run the \verb|in-out| method at more than $50$ nodes. One point of difference in our implementation of the \verb|in-out| method \citep[compared to][]{ben2007acceleration, fischetti2016redesigning} is that we use the {\color{black}true} optimal solution to Problem \eqref{saddlepointreformulation}'s continuous relaxation as a stabilization point at the root node {\color{black}(as opposed to the methods current estimate of the optimal solution)}-- this speeds up convergence {\color{black}of the \verb|in-out| method} greatly, and comes at the low price of solving an SOCO to elicit the stabilization point (we obtain the point by solving Problem \eqref{mainproblemmisocp_dualversion}; see Section \ref{ssec:socp}). {\color{black}Note that this approach may appear somewhat counter-intuitive, as it involves using the optimal solution to the relaxation to ``solve'' the relaxation. However, it is actually not, because we are interested in generating a small number of cuts which efficiently describe the surface of the relaxation near its optimal solution, rather than solving the relaxation.}

    {\color{black}Restarting mechanisms, as explored by \citet{fischetti2016benders}, may also be useful for further improving the root node relaxation, although we do not consider them in this work.}
    \item \textbf{Feasibility cuts:} The constraints $\bm{l} \leq \bm{A} \bm{x} \leq \bm{u}$ may render some vectors $\hat{\bm{z}}$ infeasible. In this case, we add the following feasibility cut which bans $\bm{\hat{z}}$ from appearing in future iterations of OA:
\begin{align}
    \sum_i \hat{z}_i (1-z_i)+\sum_i (1-\hat{z}_i)z_i \geq 1.
\end{align}
    An alternative approach is to derive constraints on $\bm{z}$ which ensure that OA never selects an infeasible $\bm{z}$. For instance, if the only constraint on $\bm{x}$ is a minimum return constraint $\bm{\mu}^\top \bm{x} \geq r$ then imposing $\sum_{i: \mu_i \geq r} z_i \geq 1$ ensures that only feasible $\bm{z}$'s are selected. Whenever eliciting these constraints is possible, we recommend imposing them, to avoid infeasible subproblems entirely.
\item \textbf{Extracting Diagonal Dominance:} In problems where $\bm{\Sigma}$ is diagonally dominant {\color{black}in the sense of \cite{barker1975cones}, i.e., $\Sigma_{i,i} \geq \sum_{j \neq i}\vert \Sigma_{i,j}\vert \ \forall i \in [n]$}, the performance of Algorithm \ref{alg:cuttingPlaneMethod} can often be substantially improved by \textit{boosting} the regularizer, i.e., selecting a diagonal matrix $\bm{D} \succeq \bm{0}$ such that $\sigma\bm{\Sigma}-\bm{D} \succeq \bm{0}$, replacing $\sigma\bm{\Sigma}$ with $\sigma\bm{\Sigma}-\bm{D}$, and using a different regularizer $\gamma_i:=\left(\frac{1}{\gamma}+D_{i,i}\right)^{-1}$ for each index $i${\color{black}, to obtain a more general problem of which \eqref{mainproblem} is a special case; note that the algorithms developed here remain valid when we have a different regularizer for each index $i$}. In general, selecting such a $\bm{D}$ involves solving a {\color{black}semidefinite optimization problem} (SDO) \citep[][]{frangioni2007sdp, zheng2014improving}, {\color{black} which is fast when $n$ is in the hundreds, but requires a prohibitive amount of memory when $n$ is in the thousands.} In the latter case, we recommend taking a second-order cone inner approximation of the {\color{black}semidefinite} cone and improving the approximation via column generation. {\color{black}Indeed, this approach has recently been shown to provide high-quality solutions to large-scale SDOs which cannot be solved by interior-point methods, due to excessive memory requirements \citep[see][]{ahmadi2015sum, bertsimas2019polyhedral}. }

\item \textbf{Copy of variables:} In problems with multiple complicating constraints, many feasibility cuts may be generated, which can hinder convergence greatly. If this occurs, we recommend introducing a copy of $\bm{x}$ in {\color{black}Problem \eqref{bilevelmainproblem}}, and imposing the following constraints:
\begin{align}
    \bm{l} \leq \bm{A}\bm{x} \leq \bm{u},\ \bm{e}^\top \bm{x}=1,\ \bm{x}\geq \bm{0},\ \bm{x} \leq \bm{z}
\end{align}
{\color{black} while keeping {\color{black}Problem \eqref{bilevelsubproblem}} the same}. This approach performs well on the highly constrained problems studied in Section \ref{sec:buyin}.
\end{itemize}

\subsection{Modeling Minimum Investment Constraints}
A frequently-studied extension to Problem {\color{black}\eqref{mainproblemnonreg}} is to impose minimum investment constraints, which control transaction fees by requiring that $x_i \in \{0\}\cup [x_{i,\min}, u_i]$. We now extend our saddle-point reformulation to cope with them.

By letting $z_i$ be a binary indicator variable which denotes whether we hold a non-zero position in the $i$th asset, we model these constraints via $z_i x_i \geq z_i x_{i,\min} \ \forall i \in [n].$ {\color{black}Moreover, we incorporate the upper bounds $u_i$ within our algorithmic framework by ``disappearing'' the constraints $x_i \leq u_i$ into the general linear constraint set $\bm{l} \leq \bm{A}\bm{x} \leq \bm{u}$.}

Moreover, by letting $\rho_i$ be the dual multiplier associated with the $i$th minimum investment constraint, and repeating the steps of our saddle-point reformulation, we retain efficient objective function and subgradient evaluations in the presence of these constraints. Specifically, including the constraints is equivalent to adding the term ${\sum_{i=1}^n \rho_i\left(z_i x_{i,\min}-z_i x_i\right)}$ to Problem \eqref{mainproblem}'s Lagrangian, which implies the saddle-point problem becomes:
\begin{equation}\label{saddlepointreformulationwithbuyin}
\begin{aligned}
    \min_{\bm{z} \in \mathcal{Z}_k^n} \ \max_{\substack{\bm{\alpha} \in \mathbb{R}^r, \ \bm{w} \in \mathbb{R}^n, \bm{\rho} \in \mathbb{R}_n^+ \\ \bm{\beta}_l,\ \bm{\beta}_u \in \mathbb{R}^m_+,\ \lambda \in \mathbb{R}}} \quad & -\frac{1}{2} \bm{\alpha}^\top \bm{\alpha} - \frac{\gamma}{2} \sum_i z_i w_i^2 +\bm{y}^\top \bm{\alpha} +\bm{\beta}_l^\top \bm{l}-\bm{\beta}_u^\top \bm{u}+\lambda+\sum_i \rho_i z_i x_{i,\min}\\
    \text{s.t.} \quad & \bm{w}\geq \bm{X}^\top \bm{\alpha}+\bm{A}^\top\left(\bm{\beta}_l-\bm{\beta}_u\right)+\lambda \bm{e}+\bm{\rho}-\bm{d}.
\end{aligned}
\end{equation}

Moreover, the subgradient with respect to each index $i$ becomes
\begin{align}\label{subgradexpressionwithbuyin}
   g_{\bm{z},i}=-\frac{\gamma}{2}w_i^{\star}(\bm{z})^2+\rho_i x_{i,\min}.
\end{align}

Finally, if $z_i=0$ then we can certainly set $\rho_i=0$ without loss of optimality. Therefore, we recommend solving a subproblem in the $k$ variables {\color{black}for which} $z_i>0$ and subsequently setting $\rho_i=0$ for the remaining variables, in the manner discussed in the previous subsection. Indeed, setting $w_i=\max\left(\bm{X}_i^\top \bm{\alpha}^\star +\bm{A}_i^\top\left(\bm{\beta}_l^\star-\bm{\beta}_u^\star\right) +\lambda^\star+\rho_i^\star -d_i,0\right)$ for each index $i$ where $z_i=0$, as discussed in the previous subsection, supplies the minimum absolute value of $w_i$.

{\color{black}
\subsection{Sensitivity Analysis}\label{ssec:sensitivityanalysis}
In this section, we study Problem \eqref{bilevelmainproblem}'s dependence on the regularization parameter $\gamma$. This is an important issue in practice, because, as suggested in the introduction, if we are interested in solving Problem {\color{black}\eqref{mainproblemnonreg}}, we can solve the regularized problem  to obtain a support vector $\bm{z}$, and subsequently resolve the unregularized problem with the support fixed to $\bm{z}$. Therefore, we are interested in the suboptimality {(\color{black}w.r.t. \eqref{mainproblemnonreg}) of an optimal solution} $\bm{z}^\star$ to Problem \eqref{mainproblem}.

We remark that the results in this section rely on basic sensitivity analysis proof techniques which can be found in several optimization textbook{\color{black}s} \citep[e.g.,][]{bertsimas1997introduction, rockafellar2009variational}. Nonetheless, we have included them, due to the central importance of regularization in this work, and because these results are not widely known.

Our first result demonstrates that the optimal support $\bm{z}$ for a larger value of $\gamma$ can serve as a high-quality warm-start for a problem with less regularization.
\begin{proposition}\label{prop:sensanal1}
Suppose that $k \in [n]$ is a fixed cardinality budget. Let $\bm{z}^\star(\gamma)$ denote an optimal solution to Problem \eqref{bilevelmainproblem} for a fixed regularizer $\gamma$, $f_\gamma(\bm{z})$ denote the optimal objective of Problem \eqref{bilevelsubproblem} for a fixed $\gamma$. Then, for any $\Delta >0$:
\begin{align}
    0 \leq f_{\gamma+\Delta}(\bm{z}^\star(\gamma))-f_{\gamma+\Delta}(\bm{z}^\star(\gamma+\Delta)) \leq \frac{1}{2\gamma}-\frac{1}{2(\gamma+\Delta)}.
\end{align}
\end{proposition}
\proof{Proof of Proposition \ref{prop:sensanal1}}
We have that:
\begin{align*}
     f_{\gamma+\Delta}(\bm{z}^\star(\gamma))-f_{\gamma+\Delta}(\bm{z}^\star(\gamma+\Delta))\leq & f_{\gamma}(\bm{z}^\star(\gamma))-f_{\gamma+\Delta}(\bm{z}^\star(\gamma+\Delta)) \\
    \leq & \left(\frac{1}{2\gamma}-\frac{1}{2(\gamma+\Delta)}\right) \Vert \bm{x}^\star(\bm{z}^\star(\gamma+\Delta))\Vert_2^2\leq \left(\frac{1}{2\gamma}-\frac{1}{2(\gamma+\Delta)}\right),
    \end{align*}
    where the first inequality holds because decreasing the amount of regularization can only lower the optimal objective value, the second inequality holds because  $\bm{x}^\star(\bm{z}^\star(\gamma+\Delta))$, an optimal choice of $\bm{x}$ with support indices $\bm{z}^\star(\gamma+\Delta)$, is a feasible solution with regularization parameter $\gamma$, specifically and the last inequality holds because all solutions $\bm{x}$ lie on the unit simplex.
    \hfill \Halmos
\endproof

Observe that, by setting $\Delta \rightarrow \infty$, Proposition \ref{prop:sensanal1} supplies a formal proof of Section $1.1$'s claim that $\bm{z}^\star(\gamma)$ is a $1/(2\gamma)$-optimal solution for $\gamma \rightarrow +\infty$.

Our next result justifies our claim in the introduction that for a sufficiently large $\gamma$ we recover the same optimal support from Problem \eqref{mainproblem} as the unregularized Problem \eqref{mainproblemnonreg}:

\begin{proposition}\label{prop:largestgamma}
Let $\bm{z}^\star(\gamma)$ denote an optimal solution to Problem \eqref{bilevelmainproblem} for a fixed regularizer $\gamma$, and $f_\gamma(\bm{z})$ denote the optimal objective of Problem \eqref{bilevelsubproblem} for a fixed $\gamma$. Then, there exists some parameter $\gamma_0 >0$ such that for any $\gamma \geq \gamma_0$, we have:
\begin{align}
    f_{\gamma}(\bm{z}^\star(\gamma_0))=f_{\gamma}(z^\star(\gamma)).
\end{align}
\end{proposition}

\proof{Proof of Proposition \ref{prop:largestgamma}}
Let us observe that, for each $\bm{z} \in \mathcal{Z}^k_n$, $f_{\gamma}(\bm{z})$ is concave in $\frac{1}{\gamma}$ as the pointwise minimum of functions which are linear in $\frac{1}{\gamma}$, and moreover $f_{\gamma}:=\min_{\bm{z} \in \mathcal{Z}^k_n} f_{\gamma}(\bm{z})$ is also a concave function in $\frac{1}{\gamma}$. By this concavity, it is a standard result from sensitivity analysis \citep[see, e.g.,][Chapter 5.6]{bertsimas1997introduction} that the set of all $\gamma$'s for which a particular $\bm{z}$ is optimal must form a (possibly open) interval. The result then follows directly from the finiteness of $\mathcal{Z}^k_n$.
\Halmos
\endproof

{\color{black}
Our final result in this section shows that the optimal support of the portfolio remains unchanged for sufficiently small $\gamma$'s (proof omitted, follows in the same fashion as Proposition \ref{prop:largestgamma}):

\begin{corollary}
Let $\bm{z}^\star(\gamma)$ denote an optimal solution to Problem \eqref{bilevelmainproblem} for a fixed regularizer $\gamma$, and $f_\gamma(\bm{z})$ denote the optimal objective of Problem \eqref{bilevelsubproblem} for a fixed $\gamma$. Then, there exists some parameter $\gamma_1 >0$ such that for any $\gamma \leq \gamma_1$, we have:
\begin{align}
    f_{\gamma}(\bm{z}^\star(\gamma_1))=f_{\gamma}(z^\star(\gamma)).
\end{align}
\end{corollary}
}

}
\subsection{Relationship with Perspective Cut Approach}\label{ssec:perspcut} 
In this section, we relate Algorithm \ref{alg:cuttingPlaneMethod} to the perspective cut approach introduced by \cite{frangioni2006perspective}. It turns out that taking the dual of the inner maximization problem in Problem \eqref{saddlepointreformulation} yields a perspective reformulation, and decomposing this reformulation in a way which retains the continuous and discrete variables in the master problem and outer-approximates the objective is precisely the perspective cut approach, as discussed in a more general setting by \citet[][Section 3.4]{bertsimas2019unified}. {\color{black} In this regard, our approach supplies a new and insightful derivation of the perspective cut scheme. In addition, our approach can easily be implemented within a modern mixed-integer optimization solver such as \verb|CPLEX| or \verb|Gurobi|, while existing implementations of the perspective cut approach often require tailored branch-and-bound schemes \citep[see][Section 3.1]{frangioni2006perspective}.

}

\section{Improving the Performance of the Cutting-Plane Method}\label{sec:improvedcuttingplanemethod}
In portfolio rebalancing applications, practitioners often require a high-quality solution to Problem \eqref{mainproblem} within a fixed time budget. Unfortunately, Algorithm \ref{alg:cuttingPlaneMethod} is ill-suited to this task: while it always identifies a certifiably optimal solution, it does not always do so within a time budget. In this section, we propose alternative techniques which sacrifice some optimality for speed, and discuss how they can be applied to improve the performance of Algorithm \ref{alg:cuttingPlaneMethod}. In Section \ref{ssec:localsearch1} we propose a warm-start heuristic which supplies a high-quality solution to Problem \eqref{mainproblem} a priori, and in Section \ref{ssec:socp} we derive a second order cone representable lower bound which is often very tight in practice. Taken together, these techniques supply a certifiably near optimal solution very quickly, which can often be further improved by running Algorithm \ref{alg:cuttingPlaneMethod} for a short amount of time.

\subsection{Improving the Upper Bound: A Warm-Start Heuristic}\label{ssec:localsearch1}
In branch-and-cut methods, a frequently observed source of inefficiency is that solvers explore highly suboptimal regions of the search space in considerable depth. To discourage this behavior, optimizers frequently supply a high-quality feasible solution (i.e., a warm-start), which is installed as an incumbent by the solver. Warm-starts are beneficial for two reasons. First, they improve Algorithm \ref{alg:cuttingPlaneMethod}'s upper bound. Second, they allow Algorithm \ref{alg:cuttingPlaneMethod} to prune vectors of partial solutions which are provably worse than the warm-start, which in turn improves Algorithm \ref{alg:cuttingPlaneMethod}'s bound quality, by reducing the set of feasible binaries which can be selected at each subsequent iteration. Indeed, by pruning suboptimal solutions, warm-starts encourage branch-and-cut methods to focus on regions of the search space which contain near-optimal solutions.

We now describe a heuristic which supplies high-quality solutions for Problem \eqref{mainproblem}, {\color{black}inspired by a heuristic due to} \citet[Algorithm 1]{bertsimas2016best}. The heuristic works under the assumption that $f(\bm{z})$ is $L$-Lipschitz continuous in $\bm{z}$, with Lipschitz {\color{black}continuous gradient} $g_{\bm{z}}$ {\color{black}such that
\begin{align*}
    \Vert \bm{g}_{\bm{z}_1}-\bm{g}_{\bm{z}_2}\Vert_2 \leq L \Vert \bm{z}_1-\bm{z}_2\Vert_2 \quad \forall \bm{z}_1, \bm{z}_2 \in \mathrm{Conv}\left(\mathcal{Z}^n_k\right).
\end{align*}
Note that this assumption} is justified whenever the optimal dual variables are bounded; see Corollary \ref{corr:lcont}. Under this assumption, the heuristic approximately minimizes $f(\bm{z})$ by iteratively minimizing a quadratic approximation of ${\color{black}f(\bm{z})}$ at ${\color{black}\bm{z}}$, namely ${\color{black}f(\bm{z})\approx\Vert \bm{z}-\bm{x}^\star(\bm{z})-\frac{1}{L}g_{\bm{z}}\Vert_2^2}$.

This idea is algorithmized as follows: given a sparsity pattern $\bm{z}_{\text{old}} \in \mathcal{Z}_k^n$ and an optimal sparse portfolio for this given sparsity pattern $\bm{x}^\star(\bm{z}_{\text{old}})$, the method iteratively solve{\color{black}s} the following problem, which ranks the differences between each security's contribution to the portfolio, $x_i^\star(\bm{z}_{\text{old}})$, and {\color{black}the associated entry $g_{\bm{z}_{i, \text{old}}}$ of the subgradient of $f$ at $\bm{z}_{\text{old}}$}:
\begin{equation}
\begin{aligned}
  \bm{z}_{\text{new}}:=\arg\min_{\bm{z} \in \mathcal{Z}_k^n}\left\Vert\bm{z}-\bm{x}^\star(\bm{z}_{\text{old}})+\frac{1}{L}\bm{g}_{\bm{z}_{\text{old}}}\right\Vert_2^2.
\end{aligned}
\end{equation}
Given $\bm{z}_{\text{old}}$, $\bm{z}_{\text{new}}$ can be obtained by setting $z_i=1$ for $k$ of the indices where $\left\vert -x_i^\star(\bm{z}_{\text{old}})+\frac{1}{L}g_{\bm{z}_{\text{old}},i}\right\vert $ is largest \citep[cf.][Proposition 3]{bertsimas2016best}. We formalize this warm-start procedure in Algorithm \ref{alg:warmstart}.
{\color{black}Note that Algorithm \ref{alg:warmstart} is a discrete first-order heuristic in the sense of \cite{nesterov2013gradient, nesterov2018lectures}, since as discussed in \citep[Section 3]{bertsimas2016best} it iteratively minimizes a first order approximation of $f(\bm{z})$ by solving for a stationary point in the first-order approximation.}

\begin{algorithm*}[h]
\caption{A discrete {\color{black}first-order} heuristic}
\label{alg:warmstart}
\begin{algorithmic}
\STATE $t \gets 1$
\STATE $\bm{z}_1 \gets $ randomly generated $k$-sparse binary vector.
{\color{black} \STATE $\bm{w}^\star \gets \bm{0}$.}
\WHILE{$\bm{z}_{t} \neq \bm{z}_{t-1}$ \textbf{and} $t < T$}
\STATE Set $\bm{w}_t$ optimal solution to:\vspace{-5mm}
\begin{align*}
     \max_{\substack{\bm{\alpha} \in \mathbb{R}^r, \ \bm{w} \in \mathbb{R}^n, \\ \bm{\beta}_l,\ \bm{\beta}_u \in \mathbb{R}^m_+,\ \lambda \in \mathbb{R}}} \quad & -\frac{1}{2} \bm{\alpha}^\top \bm{\alpha} - \frac{\gamma}{2} \sum_i {\color{black}z_{t,i}} w_i^2 +\bm{y}^\top \bm{\alpha} +\bm{\beta}_l^\top \bm{l}-\bm{\beta}_u^\top \bm{u}+\lambda\\
    \text{s.t.} \quad & \bm{w}\geq \bm{X}^\top \bm{\alpha}+\bm{A}^\top\left(\bm{\beta}_l-\bm{\beta}_u\right)+\lambda \bm{e}-\bm{d}.
\end{align*}\vspace{-10mm}
\STATE Average multipliers via $\bm{w}^\star \gets \frac{1}{t}\bm{w}_t+\frac{t-1}{t}\bm{w}^\star$.
\STATE Set $g_{\bm{z},i}= \frac{-\gamma}{2}w_i^{\star 2} \quad \forall i \in [n]$, ${\color{black}x_{t,i}}= \gamma w_i^\star \quad \forall i \in [n]: {\color{black}z_{t,i}} =1$, $\bm{z}_{t+1}=$ $\arg\min\limits_{\bm{z} \in \mathcal{Z}_k^n} \left\Vert\bm{z}-\bm{x}_{t}+\frac{1}{L}\bm{g}_{\bm{z}_t} \right\Vert_2^2$
\STATE $t \leftarrow t+1 $
\ENDWHILE
\RETURN $\bm{z}_{t}$
\end{algorithmic}
\end{algorithm*}

Some remarks on Algorithm \ref{alg:warmstart} are now in order:
\begin{itemize}
    \item In our numerical experiments, we run Algorithm \ref{alg:warmstart} from five different randomly generated $k$-sparse binary vectors, to increase the probability that it identifies a high-quality solution.
    \item Averaging the dual multipliers across iterations, as suggested in the pseudocode, improves the method's performance; note that the contribution of each $\bm{w}_t$ to $\bm{w}^\star$ is $\frac{1}{t}\prod_{i=t+1}^{T_{final}} \frac{i-1}{i}=\frac{1}{{T_{final}}}$, where ${T_{final}}$ is the total number of iterations completed by Algorithm \ref{alg:warmstart} {\color{black} and we initially set $\bm{w}^\star=\bm{0}$ in order that $\bm{w}^\star$ is defined when we perform the averaging step at the first iteration; note that when $t=1$ we have $(t-1)/t=0$ so the initialization is unimportant}.
   {\color{black} \item Each $\bm{w}_t$ is the optimal solution of a convex quadratic optimization problem which can be reformulated as a (rotated) second-order cone program. Therefore, each $\bm{w}_t$ can be obtained via a standard second-order cone solver such as \verb|CPLEX|, \verb|Gurobi| or \verb|Mosek|.}
    \item {\color{black} The Lipschitz constant $L$ is motivated as an upper bound on an entry in a subgradient of $f(\bm{z})$, $\frac{\gamma}{2}w_i^2$. However, Algorithm \ref{alg:warmstart} is ultimately a heuristic method. Therefore, we recommend picking $L$ by cross-validating to minimize the objective obtained by Algorithm \ref{alg:warmstart}. {\color{black}In practice, setting $L=10$ was sufficient to reliably obtain high-quality solutions in Section \ref{sec:compexperiments}, because Algorithm \ref{alg:refinedCuttingPlaneMethod} (see Section \ref{sec:impcutmethod}) invokes a judicious combination of outer-approximation cuts and this warm-start to convert this warm-start into an optimal solution within seconds. Therefore, we set $L=10$ throughout Section \ref{sec:compexperiments}, although it may be appropriate to cross-validate $L$ if running the method on new data.}
    }
\end{itemize}

\subsection{Improving the Lower Bound: A Second Order Cone Relaxation}
\label{ssec:socp}
In financial applications, we sometimes require a certifiably near-optimal solution quickly but do not have time to certify optimality. Therefore, we now derive near-exact lower bounds which can be computed in polynomial time. Immediately, we see that we obtain a valid lower bound by relaxing the constraint $\bm{z} \in \mathcal{Z}_k^n$ to $\bm{z} \in \mathrm{Conv}(\mathcal{Z}_k^n)$ in Problem \eqref{mainproblem}. By invoking strong duality, we now demonstrate that this lower bound can be obtained by solving a single second order cone problem.
\begin{theorem}\label{socptheorem}
Suppose that Problem \eqref{mainproblem} is feasible. Then, the following three optimization problems attain the same optimal value:
\begin{equation}\label{prob:minmaxboolrelax}
    \begin{aligned}
            \min_{\bm{z} \in \text{\rm Conv}\left(\mathcal{Z}_k^n\right)} \ \max_{\substack{\bm{\alpha} \in \mathbb{R}^r, \ \bm{w} \in \mathbb{R}^n,\\ \bm{\beta}_l,\ \bm{\beta}_u \in \mathbb{R}^m_+,\ \lambda \in \mathbb{R}}} \quad & -\frac{1}{2} \bm{\alpha}^\top \bm{\alpha} - \frac{\gamma}{2} \sum_i z_i w_i^2 +\bm{y}^\top \bm{\alpha} +\bm{\beta}_l^\top \bm{l}-\bm{\beta}_u^\top \bm{u}+\lambda\\
    {\mbox{\rm s.t.}} \quad & \bm{w}\geq \bm{X}^\top \bm{\alpha}+\lambda \bm{e}+\bm{A}^\top(\bm{\beta}_l-\bm{\beta}_u)-\bm{d}.
    \end{aligned}
\end{equation}

\begin{equation}\label{socpbound}
    \begin{aligned}
    \max_{\substack{\bm{\alpha} \in \mathbb{R}^r, \ \bm{v} \in \mathbb{R}^n_+, \ \bm{w} \in \mathbb{R}^n, \ \\ \bm{\beta}_l,\ \bm{\beta}_u \in \mathbb{R}^m_+,\ \lambda \in \mathbb{R}, \ t \in \mathbb{R}_{+}}} &
            -\frac{1}{2} \bm{\alpha}^\top \bm{\alpha}+\bm{y}^\top \bm{\alpha} +\bm{\beta}_l^\top \bm{l}-\bm{\beta}_u^\top \bm{u}+\lambda
            -\bm{e}^\top \bm{v} -k t \\
    {\mbox{\rm s.t.}} \quad & \bm{w}\geq \bm{X}^\top \bm{\alpha}+\lambda \bm{e}+\bm{A}^\top(\bm{\beta}_l-\bm{\beta}_u)-\bm{d}, \ v_i \geq \frac{\gamma}{2} w_i^2-t \quad \forall i \in [n].
    \end{aligned}
\end{equation}

\begin{equation}
\begin{aligned}
\label{mainproblemmisocp_dualversion}
    \min_{\bm{z} \in \mathrm{Conv}\left(\mathcal{Z}_k^n\right)} \ \min_{\bm{x} \in \mathbb{R}_+^n, \bm{\theta} \in \mathbb{R}_+^n}\quad & \frac{1}{2}\left\Vert \bm{X}\bm{x}-\bm{y}\right\Vert_2^2+\frac{1}{2\gamma}\bm{e}^\top \bm{\theta}+\bm{d}^\top \bm{x} \\
    {\mbox{\rm s.t.}} \quad &  \bm{l} \leq \bm{A}\bm{x} \leq \bm{u},\ \bm{e}^\top \bm{x}=1, \ x_i^2 \leq z_i \theta_i \quad \forall i \in [n].
\end{aligned}
\end{equation}

\end{theorem}
\begin{remark} We recognize Problem \eqref{mainproblemmisocp_dualversion} as a perspective relaxation of Problem \eqref{mainproblem} \citep[see][for a survey]{gunluk2012perspective}. As perspective relaxations are often near-exact in practice \citep{frangioni2006perspective, frangioni2009computational} this explains why the second-order cone bound is high-quality.
\end{remark}

{\color{black}
\proof{Proof of Theorem \ref{socptheorem}}

Problem \eqref{prob:minmaxboolrelax} is strictly feasible, since the interior of $\mathrm{Conv}(\mathcal{Z}_k^n)$ is non-empty and $\bm{w}$ can be increased without bound. Therefore, the Sion-Kakutani minimax theorem \citep[Appendix D.4.]{ben2001lectures} holds, and we can exchange the minimum and maximum operators in Problem \eqref{prob:minmaxboolrelax}, to yield:
\begin{equation}
    \begin{aligned}
    \max_{\substack{\bm{\alpha} \in \mathbb{R}^r,\ \bm{w} \in \mathbb{R}^n,\\ \bm{\beta}_l,\ \bm{\beta}_u \in \mathbb{R}^m_+,\ \lambda \in \mathbb{R}}} &
            -\frac{1}{2} \bm{\alpha}^\top \bm{\alpha}+\bm{y}^\top \bm{\alpha} +\bm{\beta}_l^\top \bm{l}-\bm{\beta}_u^\top \bm{u}+\lambda
            - \frac{\gamma}{2} \max_{\bm{z} \in \text{\rm Conv}\left(\mathcal{Z}_k^n\right)} \sum_i z_i w_i^2 \\
    \text{s.t.} \quad & \bm{w}\geq \bm{X}^\top \bm{\alpha}+\lambda \bm{e}+\bm{A}^\top(\bm{\beta}_l-\bm{\beta}_u)-\bm{d}.\\
    \end{aligned}
\end{equation}
Next, fixing $\bm{w}$ and applying strong duality between the inner primal problem
\begin{align*}
    \max_{\bm{z} \in \mathrm{Conv}\left(\mathcal{Z}_k^n\right)}\sum_i \frac{\gamma}{2} z_i w_i^2 \ = \max_{\bm{z}} \ \sum_i \frac{\gamma}{2} z_i w_i^2 \quad \text{s.t.} \quad \bm{0} \leq \bm{z} \leq \bm{e}, \ \bm{e}^\top \bm{z}\leq k,
\end{align*}
and its dual problem $$\min_{\bm{v} \in \mathbb{R}^n_+, t \in \mathbb{R}_+} \ \bm{e}^\top \bm{v} +k t \quad \text{s.t.} \quad  v_i+t \geq \frac{\gamma}{2} w_i^2 \quad \forall i \in [n]$$ proves that strong duality holds between Problems \eqref{prob:minmaxboolrelax}-\eqref{socpbound}.

Next, we observe that Problems \eqref{socpbound}-\eqref{mainproblemmisocp_dualversion} are dual, as can be seen by applying the relation $$bc \geq a^2, b, c \geq 0 \iff \left\Vert \begin{pmatrix}
                 2a\\
                 b-c
                 \end{pmatrix} \right\Vert \leq b+c$$ to rewrite Problem \eqref{socpbound} as an SOCO in standard form, and applying SOCO duality \citep[see, e.g.,][Exercise 5.43]{boyd2004convex}. Moreover, since Problem \eqref{socpbound} is strictly feasible (as $\bm{v}$, $\bm{w}$ are unbounded from above) strong duality must hold between these problems.\hfill \Halmos
\endproof
}

{\color{black} Having derived Problem \eqref{mainproblem}'s bidual, namely Problem \eqref{mainproblemmisocp_dualversion}, it follows from a direct application of convex analysis that the duality gap between Problem \eqref{mainproblem} and \eqref{mainproblemmisocp_dualversion}, $\Delta_\gamma$, decreases as we decrease $\gamma$ and becomes $0$ at some finite $\gamma>0$. Note however that this $\Delta_\gamma$ will, in general, depend upon the problem data \citep[see][Theorem XII.5.2.2]{hiriart2013convex}. 
This observation justifies our claim in the introduction that decreasing $\gamma$ makes Problem \eqref{mainproblem} easier.
}

We now derive conditions under which Problem \eqref{socpbound} provides an optimal solution to Problem \eqref{mainproblem} a priori (proof deferred to Appendix \ref{sec:ommitedproofs}). 

\begin{corollary}{\rm \textbf{A sufficient condition for support recovery}}\\
\label{corr:exactrecovery}
Let there exist some $\bm{z} \in \mathcal{Z}_k^n$ and set of dual multipliers $(\bm{v}^\star, \bm{w}^\star, \bm{\alpha}^\star, \bm{\beta}_l^\star, \bm{\beta}_u^\star, \lambda^\star)$ which solve Problem \eqref{socpbound}, such that these two quantities collectively satisfy the following conditions:

\begin{align}
\begin{aligned}
\label{conditions}
  \gamma \sum_i  z_i w_i^\star =1,\ \bm{l} \leq \gamma \sum_i  \bm{A}_i w^\star_i z_i \leq \bm{u},\ z_i w_i \geq 0 \quad \forall i \in [n],\ v_i^\star=0 \quad \forall i \in [n]\ : \ z_i=0.
\end{aligned}
\end{align}

Then, Problem \eqref{socpbound}'s lower bound is exact. Moreover, let $\vert w^\star\vert_{[k]}$ denote the $k$th largest entry in $\bm{w}^\star$ {\rm by absolute magnitude}. If $\vert w^\star\vert_{[k]} > \vert w^\star\vert_{[k+1]}$ in Problem \eqref{socpbound} then setting $$z_i=1 \quad \forall i: \vert w^\star_i\vert \geq \vert w^\star\vert_{[k]},\ z_i=0 \quad \forall i: \vert w^\star_i\vert < \vert w^\star\vert _{[k]}$$ supplies a $\bm{z} \in \mathcal{Z}_k^n$ which satisfies the above condition and hence solves Problem \eqref{mainproblem}.
\end{corollary}


{\color{black}
We now apply Theorem \ref{socptheorem} to prove that if $\bm{\Sigma}$ is a diagonal matrix, $\bm{\mu}$ is a multiple of the vector of all $1$'s and the matrix $\bm{A}$ is empty then Problem \eqref{mainproblem} is solvable in closed-form. Let us first observe that under these conditions Problem \eqref{mainproblem} is equivalent to
\begin{align*}
    \min \quad & \sum_i \frac{1}{2\gamma_i}x_i^2 \quad \text{s.t.} \quad \bm{e}^\top \bm{x}=1, \bm{x} \geq \bm{0}, \Vert \bm{x}\Vert_0 \leq k.
\end{align*}

We now have the following result:

\begin{corollary}\label{corr:polytimediag2}
Let $0 < \gamma_n \leq \gamma_{n-1} \leq ... \gamma_1$. Then, strong duality holds between the problem
\begin{align}\label{diagproblem2}
    \min \quad & \sum_i \frac{1}{2\gamma_i}x_i^2 \quad \textrm{\rm s.t.} \quad \bm{e}^\top \bm{x}=1, \bm{x} \geq \bm{0}, \Vert \bm{x}\Vert_0 \leq k
\end{align}
and its second-order cone relaxation:
\begin{equation}\label{socpboundiidstocks2}
    \begin{aligned}
    \max_{\substack{\bm{v} \in \mathbb{R}_n^+, \ \bm{w} \in \mathbb{R}^n, \ \\ \lambda \in \mathbb{R}, \ t \in \mathbb{R}_+}} &\lambda
            -\bm{e}^\top \bm{v} -k t \\
    \textrm{\rm s.t.} \quad & \bm{w}\geq \lambda \bm{e},\quad v_i \geq \frac{\gamma_i}{2} w_i^2-t \quad \forall i \in [n].
    \end{aligned}
\end{equation}

Moreover, an optimal solution to Problem \eqref{diagproblem2} is $x_i=\frac{\gamma_i}{\sum_{i=1}^k \gamma_i}$ for $i \leq k$, $x_i=0$ for $i>k$.
\end{corollary}

\proof{Proof of Corollary \ref{corr:polytimediag2}}{}
By Theorem \ref{socptheorem}, a valid lower bound to Problem \eqref{diagproblem2} is given by the SOCO:
\begin{equation}\label{socpboundiidstocks3}
    \begin{aligned}
    \max_{\substack{\bm{v} \in \mathbb{R}_n^+, \ \bm{w} \in \mathbb{R}^n, \ \\ \lambda \in \mathbb{R}, \ t \in \mathbb{R}_+}} &\lambda
            -\bm{e}^\top \bm{v} -k t \\
    \text{s.t.} \quad & \bm{w}\geq \lambda \bm{e},\quad  v_i \geq \frac{\gamma_i}{2} w_i^2-t \quad \forall i \in [n].
    \end{aligned}
\end{equation}

Let us assume that $\lambda^\star \geq 0$ (otherwise the objective value cannot exceed $0$, which is certainly suboptimal). Then, we can let the constraint $w_i \geq \lambda$ be binding without loss of optimality for each index $i$, i.e., set $\bm{w}=\lambda \bm{e}$ for some $\lambda$. This allows us to simplify this problem to:
\begin{equation}\label{socpboundiidstocks4}
    \begin{aligned}
    \max_{\substack{\bm{v} \in \mathbb{R}_n^+, \ \lambda \in \mathbb{R}, \ t \in \mathbb{R}+}} &\lambda
            -\bm{e}^\top \bm{v} -k t \\
    \text{s.t.} \quad & v_i \geq \frac{\gamma_i}{2} \lambda^2-t \quad \forall i \in [n].
    \end{aligned}
\end{equation}

The KKT conditions for max-$k$ norms \citep[see, e.g., ][Lemma 1]{zakeri2014optimization} then reveal that an optimal choice of $t$ is given by the $k$th largest value of $\frac{\gamma_i}{2}\lambda^2$, i.e., $t^\star=\frac{\gamma_k}{2}\lambda^2$ and an optimal choice of $v_i$ is given by $v_i=\max\left(\frac{\gamma_i}{2}\lambda^2-t, 0\right)$, i.e., $$v_i^\star=\begin{cases}
\frac{\gamma_i-\gamma_k}{2}\lambda^2 \quad \forall i \leq k,\\ v_i=0 \quad \quad \forall i >k.\end{cases}$$

Substituting these terms into the objective function gives an objective of $$\lambda-\sum_{i=1}^k \frac{\gamma_i}{2}\lambda^2,$$ which implies that an optimal choice of $\lambda$ is $\lambda=\frac{1}{\sum_{i=1}^k \gamma_i}$. Next, substituting the expression $\lambda=\frac{1}{\sum_{i=1}^k \gamma_i}$ into the objective function gives an objective value of $\frac{\lambda}{2}$, which implies that a lower bound on Problem \eqref{diagproblem2}'s objective is $\frac{1}{2\sum_{i=1}^k \gamma_i}$.

Finally, we construct a primal solution via $z_i=1 \forall i \leq k$, and the primal-dual KKT condition $x_i = \gamma_i z_i w_i=\gamma_i z_i \lambda=\frac{\gamma_i z_i}{\sum_{i=1}^k \gamma_i}$. This is feasible, by inspection. Moreover, it has an objective value of $$\sum_{i=1}^k \frac{1}{2\gamma_i}\left(\gamma_i \lambda\right)^2=\frac{\lambda}{2}\sum_{i=1}^k \gamma_i \lambda=\frac{\lambda}{2},$$ and therefore is optimal.

\hfill \Halmos
\endproof
}



\subsection{An Improved Cutting-Plane Method}\label{sec:impcutmethod}
We close this section by combining Algorithm \ref{alg:cuttingPlaneMethod} with the improvements discussed in this section, to obtain an efficient numerical approach to Problem \eqref{mainproblem}, which we present in Algorithm \ref{alg:refinedCuttingPlaneMethod}. Note that we use the larger of $\theta_t$ and the second-order cone lower bound in our termination criterion, as the second-order cone gap is sometimes less than $\epsilon$. 

\begin{algorithm*}
\caption{A refined cutting-plane method for Problem \eqref{mainproblem}.}
\label{alg:refinedCuttingPlaneMethod}
\begin{algorithmic}
\REQUIRE Initial warm-start solution $\bm{z}_1$
\STATE $t \leftarrow 1 $
\STATE Set $\theta_{\text{SOCO}}$ optimal objective value of Problem \eqref{socpbound}
\REPEAT
\STATE Compute $\bm{z}_{t+1}, \theta_{t+1}$ solution of \vspace{-5mm}
\begin{align*}
\min_{\bm{z} \in \mathcal{Z}_k^n, \theta} \: \theta \quad \mbox{ s.t. } \quad \theta \geq f(\bm{z}_i) + g_{\bm{z}_i}^\top (\bm{z}-\bm{z}_i) \quad \forall i \in [t].
\end{align*}\vspace{-10mm}
\STATE Compute $f(\bm{z}_{t+1})$ and $g_{\bm{z}_{t+1}} \in \partial f (\bm{z}_{t+1})$
\STATE $t \leftarrow t+1 $
\UNTIL{$ f(\bm{z}_t)-\max(\theta_t, \theta_{\text{SOCO}}) \leq \varepsilon$}
\RETURN $\bm{z}_t$
\end{algorithmic}
\end{algorithm*}

\begin{figure}[h]
    \centering
    \includegraphics[scale=0.5]{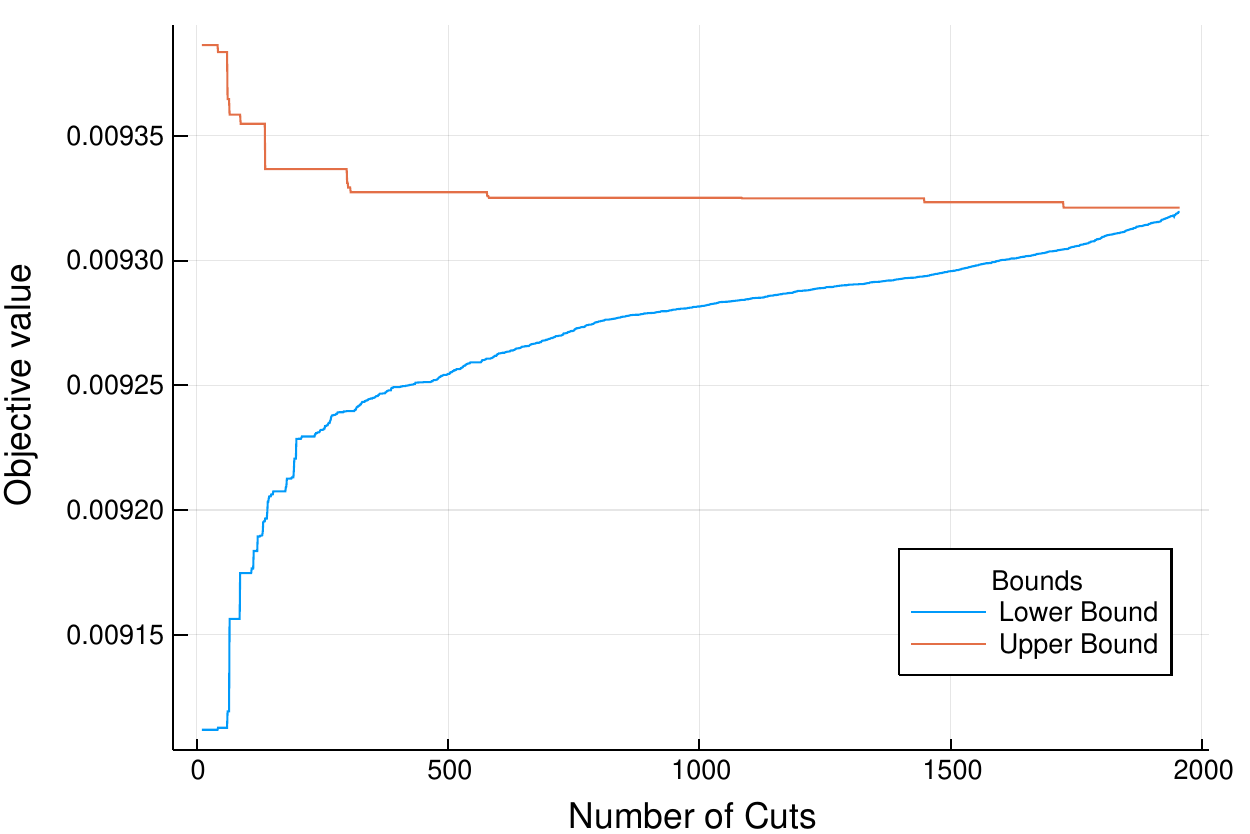}
    \caption{Convergence of Algorithm \ref{alg:refinedCuttingPlaneMethod} on the OR-library problem \textit{port2} with a minimum return constraint and a cardinality constraint $\Vert \bm{x}\Vert_0 \leq 5$. The behavior shown here is typical.}
    \label{fig:convergenceplot}
\end{figure}

Figure \ref{fig:convergenceplot} depicts the algorithm's convergence on the problem \textit{port2} with a cardinality value $k=5$ and a minimum return constraint, as described in Section \ref{ssec:compwithbigm}. Note that we did not use the second-order cone lower bound when generating this plot; the second-order cone lower bound is $0.009288$ in this instance, and Algorithm \ref{alg:refinedCuttingPlaneMethod} requires $1 225$ cuts to improve upon this bound.

\section{Computational Experiments on Real-World Data}\label{sec:compexperiments}
In this section, we evaluate our outer-approximation method, implemented in \verb|Julia| $1.1$ using the \verb|JuMP.jl| package version $0.18.5$ and solved using \verb|CPLEX| version $12.8.0$ for the master problems, and \verb|Mosek| version $9.0$ for the continuous quadratic subproblems. We compare the method against big-$M$ and MISOCO formulations of Problem \eqref{mainproblem}, solved in \verb|CPLEX|. {To bridge the gap between theory and practice, we have made our code freely available on \verb|Github| at \href{https://github.com/ryancorywright/SparsePortfolioSelection.jl}{github.com/ryancorywright/SparsePortfolioSelection.jl}. }

All experiments were performed on a MacBook Pro with a $2.9$GHz i$9$ Intel\textregistered  CPU and $16$GB DDR$4$ Memory. For simplicity, we ran all methods on one thread, using default \verb|CPLEX| parameters.

In all experiments {\color{black}which do not involve minimum investment constraints}, we solve the following optimization problem, which places a multiplier $\kappa$ on the return term but{\color{black}—since the multiplier can be absorbed into the return vector $\bm{\mu}$—}is mathematically equivalent to Problem \eqref{mainproblem}:
\begin{align}
\label{mainproblem2}
    \min_{\bm{x} \in \mathbb{R}^n_+} \ \frac{1}{2} \bm{x}^\top \bm{\Sigma} \bm{x}+\frac{1}{2\gamma}\Vert \bm{x} \Vert_2^2-\kappa\bm{\mu}^\top \bm{x} \quad \text{s.t.} \quad  \bm{l} \leq \bm{A}\bm{x} \leq \bm{u},\ \bm{e}^\top \bm{x}=1, \ \Vert \bm{x} \Vert_0 \leq k.
\end{align}
Note that we only consider cases where $\kappa=0$ or $\kappa=1$, depending on whether we are penalizing low expected return portfolios in the objective or constraining the portfolios expected return.

{\color{black}In the presence of minimum investment constraints, we augment Problem \eqref{mainproblem2} with the constraints $x_i \geq z_i x_{i, \min}$, which can be straightforwardly managed by the algorithms derived in this paper in the manner discussed in Section \ref{sec:buyin}.}

We aim to answer the following questions:
\begin{enumerate}
    \item How does Algorithm \ref{alg:refinedCuttingPlaneMethod} compare to existing {\color{black}commercial codes such as \verb|CPLEX|?}
    \item How do constraints affect Algorithm \ref{alg:refinedCuttingPlaneMethod}'s scalability?
    \item How does Algorithm \ref{alg:refinedCuttingPlaneMethod} scale as a function of the number of securities in the buyable universe?
    \item How sensitive are optimal solutions to Problem \eqref{mainproblem} to the hyperparameters $\kappa, \gamma, k$?
\end{enumerate}

\subsection{A Comparison Between Algorithm \ref{alg:refinedCuttingPlaneMethod} and Existing Commercial Codes}\label{ssec:compwithbigm}
We now present a direct comparison of Algorithm \ref{alg:refinedCuttingPlaneMethod} with \verb|CPLEX| version $12.8.0$, where \verb|CPLEX| uses both big-M and MISOCO formulations of Problem \eqref{mainproblem}. Note that the MISOCO formulation which we pass directly to \verb|CPLEX| is \citep[cf.][]{ben2001lectures, akturk2009strong}:
\begin{align}
\label{mainproblemmisocp}
    \min_{\bm{z} \in \mathcal{Z}_k^n, \bm{x} \in \mathbb{R}_+^n, \bm{\theta} \in \mathbb{R}_+^n}\ \frac{1}{2} \bm{x}^\top \bm{\Sigma} \bm{x}+\frac{1}{2\gamma}\bm{e}^\top \bm{\theta}-{\color{black}\kappa}\bm{\mu}^\top \bm{x} \quad \text{s.t.} \quad  \bm{l} \leq \bm{A}\bm{x} \leq \bm{u},\ \bm{e}^\top \bm{x}=1, \ x_i^2 \leq z_i \theta_i \quad \forall i \in [n].
\end{align}
{\color{black} We remark that \verb|CPLEX|'s MISOCO routine can be interpreted as an aggregation of the techniques reviewed in the introduction, since \verb|CPLEX| mines the optimization literature for techniques which improve its performance. For instance, as discussed by \citet{vielma2017advances}, \verb|CPLEX|'s MISOCO routine currently includes the conic disaggregation technique of \cite{vielma2017extended}. Indeed, this aggregation has been so successful that \cite{bienstock2010eigenvalue} showed that, in $2010$, \verb|CPLEX| outperform{\color{black}ed} several of the methods reviewed in the introduction.
}

We compare the three approaches in two distinct situations. First, when no constraints are applied and the system $\bm{l} \leq \bm{A}\bm{x}\leq\bm{u}$ is empty, and second when a minimum return constraint is applied, i.e.,  $\bm{\mu}^\top \bm{x}\geq \bar{r}$. In the former case we set $\kappa=1$, while in the latter case we set $\kappa=0$ and as suggested by \citet{cesarone2009efficient, zheng2014improving} we set $\bar{r}$ in the following manner: Let
\begin{align*}
    r_{\min}=&\bm{\mu}^\top \bm{x}_{\min} & \ \text{where} & \ \bm{x}_{\min}=\arg\min_{\bm{x}} \quad \frac{1}{2}\bm{x}^\top \left(\frac{1}{\gamma}\mathbb{I}+\bm{\Sigma}\right)\bm{x} \quad & \text{s.t.} \quad \bm{e}^\top \bm{x}=1, \bm{x} \geq \bm{0},\\
r_{\max}=&\bm{\mu}^\top \bm{x}_{\max} & \ \text{where} & \ \bm{x}_{\max}=\arg\max_{\bm{x}} \quad \bm{\mu}^\top \bm{x}-\frac{1}{2\gamma}\bm{x}^\top \bm{x} \quad & \text{s.t.} \quad \bm{e}^\top \bm{x}=1, \bm{x} \geq \bm{0},
\end{align*}
and set $\bar{r}=r_{\min}+0.3(r_{\max}-r_{\min})$.

Table \ref{tab:comparison} (resp. Table \ref{tab:comparisonwithminreturn}) depicts the time required for all $3$ approaches to determine an optimal allocation of funds without (resp. with) the minimum return constraint. The problem data is taken from the $5$ mean-variance portfolio optimization problems described by \cite{chang2000heuristics} and included in the OR-library test set by \citet{beasley1990or}. Note that we turned off the second-order cone lower bound for these tests, and ensured feasibility in the master problem by imposing $\sum_{i \in [n]: \mu_i \geq \bar{r}} z_i \geq 1$ when running Algorithm \ref{alg:refinedCuttingPlaneMethod} on the instances with a minimum return constraint (see Section \ref{sec:cuttingolanemethodcont}). Furthermore, as Algorithm \ref{alg:refinedCuttingPlaneMethod} is slow to converge for some instances with a return constraint, in Table \ref{tab:comparisonwithminreturn} we ran the method after applying $50$ cuts at the root node, generated using the \verb|in-out| method (see Appendix \ref{sec:auxpsueodocode}, for the relevant pseudocode).

\begin{table}[h]
\centering\footnotesize
\caption{Runtime in seconds per approach with $\kappa=1$, $\gamma=\frac{100}{\sqrt{n}}$ and no constraints in the system $\bm{l} \leq \bm{A} \bm{x} \leq \bm{u}$. We impose a time limit of $300$s and run all approaches on one thread. If a solver fails to converge, we report the number of explored nodes at the time limit.}
\begin{tabular}{@{}l l l r r r r r r r r@{}} \toprule \footnotesize
Problem &  $n$ & $k$ & \multicolumn{3}{c@{\hspace{0mm}}}{Algorithm \ref{alg:refinedCuttingPlaneMethod}} &  \multicolumn{2}{c@{\hspace{0mm}}}{CPLEX Big-M} &  \multicolumn{2}{c@{\hspace{0mm}}}{CPLEX MISOCO} \\
\cmidrule(l){4-6} \cmidrule(l){7-8} \cmidrule(l){9-10} & &  & Time & Nodes & Cuts & Time & Nodes & Time & Nodes \\\midrule
port $1$ & $31$ & $5$   & $0.17$   & $0$     & $4$      & $1.98$    & $31,640$    & $\textbf{0.03}$    & $0$ \\
         &      & $10$  & $0.16$   & $0$     & $4$      & $1.11$    & $16,890$    & $\textbf{0.01}$    & $0$ \\
         &      & $20$  & $0.14$   & $0$     & $4$      & $\textbf{0.01}$    & $108$       & $0.03$    & $0$ \\\midrule
port $2$ & $85$ & $5$   & $\textbf{0.01}$   & $0$     & $4$      & $>300$    & $1,968,000$ & $0.11$    & $0$ \\
         &      & $10$  & $\textbf{0.01}$   & $0$     & $4$      & $>300$    & $2,818,000$ & $0.12$    & $0$ \\
         &      & $20$  & $\textbf{0.01}$   & $0$     & $4$      & $>300$    & $3,152,000$ & $0.29$    & $0$ \\\midrule
port $3$ & $89$ & $5$   & $\textbf{0.01}$  & $0$      & $8$      & $>300$    & $2,113,000$ & $0.38$    & $0$ \\
         &      & $10$  & $\textbf{0.01}$  & $0$      & $4$      & $>300$    & $2,873,000$ & $0.41$    & $0$ \\
         &      & $20$  & $\textbf{0.02}$  & $0$      & $4$      & $>300$    & $2,998,000$ & $0.11$    & $0$ \\\midrule
port $4$ & $98$ & $5$   & $\textbf{0.03}$  & $0$      & $8$      & $>300$    & $1,888,000$ & $0.41$    & $0$ \\
         &      & $10$  & $\textbf{0.02}$  & $0$      & $8$      & $>300$    & $2,457,000$ & $2.74$    & $3$ \\
         &      & $20$  & $\textbf{0.03}$  & $0$      & $9$      & $>300$    & $2,454,000$ & $0.38$    & $0$ \\\midrule
port $5$ & $225$& $5$   & $\textbf{0.15}$  & $0$      & $9$      & $>300$    & $676,300$   & $11.17$   & $9$ \\
         &      & $10$  & $\textbf{0.02}$  & $0$      & $4$      & $>300$    & $926,600$   & $3.04$    & $0$ \\
         &      & $20$  & $\textbf{0.03}$  & $0$      & $7$      & $>300$    & $902,100$   & $2.88$    & $0$ \\
\bottomrule
\end{tabular}
\label{tab:comparison}
\end{table}

\begin{table}[h]
\centering\footnotesize
\caption{Runtime in seconds per approach with $\kappa=0$, $\gamma=\frac{100}{\sqrt{n}}$ and a minimum return constraint $\bm{\mu}^\top \bm{x} \geq \bar{r}$. We impose a time limit of $3600$s and run all approaches on one thread. If a solver fails to converge, we report the number of explored nodes at the time limit.}
\begin{tabular}{@{}l l l r r r r r r r r r r r@{}} \toprule
Problem &  $n$ & $k$ & \multicolumn{3}{c@{\hspace{0mm}}}{Algorithm \ref{alg:refinedCuttingPlaneMethod}} & \multicolumn{3}{c@{\hspace{0mm}}}{Algorithm \ref{alg:refinedCuttingPlaneMethod}+in-out} &  \multicolumn{2}{c@{\hspace{0mm}}}{CPLEX Big-M} &  \multicolumn{2}{c@{\hspace{0mm}}}{CPLEX MISOCO} \\
\cmidrule(l){4-6} \cmidrule(l){7-9} \cmidrule(l){10-11} \cmidrule(l){12-13} & &  & Time & Nodes & Cuts & Time & Nodes & Cuts & Time & Nodes & Time & Nodes \\\midrule
port $1$ & $31$ & $5$   & $\textbf{0.22}$  & $161$    & $32$     & $0.23$  & $113$    & $19$      & $9.32$    & $119,200$     & $0.83$    & $47$ \\
         &      & $10$  & $\textbf{0.20}$  & $159$    & $28$     & $0.25$  & $86$    & $25$      & $1970$    & $30,430,000$  & $0.84$    & $44$ \\
         &      & $20$  & $0.16$  & $0$      & $7$    & $0.16$  & $0$    & $4$     & $258.4$   & $4,966,000$   & $\textbf{0.05}$    & $0$ \\\midrule
port $2$ & $85$ & $5$   & $48.29$ & $73,850$ & $1,961$  & $\textbf{31.47}$  & $42,950$    & $1,272$      & $>3,600$   & $15,020,000$  & $91.98$   & $1,163$ \\
         &      & $10$  & $807.3$ & $243,500$& $6,433$  & $891.97$  & $255,200$    & $6,019$      & $>3,600$   & $20,890,000$  & $\textbf{82.44}$ & $902$ \\
         &      & $20$  & $\textbf{10.52}$ & $12,260$ & $1,224$  & $13.0$  & $13,650$    & $1,350$    & $>3,600$   & $21,060,000$  & $24.54$   & $210$ \\\midrule
port $3$ & $89$ & $5$   & $175.2$ & $132,700$& $3,187$  & $\textbf{151.1}$  & $96,010$    & $2,345$     & $>3,600$   & $14,680,000$  & $213.3$   & $2,528$ \\
         &      & $10$  & $>3,600$& $439,400$& $9,851$  & $>3,600$  & $490,400$    & $11,310$     & $>3,600$   & $20,710,000$  & $\textbf{531.3}$   & $5,776$ \\
         &      & $20$  & $119.5$ & $65,180$ & $4,473$  & $60.03$  & $40,240$   & $3,275$     & $>3,600$   & $22,240,000$  & $\textbf{21.32}$   & $170$ \\\midrule
port $4$ & $98$ & $5$   & $2,690$ & $479,700$& $11,320$ & $\textbf{2,475}$  & $499,700$    & $11,040$     & $>3,600$   & $12,426,000$  & $2779$    & $25,180$ \\
         &      & $10$  & $>3,600$& $311,200$& $12,400$ & $>3,600$  & $320,700$    & $14,790$     & $>3,600$   & $20,950,000$  & $>3,600$   & $30,190$ \\
         &      & $20$  & $1,638$ & $241,600$& $10,710$ & $2,067$  & $279,500$    & $12,760$     & $>3,600$   & $21,470,000$  & $\textbf{148.9}$   & $1,115$ \\\midrule
port $5$ & $225$& $5$   & $0.85$  & $1,489$  & $202$    & $\textbf{0.40}$  & $560$    & $74$     & $>3,600$   & $5,000,000$   & $28.3$    & $22$ \\
         &      & $10$  & $0.60$  & $73$     & $41$     & $\textbf{0.03}$  & $2$    & $5$     & $>3,600$   & $8,989,000$   & $3.33$    & $0$ \\
         &      & $20$  & $0.39$  & $63$     & $52$     & $\textbf{0.08}$  & $0$    & $11$     & $>3,600$   & $10,960,000$  & $115.02$  & $90$ \\
\bottomrule
\end{tabular}
\label{tab:comparisonwithminreturn}
\end{table}

Table \ref{tab:comparisonwithminreturn} indicates that {\color{black}the \verb|port|$4$ instance} cannot be solved to certifiable optimality by any approach within an hour {\color{black} when $k=10$}, in the presence of a minimum return constraint. Nonetheless, both Algorithm \ref{alg:refinedCuttingPlaneMethod} and \verb|CPLEX|'s MISOCO method obtain solutions which are certifiably within $1\%$ of optimality very quickly. Indeed, Table \ref{tab:comparisonwithminreturnbgat120s} depicts the bound gaps of all $3$ approaches at $120$s on these problems; Algorithm \ref{alg:refinedCuttingPlaneMethod} never has a bound gap larger than $0.5\%$.

\begin{table}[h]
\centering\footnotesize
\caption{{\color{black}Bound gap at $120$s per approach with $\kappa=0$, $\gamma=\frac{100}{\sqrt{n}}$ and a minimum return constraint $\bm{\mu}^\top \bm{x} \geq \bar{r}$. We run all approaches on one thread.}}
\begin{tabular}{@{}l l l r r r r r r r r r r r@{}} \toprule
Problem &  $n$ & $k$ & \multicolumn{3}{c@{\hspace{0mm}}}{Algorithm \ref{alg:refinedCuttingPlaneMethod}} & \multicolumn{3}{c@{\hspace{0mm}}}{Algorithm \ref{alg:refinedCuttingPlaneMethod}+in-out} &  \multicolumn{2}{c@{\hspace{0mm}}}{CPLEX Big-M} &  \multicolumn{2}{c@{\hspace{0mm}}}{CPLEX MISOCO} \\
\cmidrule(l){4-6} \cmidrule(l){7-9} \cmidrule(l){10-11} \cmidrule(l){12-13} & &  & Gap ($\%$) & Nodes & Cuts & Gap ($\%$) & Nodes & Cuts & Gap ($\%$) & Nodes & Gap ($\%$) & Nodes \\\midrule
port $2$ & $85$ & $5$   & $\textbf{0}$     & $73,850$  & $1,961$  & $\textbf{0}$ & $42,950$ & $1,272$ & $84.36$   & $611,500$     & $\textbf{0}$       & $1,163$ \\
         &      & $10$  & $0.26$  & $90,670$  & $3,463$  & $0.15$ & $72,240$ & $3,366$ & $425.2$   & $1,057,000$   & $\textbf{0}$       & $902$   \\
         &      & $20$  & $\textbf{0}$     & $12,260$  & $1,224$  & $\textbf{0}$ & $13,650$ & $1,350$ & $65.96$   & $1,367,000$   & $\textbf{0}$       & $210$   \\\midrule
port $3$ & $89$ & $5$   & $0.1$   & $123,100$ & $2,308$  & $\textbf{0.08}$ & $78,950$ & $2,137$ & $88.48$   & $634,300$     & $0.27$    & $1,247$ \\
         &      & $10$  & $0.29$  & $65,180$  & $4,473$ & $0.21$ & $62,840$ & $3,503$ & $452.4$   & $1,073,000$   & $\textbf{0.19}$    & $1,246$ \\
         &      & $20$  & $\textbf{0}$     & $60,090$  & $3,237$  & $\textbf{0}$ & $40,240$ & $3,275$ & $67.55$   & $1,280,000$   & $\textbf{0}$       & $170$   \\\midrule
port $4$ & $98$ & $5$   & $\textbf{0.18}$  & $55,460$  & $3,419$  & $0.37$ & $53,780$ & $3,648$ & $87.67$   & $541,800$     & $0.60$    & $888$   \\
         &      & $10$  & $0.46$  & $51,500$  & $3,704$  & $\textbf{0.29}$ & $46,700$ & $3,241$ & $84.22$   & $1,018,000$   & $\textbf{0.29}$    & $977$   \\
         &      & $20$  & $0.17$  & $57,990$  & $3,393$  & $0.13$ & $59,820$ & $3,886$ & $71.42$   & $1,163,000$   & $\textbf{0.05}$    & $846$   \\
\bottomrule
\end{tabular}
\label{tab:comparisonwithminreturnbgat120s}
\end{table}

The experimental results illustrate that our approach is several orders of magnitude more efficient than the big-$M$ approach on all problems considered, and {\color{black}performs comparably with the} MISOCO approach. Moreover, our approach's edge over \verb|CPLEX| increases with the problem size.

Our main findings from this set of experiments are as follows:
\begin{enumerate}
    \item For problems with unit simplex constraints, big-$M$ approaches do not scale to real-world problem sizes in the presence of ridge regularization, because they cannot exploit the ridge regularizer and therefore obtain low-quality lower bounds, even after expanding a large number of nodes. This poor performance is due to the ridge regularizer; the big-$M$ approach typically exhibits better performance than this in numerical studies done without a regularizer.
    \item MISOCO approaches perform competitively, and are often a computationally reasonable approach for small to medium sized instances of Problem \eqref{mainproblem}, as they are easy to implement and typically have bound gaps of $<1\%$ in instances where they fail to converge within the time budget.
    \item Varying the cardinality of the optimal portfolio does not affect solve times substantially without a minimum return constraint, although it has a nonlinear effect with this constraint.
\end{enumerate}

For the rest of the paper, we do not consider big-$M$ formulations of Problem \eqref{mainproblem}, as they do not scale to larger problems with $200$ or more securities in the universe of buyable assets.

\subsection{Benchmarking Algorithm \ref{alg:refinedCuttingPlaneMethod} in the Presence of Minimum Investment Constraints}\label{sec:buyin}
In this section, we explore Algorithm \ref{alg:refinedCuttingPlaneMethod}'s scalability in the presence of minimum investment constraints, by solving the problems generated by \citet{frangioni2006perspective} and subsequently solved by \cite{frangioni2007sdp, frangioni2009computational, zheng2014improving} among others\footnote{This problem data is available at \href{http://www.di.unipi.it/optimize/Data/MV.html}{www.di.unipi.it/optimize/Data/MV.html}}. These problems have minimum investment, maximum investment, and minimum return constraints, which render many entries in $\mathcal{Z}^n_k$ infeasible. Therefore, to avoid generating an excessive number of feasibility cuts, we use the copy of variables technique suggested in Section \ref{sec:cuttingolanemethodcont}.

Additionally, as the covariance matrices in these problems are highly diagonally dominant (with much larger on-diagonal entries than off-diagonal entries), the method does not converge quickly if we do not extract any diagonal dominance. Indeed, Appendix \ref{sec:suppbuyin} shows that the method often fails to converge within $600$s for the problems studied in this section when we do not extract a diagonally dominant term. Therefore, we first preprocess the covariance matrices to extract more diagonal dominance, as discussed in Section \ref{sec:cuttingolanemethodcont}. Note that we need not actually solve any SDOs to preprocess the data, as high quality diagonal matrices for this problem data have been made publicly available by \citet{frangioni2017improving} at \href{http://www.di.unipi.it/optimize/Data/MV/diagonals.tgz}{http://www.di.unipi.it/optimize/Data/MV/diagonals.tgz} (specifically, we use the entries in the ``s'' folder of this repository). After reading in their diagonal matrix $\bm{D}$, we replace $\bm{\Sigma}$ with $\bm{\Sigma}-\bm{D}$ and use the regularizer $\gamma_i$ for each index $i$, where $\gamma_i=\left(\frac{1}{\gamma}+D_{i,i} \right)^{-1}$ {\color{black}(we do this for all approaches benchmarked here, including \verb|CPLEX|'s MISOCO routine).}

We now compare the times for Algorithm \ref{alg:refinedCuttingPlaneMethod} and \verb|CPLEX|'s MISOCO routines to solve the diagonally dominant instances in the dataset generated by \citet{frangioni2006perspective}, along with a variant of Algorithm \ref{alg:refinedCuttingPlaneMethod} where we use the \verb|in-out| method at the root node, and another variant where we apply the \verb|in-out| method at both the root node and $50$ additional nodes. In all cases, we take $\gamma=\frac{1000}{n}$, which ensures that $\gamma_i \approx \frac{1}{D_{i,i}}${\color{black}, since on this dataset $\frac{1}{2\gamma}$ is around $5$ orders of magnitude smaller than $D_{i,i}$ and thus the net contribution of the regularization term to the objective is negligible}. Table \ref{tab:comparisonwithminreturn_fullpard} depicts the average time taken by each approach, and demonstrates that Algorithm \ref{alg:refinedCuttingPlaneMethod} substantially outperforms \verb|CPLEX|, particularly for problems without a{\color{black}n explicit} cardinality constraint{\color{black}, but with an implicit cardinality constraint of $\Vert \bm{x}\Vert_0 \leq k$ for $k \leq 13$ due to the minimum investment constraints $x_i \geq l_i z_i$ where $l_i$ is uniformly drawn from $[0.075, 0.125]$ \citep{frangioni2007sdp}}. We provide the full instance-wise results in Appendix \ref{sec:suppbuyin}{\color{black}, and also consider the possibility of supplying the cuts from the \verb|in-out| method to \verb|CPLEX| before running the MISOCP method therein.}

\begin{table}[h]
\centering\footnotesize
\caption{Average runtime in seconds per approach with $\kappa=0$, $\gamma=\frac{1000}{n}$ for the problems generated by \citet{frangioni2006perspective}. We impose a time limit of $600$s and run all approaches on one thread. If a solver fails to converge, we report the number of explored nodes at the time limit, use $600$s in lieu of the solve time, and report the number of failed instances (out of $10$) next to the solve time in brackets. {\color{black}Note that the minimum investment constraints impose an implicit cardinality constraint with $k \approx {\color{black}13}$.}}
\begin{tabular}{@{}l l r r r r r r r r r r r r@{}} \toprule
Problem & $k$ & \multicolumn{3}{c@{\hspace{0mm}}}{Algorithm \ref{alg:refinedCuttingPlaneMethod}} & \multicolumn{3}{c@{\hspace{0mm}}}{Algorithm \ref{alg:refinedCuttingPlaneMethod} $+$ in-out} & \multicolumn{3}{c@{\hspace{0mm}}}{Algorithm \ref{alg:refinedCuttingPlaneMethod} in-out $+$ 50} &  \multicolumn{2}{c@{\hspace{0mm}}}{CPLEX MISOCO} \\
\cmidrule(l){3-5} \cmidrule(l){6-8} \cmidrule(l){9-11} \cmidrule(l){12-13}  & & Time & Nodes & Cuts & Time & Nodes & Cuts & Time & Nodes & Cuts & Time & Nodes \\\midrule
200+ & 6 & \textbf{1.55} & 1,298 & 236.3 & 1.77 & 1,262 & 209.4 & 7.4 & 910.4 & 118 & 87.74 (0)& 95.3\\
200+ & 8 & \textbf{1.95} & 1,968 & 260.3 & 2.30 & 1,626 & 217 & 7.97 & 949.1 & 97.3 & 73.42 (0)& 79.8\\
200+ & 10 & 7.74 & 7,606 & 509.7 & \textbf{4.33} & 3,686 & 298.9 & 10.35 & 2,066 & 175.5 & 161.9 (0)& 184\\
200+ & 12 & 25.57 & 28,830 & 203.8 & \textbf{2.06} & 1,764 & 71.6 & 9.04 & 1,000 & 33.9 & 353.1 (4) & 398.1\\
200+ & 200 & 18.71 & 23,190 & 208.4 & \textbf{2.79} & 2,288 & 92 & 10 & 1,394 & 56.1 & 599.3 (9) & 735.1\\\midrule
300+ & 6 & \textbf{16.83} & 9,141 & 974.2 & 23.59 & 8,025 & 864.1 & 29.92 & 5,738 & 565.9 & 434.5 (3)& 157.6\\
300+ & 8 & \textbf{44.68} & 21,050 & 1,577 & 64.46 & 19,682 & 1457.8 & 61.0 & 14,236 & 1,036 & 489.5 (5)& 174.0\\
300+ & 10 & 88.57 & 44,160 & 1,901 & \textbf{78.05} & 33,253 & 1438.4 & 110.2 & 24,487 & 971.5 & 472.0 (5)& 171.9\\
300+ & 12 & 16.16 & 13,880 & 262.7 & \textbf{4.65} & 3,181 & 127.4 & 15.94 & 1475 & 66.7 & 401.5 (4)& 158.2\\
300+ & 300 & 21.36 & 18,140 & 262.1 & \textbf{9.24} & 6,288 & 191.9 & 24.33 & 5,971 & 168.4 & 600.0 (10)& 219.2\\\midrule
400+ & 6 & \textbf{54.47} & 13,330 & 1,717 &66.52 & 12,160 & 1,619 & 85.51 & 11,070 & 1,402 & 531.7 (8)& 84.0\\
400+ & 8 & 173.8 & 35,390 & 2,828 & \textbf{160.9} & 32,930 & 2,709 & 163.3 & 28,020 & 2,363 & 534.0 (8)& 80.8\\
400+ & 10 & 158.0 & 55,490 & 1,669 & 104.5 & 32,314 & 1369.7 & \textbf{81.48} & 22,130 & 824.9 & 517.9 (8)& 74.8\\
400+ & 12 & 3.97 & 4,324 & 116.6 & \textbf{1.9} & 1,214 & 48.6 & 15.67 & 627.4 & 29.8 & 478.0 (4)& 75.3\\
400+ & 400 & 8.68 & 7,540 & 120.5 & \textbf{5.19} & 3,539 & 88.8 & 21.31 & 3,210 & 79.4 & 600.0 (10)& 74.2\\
\bottomrule
\end{tabular}
\label{tab:comparisonwithminreturn_fullpard}
\end{table}

Our main findings from this experiment are as follows:
\begin{itemize}
    \item Algorithm \ref{alg:refinedCuttingPlaneMethod} outperforms \verb|CPLEX| in the presence of minimum investment constraints. {\color{black}Interestingly, the log output indicates that \verb|CPLEX|'s custom heuristics typically obtain incumbent solutions which are better than those obtained by Algorithm \ref{alg:warmstart} and supplied to Algorithm \ref{alg:refinedCuttingPlaneMethod} as a warm-start. This suggests that Algorithm \ref{alg:refinedCuttingPlaneMethod}'s superior numerical performance arises because it can develop high-quality (i.e., better than the SOCO bound) lower bounds more quickly. This shouldn't be too surprising, since beating a continuous relaxation bound involves branching (or cutting), which is cheaper to perform on a pure binary formulation than a MISOCO formulation. }
    \item Running the \verb|in-out| method at the root node improves solve times when $k \geq 10$, but does more harm than good when $k<10$, because in the latter case Algorithm \ref{alg:refinedCuttingPlaneMethod} already performs well.
    \item Running the \verb|in-out| method at non-root nodes does more harm than good for easy problems, but improves solve times for larger problems ($400+$ with $k \in \{8, 10\}$), as reported in Appendix \ref{sec:suppbuyin}.
    \item With a cardinality constraint, Algorithm \ref{alg:refinedCuttingPlaneMethod}'s solve times are comparable to those reported by \citet{zheng2014improving, frangioni2016approximated} {\color{black}(note however that this comparison is imperfect since all three approaches were run on different machines; the source code for the other two approaches is unavailable and hence we cannot obtain truly comparable data}). {\color{black}This can be explained by the fact that all three methods solve these problems in $10$s of seconds, and thus these problems can be viewed as ``easy''.} Without an {\color{black}explicit} cardinality constraint {\color{black}(but with minimum investment constraints which impose an implicit cardinality constraint {\color{black}of $k \approx 13$ \cite{frangioni2007sdp}})}, our solve times are two orders of magnitude faster than those reported by \citet{zheng2014improving}'s {\color{black} (an average of $580$s for $400+$)}, and an order of magnitude faster than those by \citet{frangioni2016approximated} {\color{black} (an average of $52$s for $400+$)}.
    \item As shown in Appendix \ref{sec:suppbuyin}, applying the diagonal dominance preprocessing technique proposed by \citet{frangioni2007sdp} yields faster solve times than applying the technique proposed by \citet{zheng2014improving}, even though the latter technique yields tighter continuous relaxations \citep{zheng2014improving}. This might occur because \citet{frangioni2007sdp}'s technique prompts our approach to make better branching decisions and/or \citet{zheng2014improving}'s approach is only guaranteed to yield tighter continuous relaxations before (i.e., not after) branching.
     \item {\color{black}As shown in Appendix \ref{sec:suppbuyin}, passing the cuts generated by the \verb|in-out| method to \verb|CPLEX| does more harm than good. This can be explained by the fact that the SOCP relaxation, which is very affordable to compute, is at least as strong as the lower approximation generated by the \verb|in-out| method (indeed, the \verb|in-out| method can be interpreted as a linearization of the SOCP relaxation about a ``cleverly'' sampled set of points), while providing the in-out method's lower approximation to \verb|CPLEX| increases the amount of work which \verb|CPLEX| needs to perform at each node in the branch-and-bound tree.}
\end{itemize}

\FloatBarrier
\subsection{Exploring the Scalability of Algorithm \ref{alg:refinedCuttingPlaneMethod}}
In this section, we explore Algorithm \ref{alg:refinedCuttingPlaneMethod}'s scalability with respect to the number of securities in the buyable universe, by measuring the time required to solve several large-scale sparse portfolio selection problems to provable optimality: the S$\&$P $500$, the Russell $1000$, and the Wilshire $5000$. In all three cases, the problem data is taken from daily closing prices from January $3$ $2007$ to December $29$ $2017$, which are obtained from Yahoo! Finance via the R package \textit{quantmod} (see \cite{ryan2018quantmod}), and rescaled to correspond to a holding period of one month. We apply Singular Value Decomposition to obtain low-rank estimates of the correlation matrix, and rescale the low-rank correlation matrix by each asset's variance to obtain a low-rank covariance matrix $\bm{\Sigma}$. We also omit days with a greater than $20\%$ change in closing prices when computing the mean and covariance for the Russell $1000$ and Wilshire $5000$, since these changes occur on low-volume trading and typically reverse the next day.

Tables \ref{tab:benchmarksp500}--\ref{tab:benchmarkwilshire5000} depict the times required for Algorithm \ref{alg:refinedCuttingPlaneMethod} and \verb|CPLEX MISOCO| to solve the problem to provable optimality for different choices of $\gamma$, $k$, and Rank$(\bm{\Sigma})$. In particular, they depict the time taken to solve (a) an unconstrained problem where $\kappa=1$ and (b) a constrained problem where $\kappa=0$ containing a minimum return constraint computed in the same fashion as in Section \ref{ssec:compwithbigm}.

\begin{table}[h]
\centering\footnotesize
\caption{Runtimes in seconds per approach for the S$\&$P $500$ with $\kappa=1$ (left); $\kappa=0$ and a minimum return constraint (right), a one-month holding period and a runtime limit of $600$s. For instances with a minimum return constraint where $\gamma=\frac{100}{\sqrt{n}}$, we run the in-out method at the root node before running Algorithm \ref{alg:refinedCuttingPlaneMethod}. We run all approaches on one thread. When a method fails to converge, we report the bound gap at $600$s.}
\begin{tabular}{@{}l l r r r r r r r r r r r@{}} \toprule
$\gamma$ & Rank$(\bm{\Sigma})$ & $k$ & \multicolumn{3}{c@{\hspace{0mm}}}{Algorithm \ref{alg:refinedCuttingPlaneMethod}} &  \multicolumn{2}{c@{\hspace{0mm}}}{CPLEX MISOCO} & \multicolumn{3}{c@{\hspace{0mm}}}{Algorithm \ref{alg:refinedCuttingPlaneMethod}} &  \multicolumn{2}{c@{\hspace{0mm}}}{CPLEX MISOCO} \\
\cmidrule(l){4-6} \cmidrule(l){7-8} \cmidrule(l){9-11} \cmidrule(l){12-13}  & & & Time & Nodes & Cuts & Time & Nodes & Time & Nodes & Cuts & Time & Nodes \\\midrule
$\frac{1}{\sqrt{n}}$ & $50$ & $10$ & $0.01$ & $0$ & $4$& $0.54$ & $0$ & $0.01$ & $0$ & $3$& $73.28$ & $210$\\
& & $50$ & $0.02$ & $0$ & $4$& $0.49$ & $0$ & $0.28$ & $108$ & $45$& $78.59$ & $499$\\
& & $100$ & $0.03$ & $0$ & $4$ & $1.00$ & $0$ & $0.05$ & $7$ & $7$& $0.97$ & $0$\\
& & $200$ & $0.06$ & $0$ & $4$& $0.86$ & $0$ & $0.08$ & $1$ & $5$& $53.53$ & $300$\\\midrule
$\frac{1}{\sqrt{n}}$ & $100$ & $10$ & $0.01$ & $0$ & $4$& $1.49$ & $0$ & $2.01$ & $972$ & $344$& $339.8$ & $420$\\
& & $50$ & $0.02$ & $0$ & $4$& $1.36$ & $0$ & $0.32$ & $104$ & $41$& $283.8$ & $410$\\
& & $100$ & $0.04$ & $0$ & $4$& $1.30$ & $0$ & $0.06$ & $5$ & $7$& $286.2$ & $520$\\
& & $200$ & $0.09$ & $0$ & $4$& $3.10$ & $0$ & $0.06$ & $0$ & $3$& $472.7$ & $990$\\\midrule
$\frac{1}{\sqrt{n}}$ & $150$ & $10$ & $0.01$ & $0$ & $4$& $2.61$ & $0$ & $3.96$ & $1,633$ & $410$& $268.3$ & $157$\\
& & $50$ & $0.03$ & $0$ & $4$& $2.23$ & $0$ & $0.29$ & $62$ & $33$& $265.6$ & $200$\\
& & $100$ & $0.06$ & $0$ & $4$& $4.71$ & $0$ & $0.07$ & $0$ & $6$& $394.9$ & $340$\\
& & $200$ & $0.14$ & $0$ & $4$& $4.80$ & $0$ & $0.13$ & $0$ & $3$& $6.20$ & $0$\\\midrule
$\frac{1}{\sqrt{n}}$ & $200$ & $10$ & $0.01$ & $0$ & $4$& $2.74$ & $0$ & $5.20$ & $2,804$ & $450$& $345.0$ & $171$\\
& & $50$ & $0.03$ & $0$ & $4$& $3.14$ & $0$ & $0.49$ & $86$ & $47$& $337.7$ & $210$\\
& & $100$ & $0.06$ & $0$ & $4$& $17.27$ & $3$ & $0.15$ & $5$ & $8$& $104.2$ & $40$\\
& & $200$ & $0.13$ & $0$ & $4$& $105.2$ & $60$ & $0.10$ & $0$ & $3$& $46.18$ & $10$\\\midrule
$\frac{100}{\sqrt{n}}$ & $50$ & $10$ & $0.01$ & $0$ & $4$& $0.51$ & $0$ & $0.09\%$ & $70,200$ & $3,855$& $0.10\%$ & $1,600$\\
& & $50$ & $0.02$ & $0$ & $4$ & $1.14$ & $0$ & $0.77$ & $309$ & $113$& $268.5$ & $841$\\
& & $100$ & $0.03$ & $0$ & $4$ & $0.68$ & $0$ & $0.09$ & $0$ & $8$& $1.66$ & $0$\\
& & $200$ & $0.07$ & $0$ & $4$& $1.04$ & $0$ & $0.16$ & $0$ & $4$& $15.26$ & $10$\\ \midrule
$\frac{100}{\sqrt{n}}$ & $100$ & $10$ & $0.01$ & $0$ & $4$ & $1.30$ & $0$ & $0.26\%$ & $54,000$ & $4,721$& $0.24\%$ & $598$\\
& & $50$ & $0.03$ & $0$ & $4$& $3.48$ & $0$ & $0.07$ & $1$ & $7$& $0.28\%$ & $291$\\
& & $100$ & $0.06$ & $0$ & $5$& $5.93$ & $0$ & $0.11$ & $0$ & $4$& $0.29\%$ & $352$\\
& & $200$ & $0.13$ & $0$ & $5$& $2.16$ & $0$ & $0.14$ & $0$ & $5$& $301.3$ & $380$\\\midrule
$\frac{100}{\sqrt{n}}$ & $150$ & $10$ & $0.02$ & $0$ & $4$& $2.00$ & $0$ & $0.33\%$ & $46,720$ & $4,437$& $0.28\%$ & $345$\\
& & $50$ & $0.04$ & $0$ & $4$& $34.57$ & $20$ & $0.09$ & $0$ & $6$& $0.28\%$ & $291$\\
& & $100$ & $0.06$ & $0$ & $4$& $27.76$ & $20$ & $0.11$ & $0$ & $4$& $0.29\%$ & $352$\\
& & $200$ & $0.10$ & $0$ & $4$& $26.93$ & $20$ & $0.35$ & $0$ & $6$& $344.3$ & $270$\\\midrule
$\frac{100}{\sqrt{n}}$ & $200$ & $10$ & $0.02$ & $0$ & $4$& $7.77$ & $0$ & $0.45\%$ & $56,100$ & $4,336$& $0.36\%$ & $280$\\
& & $50$ & $0.04$ & $0$ & $4$& $48.75$ & $20$ & $0.20$ & $1$ & $19$& $0.35\%$ & $256$\\
& & $100$ & $0.10$ & $0$ & $5$& $44.02$ & $20$ & $0.15$ & $0$ & $5$& $104.2$ & $40$\\
& & $200$ & $0.16$ & $0$ & $4$& $36.57$ & $20$ & $0.18$ & $0$ & $4$& $76.80$ & $10$\\
\bottomrule
\end{tabular}
\label{tab:benchmarksp500}
\end{table}

\begin{table}[h]
\centering\footnotesize
\caption{
Runtimes in seconds per approach for the Russell $1000$ with $\kappa=1$ (left); $\kappa=0$ and a minimum return constraint (right), a one-month holding period and a runtime limit of $600$s. For instances with a minimum return constraint, we run the in-out method at the root node before running Algorithm \ref{alg:refinedCuttingPlaneMethod}. We run all approaches on one thread. When a method fails to converge, we report the bound gap at $600$s.}
\begin{tabular}{@{}l l r r r r r r r r r r r@{}} \toprule
$\gamma$ & Rank$(\bm{\Sigma})$ & $k$ & \multicolumn{3}{c@{\hspace{0mm}}}{Algorithm \ref{alg:refinedCuttingPlaneMethod}} &  \multicolumn{2}{c@{\hspace{0mm}}}{CPLEX MISOCO} & \multicolumn{3}{c@{\hspace{0mm}}}{Algorithm \ref{alg:refinedCuttingPlaneMethod}} & \multicolumn{2}{c@{\hspace{0mm}}}{CPLEX MISOCO} \\
\cmidrule(l){4-6} \cmidrule(l){7-8} \cmidrule(l){9-11} \cmidrule(l){12-13}  & & & Time & Nodes & Cuts & Time & Nodes & Time & Nodes & Cuts & Time & Nodes \\\midrule
$\frac{1}{\sqrt{n}}$ & $50$ & $10$ & $0.02$ & $0$ & $6$& $7.77$ & $3$ & $12.38$ & $2467$ & $316$ & $0.01\%$ & $545$\\
& & $50$ & $0.06$ & $0$ & $12$ & $9.19$ & $7$ & $0.32$ & $0$ & $7$ & $0.01\%$ & $900$\\
& & $100$ & $0.04$ & $0$ & $5$& $0.92$ & $0$ & $0.81$ & $10$ & $14$ & $0.01\%$ & $1,048$\\
& & $200$ & $0.07$ & $0$ & $5$& $1.83$ & $0$ & $0.46$ & $1$ & $14$ & $0.01\%$ & $1,043$\\\midrule
$\frac{1}{\sqrt{n}}$ & $100$ & $10$ & $0.03$ & $3$ & $7$ & $13.27$ & $5$ & $272.3$ & $49,200$ & $1,266$& $0.01\%$ & $400$\\
& & $50$ & $0.09$ & $2$ & $12$ & $154.0$ & $90$ & $1.32$ & $10$ & $13$& $0.01\%$ & $400$\\
& & $100$ & $0.05$ & $0$ & $5$& $2.83$ & $0$ & $15.52$ & $5,271$ & $250$& $0.01\%$ & $599$\\
& & $200$ & $0.14$ & $0$ & $8$& $335.9$ & $260$ & $2.12$ & $111$ & $64$& $0.01\%$ & $399$\\\midrule
$\frac{1}{\sqrt{n}}$ & $200$ & $10$ & $0.05$ & $2$ & $10$ & $41.08$ & $7$ & $24.14$ & $8,200$ & $318$& $31.20\%$ & $138$\\
& & $50$ & $0.16$ & $8$ & $14$& $344.0$ & $60$ & $0.01\%$ & $86,020$ & $1,433$& $0.02\%$ & $100$\\
& & $100$ & $0.08$ & $0$ & $5$ & $60.54$ & $10$ & $4.88$ & $131$ & $41$& $0.01\%$ & $100$\\
& & $200$ & $0.16$ & $0$ & $5 $& $6.79$ & $0$ & $0.01\%$ & $64,600$ & $1,049$ & $4.00\%$ & $100$\\\midrule
$\frac{1}{\sqrt{n}}$ & $300$ & $10$ & $0.06$ & $1$ & $10$& $175.9$ & $15$ & $9.00$ & $1,200$ & $246$& $0.01\%$ & $70$\\
& & $50$ & $0.16$ & $2$ & $13$& $323.3$ & $31$ & $0.02\%$ & $61,100$ & $1,227$& $63.35\%$ & $78$\\
& & $100$ & $0.16$ & $0$ & $8$ & $260.4$ & $30$ & $0.02\%$ & $48,550$ & $856$& $3.01\%$ & $64$\\
& & $200$ & $0.31$ & $0$ & $8$ & $0.01\%$ & $464$ & $0.01\%$ & $29,480$ & $786$& $6.00\%$ & $75$\\\midrule
$\frac{100}{\sqrt{n}}$ & $50$ & $10$ & $0.04$ & $1$ & $11$& $7.58$ & $3$ & $0.59\%$ & $62,050$ & $2,553$& $0.39\%$ & $700$\\
& & $50$ & $0.04$ & $0$ & $9$ & $2.34$ & $0$ & $0.61\%$ & $114,000$ & $1,531$& $0.32\%$ & $837$\\
& & $100$ & $0.36$ & $0$ & $4$ & $2.57$ & $0$ & $112.1$ & $42,661$ & $787$ & $0.12\%$ & $993$\\
& & $200$ & $0.09$ & $0$ & $5$ & $1.25$ & $0$ & $0.40$ & $0$ & $19$& $135.4$ & $220$\\\midrule
$\frac{100}{\sqrt{n}}$ & $100$ & $10$ & $0.03$ & $2$ & $9$ & $11.50$ & $5$ & $0.77\%$ & $69,800$ & $2,599$& $0.55\%$ & $400$\\
& & $50$ & $0.06$ & $0$ & $8$& $65.82$ & $40$ & $0.68\%$ & $93,580$ & $1,472$& $0.38\%$ & $400$\\
& & $100$ & $0.06$ & $0$ & $5$& $22.82$ & $10$ & $0.43\%$ & $82,500$ & $1,359$& $0.46\%$ & $470$\\
& & $200$ & $0.11$ & $0$ & $5$& $42.68$ & $30$ & $0.34$ & $0$ & $13$& $0.37\%$ & $400$\\\midrule
$\frac{100}{\sqrt{n}}$ & $200$ & $10$ & $0.06$ & $1$ & $10$ & $31.33$ & $0$ & $0.84\%$ & $84,617$ & $2,183$& $1.23\%$ & $126$\\
& & $50$ & $0.10$ & $1$ & $10$ & $164.4$ & $30$ & $1.55\%$ & $99,600$ & $1,576$& $1.31\%$ & $100$\\
& & $100$ & $0.08$ & $0$ & $4$ & $71.91$  & $10$ & $0.92\%$ & $69,850$ & $1,279$& $9.36\%$ & $118$\\
& & $200$ & $0.16$ & $0$ & $5$ & $50.98$ & $10$ & $0.35$ & $1$ & $10$& $2.20\%$ & $123$\\\midrule
$\frac{100}{\sqrt{n}}$ & $300$ & $10$ & $0.06$ & $1$ & $12$& $134.1$ & $15$ & $0.94\%$ & $72,400$ & $2,759$ & $1.00\%$ & $65$\\
& & $50$ & $0.10$ & $0$ & $8$ & $207.7$ & $20$ & $1.78\%$ & $58,740$ & $1,363$& $48.31\%$ & $62$\\
& & $100$ & $0.10$ & $0$ & $4$ & $544.3$ & $50$ & $1.16\%$ & $55,810$ & $1,027$& $14.40\%$ & $61$\\
& & $200$ & $0.20$ & $0$ & $5$ & $221.4$ & $30$ & $1.04\%$ & $61,410$ & $1,122$& $2.92\%$ & $61$\\
\bottomrule
\end{tabular}
\label{tab:benchmarkruss1000}
\end{table}

\begin{table}[h]
\centering\footnotesize
\caption{Runtimes in seconds per approach for the Wilshire $5000$ with $\kappa=1$ (left); $\kappa=0$ and a minimum return constraint (right), a one-month holding period and a runtime limit of $600$s. For instances with a minimum return constraint where $\gamma=\frac{100}{\sqrt{n}}$, we run the in-out method at the root node before running Algorithm \ref{alg:refinedCuttingPlaneMethod}. We run all approaches on one thread. When a method fails to converge, we report the bound gap at $600$s (using the symbol ``$-$'' to denote that a method failed to produce a feasible solution).}
\begin{tabular}{@{}l l r r r r r r r r r r r@{}} \toprule
$\gamma$ & Rank$(\bm{\Sigma})$ & $k$ & \multicolumn{3}{c@{\hspace{0mm}}}{Algorithm \ref{alg:refinedCuttingPlaneMethod}} &  \multicolumn{2}{c@{\hspace{0mm}}}{CPLEX MISOCO} & \multicolumn{3}{c@{\hspace{0mm}}}{Algorithm \ref{alg:refinedCuttingPlaneMethod}} &  \multicolumn{2}{c@{\hspace{0mm}}}{CPLEX MISOCO} \\
\cmidrule(l){4-6} \cmidrule(l){7-8} \cmidrule(l){9-11} \cmidrule(l){12-13}  & & & Time & Nodes & Cuts & Time & Nodes & Time & Nodes & Cuts & Time & Nodes \\\midrule
$\frac{1}{\sqrt{n}}$ & $100$ & $10$ & $0.04$ & $0$ & $4$ & $15.07$ & $0$ & $1.95$ & $0$ & $2$ & $50.0\%$ & $122$\\
& & $50$ & $0.07$ & $0$ & $10$ & $244.8$ & $29$ & $2.32$ & $0$ & $2$ & $32.0\%$ & $132$ \\
& & $100$ & $0.22$ & $0$ & $4$ & $30.08$ & $2$ & $0.59$ & $10$ & $9$ & $62.0\%$ & $127$ \\
& & $200$ & $0.24$ & $0$ & $4$ & $40.54$ & $3$ & $0.27$ & $0$ & $6$ & $44.5\%$ & $100$\\\midrule
$\frac{1}{\sqrt{n}}$ & $200$ & $10$ & $0.07$ & $0$ & $6$ & $70.12$ & $2$ & $2.34$ & $0$ & $8$ & $30.0\%$ & $43$ \\
& & $50$ & $0.08$ & $0$ & $8$ & $392.5$ & $25$ & $5.54$ & $31$ & $17$ & $74.0\%$ & $40$\\
& & $100$ & $0.08$ & $0$ & $5$ & $0.01\%$ & $91$ & $1.70$ & $0$ & $12$ & $62.0\%$ & $44$ \\
& & $200$ & $0.25$ & $0$ & $4$ & $49.26$ & $0$ & $0.01\%$ & $8,451$ & $361$ & $41.5\%$ & $40$ \\\midrule
$\frac{1}{\sqrt{n}}$ & $500$ & $10$ & $0.15$ & $10$ & $11$ & $0.01\%$ & $13$ & $10.53$ & $300$ & $66$ & $-$ & $5$ \\
& & $50$ & $0.54$ & $0$ & $8$ & $0.01\%$ & $30$ & $0.01\%$ & $49,000$ & $805$ & $-$ & $6$ \\
& & $100$ & $0.20$ & $0$ & $6$ & $492.7$ & $3$ & $0.01\%$ & $36,670$ & $1,068$ & $-$ & $5$ \\
& & $200$ & $0.57$ & $0$ & $5$ & $0.01\%$ & $20$ & $3.41$ & $0$ & $14$ & $-$ & $5$ \\\midrule
$\frac{1}{\sqrt{n}}$ & $1,000$ & $10$ & $0.48$ & $35$ & $28$ & $0.01\%$ & $9$ & $0.01\%$ & $40,500$ & $1,130$ & $-$ & $2$ \\
& & $50$ & $1.08$ & $20$ & $29$ & $0.01\%$ & $9$ & $0.02\%$ & $56,800$ & $937$ & $-$ & $2$ \\
& & $100$ & $0.44$ & $0$ & $7$ & $0.01\%$ & $11$ & $0.02\%$ & $25,040$ & $523$ & $-$ & $2$ \\
& & $200$ & $0.56$ & $0$ & $4$ & $0.01\%$ & $10$ & $2.61$ & $1$ & $12$ & $-$ & $2$ \\\midrule
$\frac{100}{\sqrt{n}}$ & $100$ & $10$ & $0.03$ & $0$ & $5$ & $29.30$ & $2$ & $0.28\%$ & $24,870$ & $1,178$ & $50.1\%$ & $91$ \\
& & $50$ & $0.04$ & $0$ & $5$ & $39.08$ & $3$ & $0.38\%$ & $45,810$ & $636$ & $62.1\%$ & $82$ \\
& & $100$ & $0.07$ & $0$ & $5$ & $200.0$ & $11$ & $0.12\%$ & $55,700$ & $912$ & $45.1\%$ & $80$ \\
& & $200$ & $0.10$ & $0$ & $10$ & $99.07$ & $10$ & $0.49$ & $0$ & $10$ & $22.1\%$ & $91$ \\\midrule
$\frac{100}{\sqrt{n}}$ & $200$ & $10$ & $0.06$ & $0$ & $8$ & $56.48$ & $2$ & $0.38\%$ & $34,100$ & $1,071$ & $-$ & $29$ \\
& & $50$ & $0.08$ & $0$ & $7$ & $78.00$ & $3$ & $0.47\%$ & $40,340$ & $1,034$ & $66.1\%$ & $30$ \\
& & $100$ & $0.20$ & $0$ & $5$ & $0.01\%$ & $20$ & $0.43\%$ & $15,010$ & $325$ & $45.1\%$ & $33$ \\
& & $200$ & $0.14$ & $0$ & $4$ & $224.4$ & $10$ & $0.98$ & $6$ & $10$ & $20.1\%$ & $30$ \\\midrule
$\frac{100}{\sqrt{n}}$ & $500$ & $10$ & $0.15$ & $5$ & $13$ & $0.01\%$ & $16$ & $0.52\%$ & $32,920$ & $1,235$ & $-$ & $4$ \\
& & $50$ & $0.08$ & $0$ & $4$ & $0.01\%$ & $6$ & $1.11\%$ & $65,560$ & $771$ & $-$ & $3$ \\
& & $100$ & $0.27$ & $0$ & $8$ & $0.01\%$ & $10$ & $0.79\%$ & $23,540$ & $651$ & $-$ & $2$ \\
& & $200$ & $0.30$ & $0$ & $4$ & $0.01\%$ & $10$ & $0.52$ & $0$ & $6$ & $-$ & $2$ \\\midrule
$\frac{100}{\sqrt{n}}$ & $1,000$ & $10$ & $0.48$ & $29$ & $32$ & $0.17\%$ & $10$ & $1.02\%$ & $6,7600$ & $1,108$ & $-$ & $2$ \\
& & $50$ & $0.26$ & $0$ & $8$ & $0.01\%$ & $9$ & $1.74\%$ & $33,930$ & $1,122$ & $-$ & $0$ \\
& & $100$ & $0.56$ & $0$ & $8$ & $0.01\%$ & $9$ & $1.85\%$ & $53,500$ & $804$ & $-$ & $2$ \\
& & $200$ & $0.66$ & $0$ & $4$ & $0.01\%$ & $9$ & $1.28$ & $1$ & $7$ & $-$ & $2$ \\
\bottomrule
\end{tabular}
\label{tab:benchmarkwilshire5000}
\end{table}

Our main finding from this set of experiments is that Algorithm \ref{alg:refinedCuttingPlaneMethod} is substantially faster than \verb|CPLEX|'s MISOCO routine, particularly as the rank of $\bm{\Sigma}$ increases. The relative numerical success of Algorithm \ref{alg:refinedCuttingPlaneMethod} in this section, compared to the previous section, can be explained by the differences in the problems solved: (a) in this section, we optimize over a sparse unit simplex, while in the previous section we optimized over minimum-return and minimum-investment constraints, (b) in this section, we use data taken directly from stock markets, while in the previous section we used less realistic synthetic data, which evidently made the problem harder.

\subsection{Exploring Sensitivity to Hyperparameters}
Our next set of experiments explores Problem \eqref{mainproblem}'s stability to changes in its hyperparameter $\gamma$.

We first explore optimizing over a rank--$300$ approximation of the Russell $1000$ with a one month holding period, a sparsity budget $k=10$ and a weight $\kappa=1$.

\begin{figure}[h]
\centering
\caption{Sensitivity to $\gamma$ for the Russell $1000$ with $\kappa=1$ and $k=10$. The optimal security indices $\bm{z}^\star$ changed once over the entire range of $\gamma$.}
\includegraphics[scale=0.4]{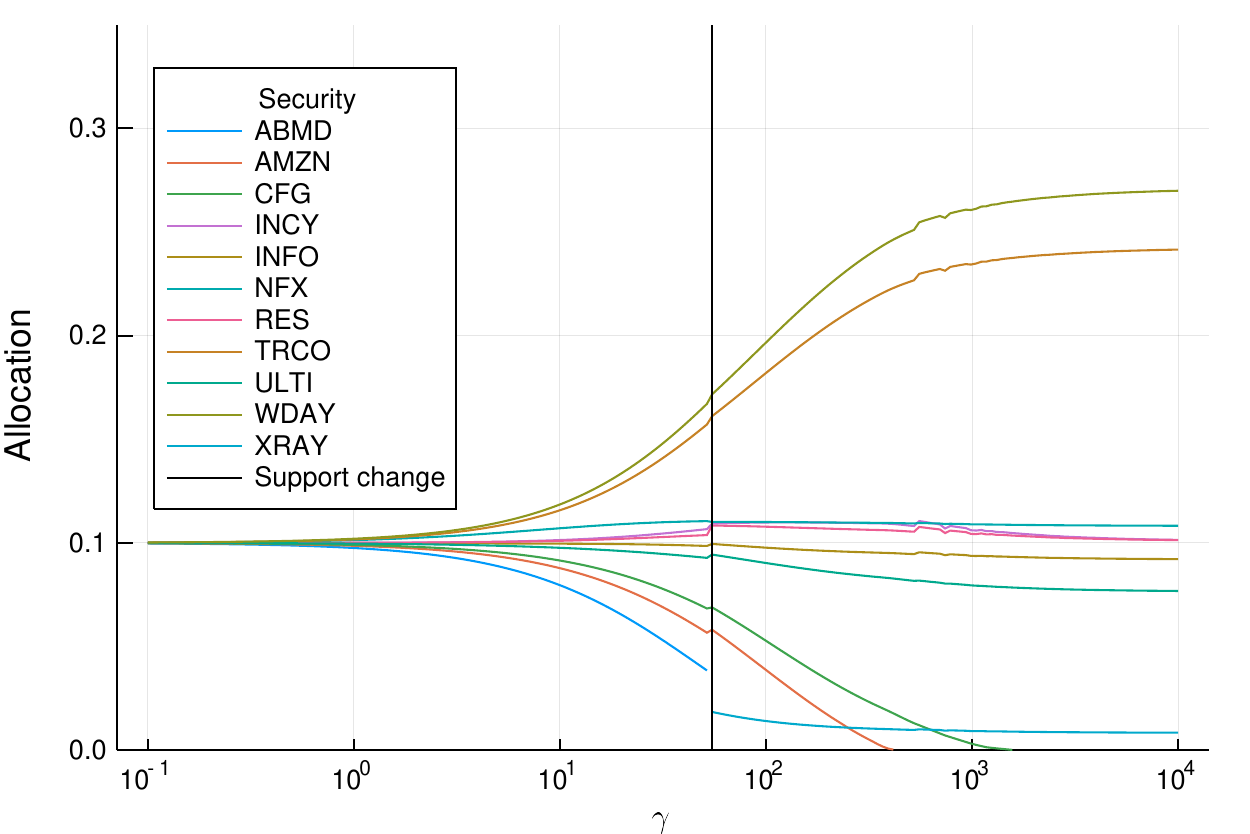}
\label{fig:russ1000sensitivity}
\end{figure}

Figure \ref{fig:russ1000sensitivity} depicts the relationship between $\bm{x}^\star$ and $\gamma$ for this set of hyperparameters, and indicates that $\bm{x}^\star$ is stable with respect to small changes in $\gamma$. Moreover, the optimal support indices when $\gamma$ is small {\color{black}(to the left of the vertical line)} are near-optimal when $\gamma$ is large {\color{black}, while the optimal support indicies when $\gamma$ is large {\color{black}(to the right of the vertical line)} are optimal for Problem \eqref{mainproblemnonreg}, which echoes our sensitivity analysis findings and particularly Proposition \ref{prop:largestgamma}}. This suggests that a good strategy for cross-validating $\gamma$ could be to solve Problem \eqref{mainproblem} to certifiable optimality for one value of $\gamma$, find the best value of $\gamma$ conditional on using these support indices, and finally resolve Problem \eqref{mainproblem} with the optimal $\gamma$.

{\color{black} Our final experiment studies the impact of the regularizer $\gamma$ on both solve times and the number of cuts generated, in order to justify our assertion in the introduction that increasing the amount of regularization in the problem makes the problem easier. In this direction, we solve the ten $300+$ and $400+$ instances with minimum investment and minimum return constraints studied in Section \ref{sec:buyin} for different values of $\gamma$ (with the copy of variables technique and in-out method on). We report the average runtime and number of cuts generated by Algorithm \ref{alg:refinedCuttingPlaneMethod} across the $10$ instances for each $n$ in Figures \ref{fig:gammavariationpard300} ($n=300$) and \ref{fig:gammavariationpard400} ($n=400$).

Observe that the average runtime is essentially non-decreasing in $\gamma$, and for both {\color{black}$n=300$ and $n=400$} there exists a finite $\gamma>0$ at which all instances can be solved using a single cut. This empirically verifies Section \ref{ssec:sensitivityanalysis}'s sensitivity analysis findings.

\begin{figure}[h!]
\centering
\begin{subfigure}[t]{.47\linewidth}
	\centering
	\includegraphics[width=\linewidth]{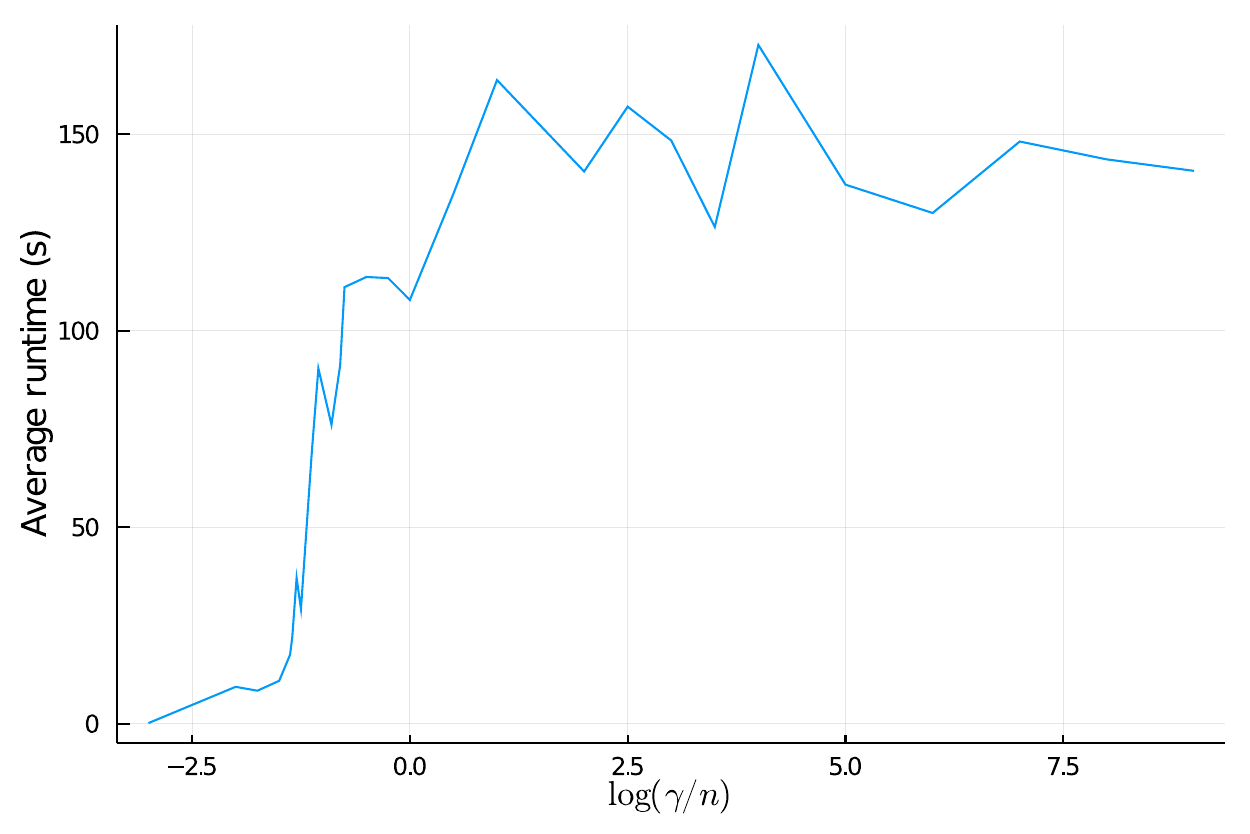}
\end{subfigure} %
\begin{subfigure}[t]{.47\linewidth}
	\centering
	\includegraphics[width=\linewidth]{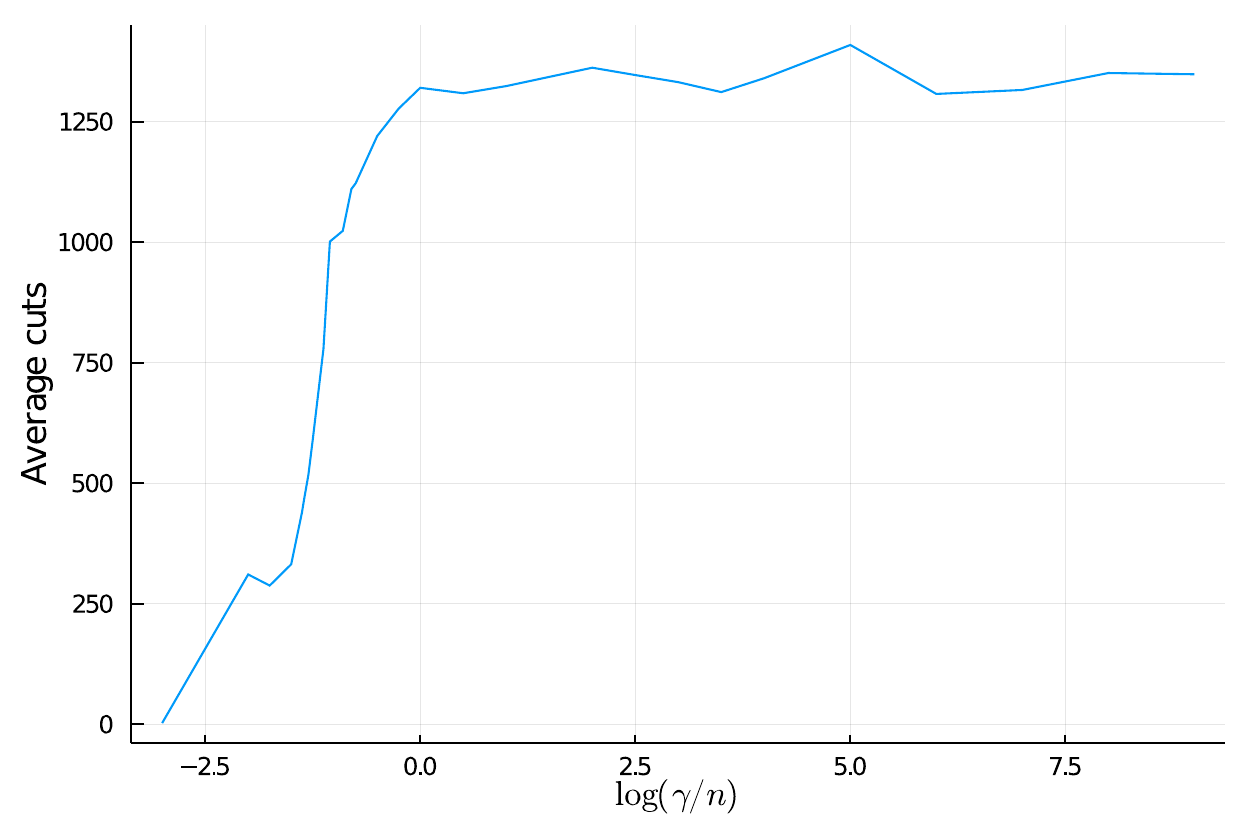}
\end{subfigure}
\caption{{\color{black}Average runtime (left) and number of cuts (right) vs. $\log(\gamma)$ for the $300+$ instances with buy-in and minimum return constraints with a cardinality budget of $k=10$ ($n=300$).}}
\label{fig:gammavariationpard300}
\end{figure}

\begin{figure}[h!]
\centering
\begin{subfigure}[t]{.47\linewidth}
	\centering
	\includegraphics[width=\linewidth]{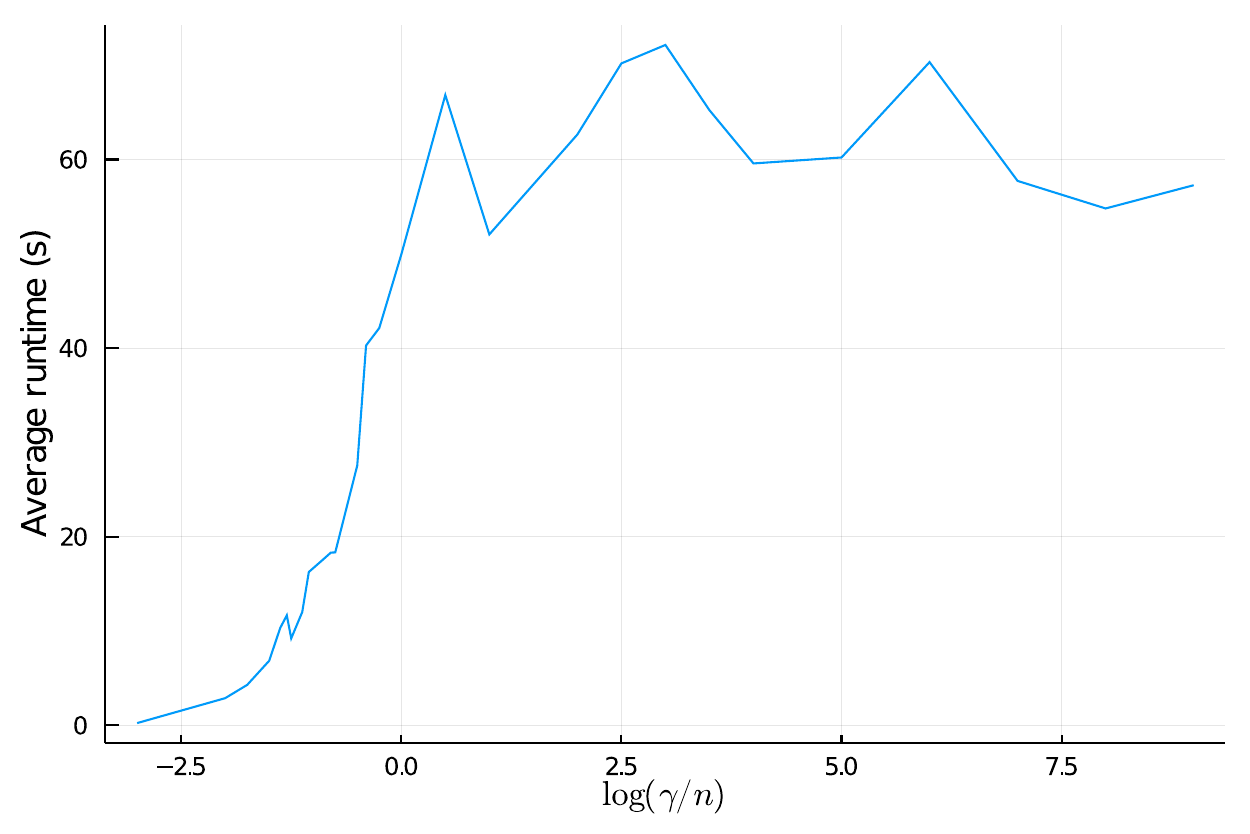}
\end{subfigure} %
\begin{subfigure}[t]{.47\linewidth}
	\centering
	\includegraphics[width=\linewidth]{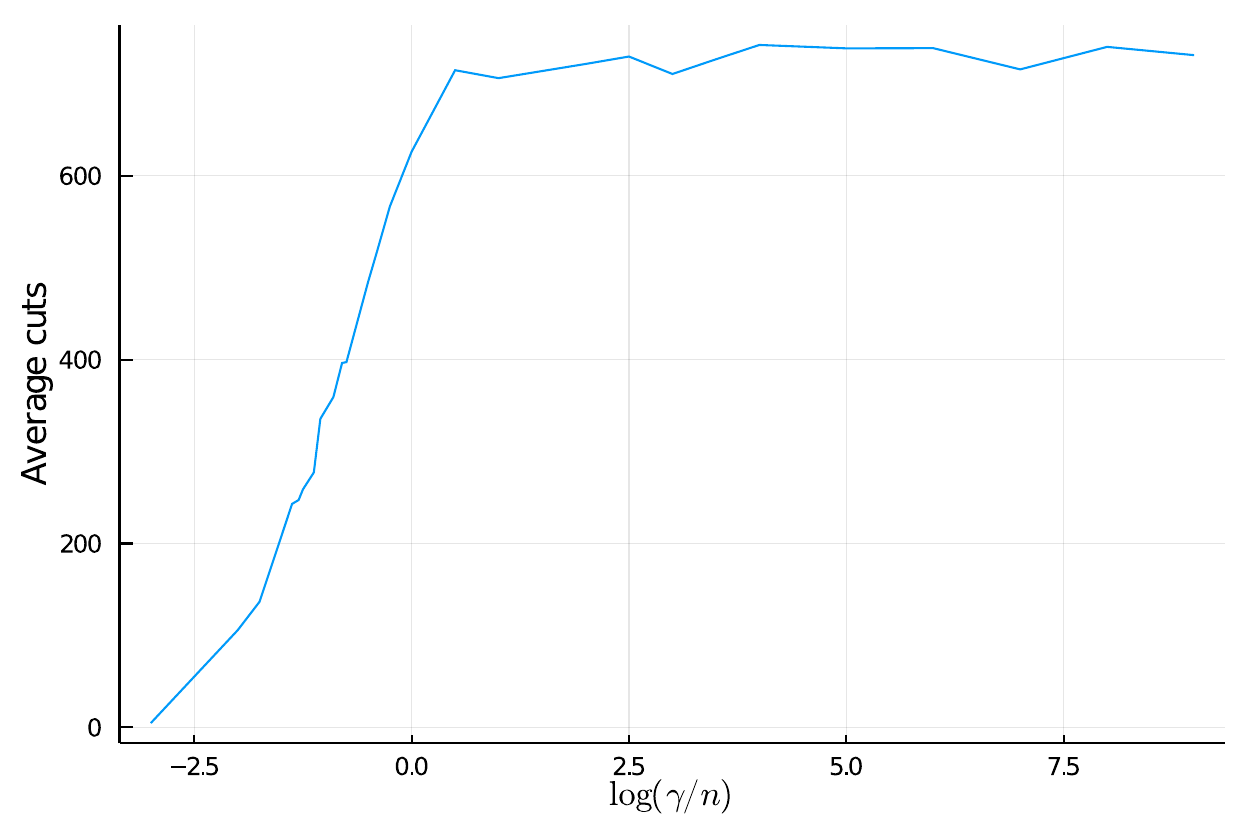}
\end{subfigure}
\caption{{\color{black}Average runtime (left) and number of cuts (right) vs. $\log(\gamma)$ for the $400+$ instances with buy-in and minimum return constraints with a cardinality budget of $k=10$ ($n=400$).}}
\label{fig:gammavariationpard400}
\end{figure}

}

\subsection{Summary of Findings From Numerical Experiments}
We are now in a position to answer the four questions introduced at the start of this section.
Our findings are as follows:
\begin{enumerate}
    \item In the absence of complicating constraints, Algorithm \ref{alg:refinedCuttingPlaneMethod} is substantially more efficient than {\color{black}commercial} MIQO solvers such as \verb|CPLEX|. This efficiency improvement can be explained by (a) our ability to generate stronger and more informative lower bounds via dual subproblems, and (b) our dual representation of the problems' subgradients. Indeed, the method did not require more than one second to solve any of the constraint-free problems considered here, although this phenomenon can be partially attributed to the problem data used. {\color{black}Moreover, the solve times are comparable to those reported by other recent methods such as \cite{zheng2014improving, frangioni2016approximated, frangioni2017improving}; note that this comparison is imperfect since the cited approaches were run on different machines and their source code is unavailable, which precludes a truly fair comparison.}
    \item Although imposing complicating constraints, such as minimum investment constraints, slows Algorithm \ref{alg:refinedCuttingPlaneMethod}, the method performs competitively in the presence of these constraints. Moreover, running the \verb|in-out| cutting-plane method at the root node substantially reduces the initial bound gap, and allows the method to supply a certifiably near-optimal (if not optimal) solution in seconds. {\color{black} This suggests that running the \verb|in-out| method at the root node should be considered as a viable and more scalable alternative to existing root node techniques, particularly in the presence of complicating constraints such as minimum investment constraints, or if the cardinality budget is at least $10$ (although it can do more harm than good for easier problems).}
    \item Algorithm \ref{alg:refinedCuttingPlaneMethod} scales to solve real-world problem instances which comprise selecting assets from universes with thousands of securities, such as the Russell $1000$ and the Wilshire $5000$, while existing {\color{black}commercial solvers} such as \verb|CPLEX| either solve these problems much more slowly or do not successfully solve them, because they cannot attain sufficiently strong lower bounds quickly.
    \item Solutions to Problem \eqref{mainproblem} are stable with respect to the hyperparameter $\gamma$.
\end{enumerate}

\section*{Acknowledgments} \footnotesize We are grateful to the anonymous referees of this and previous versions of the paper for insightful comments which improved the quality of the manuscript, Brad Sturt for editorial comments on a previous version of the paper, and Jean Pauphilet for providing a \verb|Julia| implementation of the \verb|in-out| method.
{
\bibliographystyle{informs2014}

}
\normalsize
\begin{APPENDICES}

\section{Omitted Proofs}\label{sec:ommitedproofs}
In this section, we supply the omitted proofs of results stated in the manuscript, in the order in which the results were stated.

\subsection{Proof of Corollary \ref{corr:lcont}}
\proof{Proof of Corollary \ref{corr:lcont}}
{\color{black} This result is a simple consequence of the subgradient inequality for convex functions \citep[see, e.g.,][]{boyd2004convex}. For completeness, we now provide a proof from first principles:
}

The result follows from applying the lower approximation $$f(\hat{\bm{z}}) \geq f(\bm{z})+\bm{g}_{\bm{z}}^\top(\bm{\hat{z}}-\bm{z}), $$re-arranging to yield $$f(\bm{z})-f(\bm{\hat{z}}) \leq -\bm{g}_{\bm{z}}^\top (\bm{\hat{z}}-\bm{z})$$ and invoking Corollary \ref{subgradientvalues} to rewrite the right-hand-side in the desired form. \hfill \Halmos
\endproof



\subsection{Proof of Corollary \ref{corr:exactrecovery}}
\proof{Proof of Corollary \ref{corr:exactrecovery}}
Let there exist some $(\bm{v}^\star, \bm{w}^\star, \bm{\alpha}^\star, \bm{\beta}_l^\star, \bm{\beta}_u^\star, \lambda^\star)$ which solve Problem \eqref{socpbound}, and binary vector $\bm{z} \in \mathcal{Z}_k^n$, such that these two quantities collectively satisfy the conditions encapsulated in Expression \eqref{conditions}. Then, this optimal solution to Problem \eqref{socpbound} provides the following lower bound for Problem \eqref{mainproblem}:
\begin{align*}
    -\frac{1}{2} \bm{\alpha}^{*\top} \bm{\alpha^\star}+\bm{y}^\top \bm{\alpha^\star} +\bm{\beta}_l^{*\top} \bm{l}-\bm{\beta}_u^{*\top} \bm{u}+\lambda^\star-\bm{e}^\top \bm{v^\star} -k t^\star.
\end{align*}
Moreover, let $\bm{\hat{x}}$ be a candidate solution to Problem \eqref{mainproblem} defined by $\hat{x}_i:=\gamma w_i z_i$. Then, $\bm{\hat{x}}$ is feasible for Problem \eqref{mainproblem}, since $\bm{l} \leq A\bm{\hat{x}} \leq \bm{u}$, $\bm{e}^\top \hat{\bm{x}}=1$, $\bm{\hat{x}} \geq \bm{0}$ and $\Vert \hat{\bm{x}} \Vert_0 \leq k$ by Expression \eqref{conditions} and the definition of $\bm{z}$. Additionally, since an optimal choice of $t$ is the $k$th largest value of $\frac{\gamma}{2}w_i^2$, i.e., $\frac{\gamma}{2}w_{[k]}^2$ \citep[see][Lemma 1]{zakeri2014optimization}, at optimality we have that $\bm{e}^\top \bm{v}+kt=\frac{1}{2\gamma}\hat{\bm{x}}^\top \hat{\bm{x}}$. Therefore, Problem \eqref{mainproblem}'s objective when $\bm{x}=\hat{\bm{x}}$  is given by:
\begin{align*}
    -\frac{1}{2} \bm{\alpha}^{*\top} \bm{\alpha^\star}+\bm{y}^\top \bm{\alpha^\star} +\bm{\beta}_l^{*\top} \bm{l}-\bm{\beta}_u^{*\top} \bm{u}+\lambda^\star-\frac{1}{2\gamma}\bm{\hat{x}}^\top \bm{\hat{x}},
\end{align*}
which is less than or equal to Problem \eqref{socpbound}'s objective, since $v_i^\star=0 \quad \forall i \in [n]\ : \ z_i=0$.

Finally, let $\vert w^\star\vert_{[k]}> \vert w^\star\vert_{[k+1]}$ and let $S$ denote the set of indices such that $\vert w^\star_{i}\vert \geq \vert w^\star\vert_{[k]}$. Then, as the primal-dual KKT conditions for max-$k$ norms \citep[see, e.g., ][Lemma 1]{zakeri2014optimization} imply that an optimal choice of $t$ is given by $t^\star=\frac{\gamma}{2}w^{\star 2}_{[k]}$, we can set $t^\star=\frac{\gamma}{2}w^{\star 2}_{[k]}$ without loss of generality. Note that, in general, this choice is not unique. Indeed, any $t \in [\frac{\gamma}{2} w_{[k+1]}^\star, \frac{\gamma}{2} w_{[k]}^\star]$ constitutes an optimal choice \citep{zakeri2014optimization}.

We then have that $v_i^\star=0 \forall i \notin S$, which implies that the constraint $v_i+t \geq \frac{\gamma}{2}w_i^2$ holds strictly for any $i \notin S$. Therefore, the dual multipliers associated with these constraints must take value $0$. But these constraints' dual multipliers are precisely $\bm{z} \in \mathrm{Conv}(\mathcal{Z}_k^n)$, which implies that $z_i=1 \forall i \in S$ gives a valid set of dual multipliers. Moreover, by Equation \eqref{primaldualkkt}, setting $\bm{x}_i=\gamma z_i w_i^\star$ supplies an optimal (and thus feasible) choice of $\bm{x}$ for this fixed $\bm{z}$. Therefore, this primal-dual pair satisfies \eqref{conditions}. \hfill \Halmos
\endproof



\section{Supplementary Experimental Results}
\subsection{Supplementary Results on Problems With Minimum Investment Constraints}\label{sec:suppbuyin}
We now present the instance-wise runtimes (in seconds) for all instances generated by \citet{frangioni2006perspective}, in Tables \ref{tab:comparisonwithmisocp-pard200plus}-\ref{tab:comparisonwithmisocp-pard400plus}.

\begin{table}[h]
\centering\footnotesize
\caption{Performance of the outer-approximation method vs. CPLEX's MISOCO method on the $200^+$ instances generated by \citet{frangioni2006perspective}, with a time budget of $600$s per approach, $\kappa=0$, $\gamma=\frac{1000}{n}$. We run all approaches on one thread. Note that ``nc'' refers to an instance without an explicit cardinality constraint.}
\begin{tabular}{@{}l l r r r r r r r r r r r r@{}} \toprule
Problem & $k$ & \multicolumn{3}{c@{\hspace{0mm}}}{Algorithm \ref{alg:refinedCuttingPlaneMethod}} & \multicolumn{3}{c@{\hspace{0mm}}}{Algorithm \ref{alg:refinedCuttingPlaneMethod} $+$ in-out} & \multicolumn{3}{c@{\hspace{0mm}}}{Algorithm \ref{alg:refinedCuttingPlaneMethod} $+$ in-out $+$ 50} &  \multicolumn{2}{c@{\hspace{0mm}}}{CPLEX MISOCO} \\
\cmidrule(l){3-5} \cmidrule(l){6-8} \cmidrule(l){9-11} \cmidrule(l){12-13}  & & Time & Nodes & Cuts & Time & Nodes & Cuts & Time & Nodes & Cuts & Time & Nodes \\\midrule
pard200-1 & 6 & 0.74 & 674 & 126 & 1.61 & 832 & 162 & 10.80 & 248 & 36 & 28.54 & 33\\
pard200-1 & 8 & 1.10 & 1,038 & 188 & 2.00 & 1,140 & 220 & 12.95 & 488 & 60 & 30.63 & 33\\
pard200-1 & 10 & 8.73 & 8,844 & 654 & 8.58 & 6,582 & 494 & 12.44 & 2,838 & 180 & 65.04 & 69\\
pard200-1 & 12 & 1.37 & 1,902 & 136 & 1.05 & 1,079 & 74 & 7.84 & 532 & 29 & 130.9 & 122\\
pard200-1 & nc & 1.66 & 3,026 & 123 & 1.58 & 1,182 & 101 & 9.68 & 928 & 79 & $>600$& 624\\\midrule
pard200-2 & 6 & 0.13 & 141 & 24 & 0.16 & 81 & 10 & 3.21 & 42 & 5 & 33.26 & 37\\
pard200-2 & 8 & 0.36 & 327 & 59 & 0.28 & 117 & 23 & 5.19 & 64 & 11 & 35.90 & 29\\
pard200-2 & 10 & 4.19 & 6,313 & 317 & 4.07 & 4,449 & 243 & 10.62 & 2,530 & 144 & 278.6 & 331\\
pard200-2 & 12 & 0.65 & 1,953 & 20 & 0.42 & 187 & 8 & 7.94 & 182 & 7 & $>600$& 775\\
pard200-2 & nc & 0.7 & 1,716 & 24 & 0.53 & 233 & 11 & 7.37 & 184 & 10 & $>600$ & 800\\\midrule
pard200-3 & 6 & 0.87 & 904 & 158 & 0.81 & 740 & 96 & 6.55 & 413 & 40 & 103.7 & 103\\
pard200-3 & 8 & 0.67 & 818 & 98 & 0.82 & 671 & 84 & 6.54 & 343 & 40 & 85.76 & 81\\
pard200-3 & 10 & 2.45 & 4,584 & 189 & 1.45 & 1,444 & 90 & 8.28 & 840 & 42 & 210.9 & 215\\
pard200-3 & 12 & 1.71 & 3,034 & 42 & 0.78 & 923 & 23 & 7.75 & 461 & 9 & $>600$& 648\\
pard200-3 & nc & 1.21 & 2,803 & 38 & 0.65 & 781 & 21 & 8.68 & 416 & 9 & $>600$ & 688\\\midrule
pard200-4 & 6 & 1.53 & 2,096 & 262 & 2.31 & 1,740 & 227 & 8.34 & 1,529 & 193 & 230.3 & 248\\
pard200-4 & 8 & 2.73 & 3,820 & 343 & 2.82 & 3,055 & 260 & 9.50 & 1,740 & 139 & 176.1 & 194\\
pard200-4 & 10 & 10.83 & 14,200 & 647 & 10.82 & 11,380 & 522 & 17.9 & 7,912 & 446 & 527.5 & 617\\
pard200-4 & 12 & 0.98 & 2,332 & 22 & 0.32 & 272 & 12 & 7.83 & 226 & 11 & 581.4 & 643\\
pard200-4 & nc & 0.86 & 2,315 & 22 & 0.39 & 251 & 12 & 7.65 & 206 & 11 & 592.7 & 648\\\midrule
pard200-5 & 6 & 0.44 & 225 & 79 & 0.41 & 147 & 49 & 5.06 & 69 & 15 & 33.6 & 31\\
pard200-5 & 8 & 0.66 & 407 & 112 & 0.53 & 253 & 65 & 5.71 & 93 & 16 & 36.0 & 34\\
pard200-5 & 10 & 2.86 & 3,577 & 322 & 1.76 & 1,644 & 149 & 7.62 & 453 & 48 & 84.98 & 79\\
pard200-5 & 12 & 135.9 & 171,500 & 722 & 12.13 & 10,700 & 294 & 18.45 & 6,098 & 171 & $>600$& 686\\
pard200-5 & nc & 120.4 & 131,100 & 777 & 15.16 & 11,960 & 285 & 17.06 & 5,866 & 142 & $>600$& 818\\\midrule
pard200-6 & 6 & 7.15 & 4,635 & 933 & 7.44 & 5,148 & 902 & 16.1 & 4,875 & 675 & 172.2 & 199\\
pard200-6 & 8 & 6.05 & 5,985 & 777 & 8.38 & 5,885 & 733 & 12.25 & 4,120 & 349 & 112.6 & 135\\
pard200-6 & 10 & 2.64 & 2,305 & 283 & 2.05 & 1,172 & 206 & 7.48 & 514 & 107 & 82.61 & 103\\
pard200-6 & 12 & 1.08 & 1,934 & 81 & 0.54 & 461 & 37 & 7.18 & 409 & 23 & 189.0 & 211\\
pard200-6 & nc & 1.1 & 1,737 & 76 & 0.74 & 799 & 48 & 8.22 & 483 & 32 & $>600$& 700\\\midrule
pard200-7 & 6 & 0.64 & 687 & 122 & 0.74 & 602 & 97 & 6.77 & 291 & 54 & 112.7 & 119\\
pard200-7 & 8 & 0.36 & 431 & 59 & 0.32 & 207 & 31 & 7.12 & 88 & 16 & 59.57 & 65\\
pard200-7 & 10 & 2.21 & 3,570 & 216 & 1.15 & 1,205 & 105 & 8.13 & 640 & 46 & 83.51 & 75\\
pard200-7 & 12 & 12.03 & 15,000 & 76 & 1.17 & 1,185 & 20 & 9.17 & 725 & 13 & $>600$ & 648\\
pard200-7 & nc & 8.77 & 12,930 & 72 & 1.32 & 1,464 & 23 & 9.39 & 841 & 16 & $>600$ & 802\\\midrule
pard200-8 & 6 & 0.2 & 97 & 43 & 0.19 & 75 & 25 & 2.57 & 41 & 11 & 20.55 & 19\\
pard200-8 & 8 & 0.36 & 199 & 68 & 0.51 & 124 & 36 & 3.37 & 47 & 12 & 32.90 & 30\\
pard200-8 & 10 & 3.21 & 3,635 & 295 & 0.78 & 442 & 88 & 7.26 & 151 & 20 & 40.55 & 37\\
pard200-8 & 12 & 96.29 & 82,400 & 581 & 3.09 & 2,237 & 171 & 8.71 & 1,080 & 56 & 185.2 & 200\\
pard200-8 & nc & 45.68 & 66,790 & 574 & 2.96 & 2,455 & 160 & 10.21 & 1,999 & 67 & $>600$& 964\\\midrule
pard200-9 & 6 & 2.6 & 2,404 & 390 & 2.62 & 2,262 & 338 & 7.75 & 1,211 & 110 & 79.61 & 91\\
pard200-9 & 8 & 5.62 & 5,052 & 657 & 5.83 & 3,814 & 540 & 10.19 & 2,093 & 289 & 108.1 & 136\\
pard200-9 & 10 & 3.76 & 3,582 & 403 & 3.09 & 1,817 & 259 & 9.53 & 754 & 125 & 65.21 & 82\\
pard200-9 & 12 & 1.98 & 2,535 & 147 & 0.52 & 296 & 42 & 7.76 & 134 & 10 & 23.70 & 23\\
pard200-9 & nc & 1.96 & 2,473 & 148 & 1.65 & 1,675 & 95 & 9.89 & 899 & 78 & $>600$& 675\\\midrule
pard200-10 & 6 & 1.2 & 1,122 & 226 & 1.42 & 992 & 188 & 6.86 & 385 & 41 & 63.01 & 73\\
pard200-10 & 8 & 1.62 & 1,599 & 242 & 1.48 & 992 & 178 & 6.91 & 415 & 41 & 56.54 & 61\\
pard200-10 & 10 & 36.51 & 25,450 & 1,771 & 9.51 & 6,730 & 833 & 14.27 & 4,025 & 597 & 178.0 & 232\\
pard200-10 & 12 & 3.8 & 5,711 & 211 & 0.58 & 300 & 35 & 7.78 & 152 & 10 & 20.48 & 25\\
pard200-10 & nc & 4.75 & 7,010 & 230 & 2.93 & 2,085 & 164 & 11.84 & 2115 & 117 & $>600$ & 632\\
\bottomrule
\end{tabular}
\label{tab:comparisonwithmisocp-pard200plus}
\end{table}

\begin{table}[h]
\centering\footnotesize
\caption{Performance of the outer-approximation method vs. CPLEX's MISOCO method on the $300^+$ instances generated by \citet{frangioni2006perspective}, with a time budget of $600$s per approach, with $\kappa=0$, $\gamma=\frac{1000}{n}$. We run all approaches on one thread. Note that ``nc'' refers to an instance without an explicit cardinality constraint.}
\begin{tabular}{@{}l l r r r r r r r r r r r r@{}} \toprule
Problem & $k$ & \multicolumn{3}{c@{\hspace{0mm}}}{Algorithm \ref{alg:refinedCuttingPlaneMethod}} & \multicolumn{3}{c@{\hspace{0mm}}}{Algorithm \ref{alg:refinedCuttingPlaneMethod} $+$ in-out} & \multicolumn{3}{c@{\hspace{0mm}}}{Algorithm \ref{alg:refinedCuttingPlaneMethod} $+$ in-out $+$ 50} &  \multicolumn{2}{c@{\hspace{0mm}}}{CPLEX MISOCO} \\
\cmidrule(l){3-5} \cmidrule(l){6-8} \cmidrule(l){9-11} \cmidrule(l){12-13}  & & Time & Nodes & Cuts & Time & Nodes & Cuts & Time & Nodes & Cuts & Time & Nodes \\\midrule
pard300-1 & 6 & 36.56 & 15,870 & 1,974 & 68.44 & 14,130 & 1,864 & 67.68 & 10,920 & 1,295 & $>600$ & 210\\
pard300-1 & 8 & 158.8 & 39,650 & 3,412 & 238.9 & 37,830 & 3,320 & 193.8 & 29,560 & 2,508 & $>600$& 190\\
pard300-1 & 10 & 108.3 & 60,560 & 2,593 & 94.63 & 46,350 & 2,053 & 59.74 & 23,270 & 1,243 & $>600$& 230\\
pard300-1 & 12 & 17.93 & 17,070 & 523 & 8.23 & 5,390 & 219 & 16.30 & 1,567 & 73 & 261.3 & 101\\
pard300-1 & nc & 33.85 & 29,030 & 483 & 26.44 & 15,490 & 419 & 44.41 & 15,060 & 418 & $>600$& 206\\\midrule
pard300-2 & 6 & 13.87 & 7,583 & 935 & 22.34 & 6,805 & 819 & 28.03 & 4,563 & 496 & 346.5 & 117\\
pard300-2 & 8 & 37.76 & 28,910 & 1,962 & 68.53 & 21,470 & 1,753 & 66.81 & 16,060 & 1,254 & 583.0 & 233\\
pard300-2 & 10 & 64.51 & 41,320 & 2,247 & 44.76 & 30,230 & 1,237 & 439.9 & 15,460 & 658 & 562.7 & 216\\
pard300-2 & 12 & 2.76 & 4,355 & 115 & 0.77 & 532 & 26 & 12.22 & 194 & 4 & 172.3 & 57\\
pard300-2 & nc & 4.91 & 7,423 & 123 & 1.91 & 1,440 & 63 & 13.29 & 739 & 38 & $>600$ & 250\\\midrule
pard300-3 & 6 & 33.24 & 21,540 & 1,205 & 47.62 & 20,050 & 1,137 & 52.85 & 13,420 & 793 & $>600$ & 210\\
pard300-3 & 8 & 34.00 & 37,920 & 1,410 & 52.85 & 35,640 & 1,293 & 38.12 & 21,620 & 668 & $>600$ & 206\\
pard300-3 & 10 & 254.3 & 128,500 & 3,526 & 283.2 & 103,700 & 3,114 & 288.0 & 112,000 & 2,502 & $>600$ & 210\\
pard300-3 & 12 & 81.05 & 59,470 & 342 & 4.5 & 3,607 & 84 & 18.83 & 1,144 & 47 & $>600$ & 295\\
pard300-3 & nc & 77.34 & 58,550 & 328 & 8.26 & 5,867 & 116 & 25.16 & 5,137 & 74 & $>600$ & 206\\\midrule
pard300-4 & 6 & 51.46 & 18,810 & 2,255 & 55.49 & 18,400 & 1,849 & 64.77 & 15,930 & 1,539 & $>600$& 225\\
pard300-4 & 8 & 75.72 & 39,830 & 3,192 & 161.3 & 45,040 & 3,108 & 154.5 & 36,930 & 2,611 & $>600$ & 224\\
pard300-4 & 10 & 168.0 & 58,710 & 3,968 & 195.2 & 58,490 & 3,672 & 168.2 & 44,360 & 3,048 & $>600$& 238\\
pard300-4 & 12 & 19.98 & 18,650 & 267 & 8.93 & 5,509 & 187 & 22.75 & 3,707 & 127 & $>600$ & 229\\
pard300-4 & nc & 27.19 & 22,010 & 284 & 16.06 & 12,240 & 215 & 31.03 & 9,941 & 199 & $>600$ & 257\\\midrule
pard300-5 & 6 & 3.99 & 2,670 & 425 & 5.49 & 2,295 & 351 & 12.05 & 934 & 104 & 358.1 & 131\\
pard300-5 & 8 & 6.88 & 5,509 & 491 & 6.53 & 4,227 & 357 & 11.90 & 1,100 & 73 & 192.5 & 68\\
pard300-5 & 10 & 13.5 & 11,890 & 790 & 6.86 & 4,610 & 385 & 14.23 & 1,419 & 140 & 330.3 & 123\\
pard300-5 & 12 & 1.07 & 1,141 & 81 & 0.43 & 174 & 30 & 10.71 & 120 & 18 & 247.7 & 92\\
pard300-5 & nc & 1.05 & 813 & 82 & 1.03 & 511 & 50 & 12.24 & 360 & 31 & $>600$& 224\\\midrule
pard300-6 & 6 & 3.66 & 2,478 & 420 & 4.00 & 2,341 & 353 & 11.47 & 1,089 & 101 & 155.5 & 55\\
pard300-6 & 8 & 10.35 & 8,220 & 771 & 9.71 & 6,503 & 635 & 15.40 & 3,518 & 295 & 307.6 & 110\\
pard300-6 & 10 & 26.95 & 15,450 & 1,151 & 19.09 & 12,090 & 927 & 21.74 & 6,265 & 402 & $>600$& 209\\
pard300-6 & 12 & 6.62 & 7,094 & 275 & 3.74 & 1,118 & 136 & 19.10 & 794 & 79 & 121.4 & 47\\
pard300-6 & nc & 7.65 & 9,391 & 257 & 5.17 & 3,233 & 195 & 19.06 & 3,056 & 158 & $>600$& 214\\\midrule
pard300-7 & 6 & 2.82 & 2,009 & 323 & 3.89 & 1,413 & 285 & 12.86 & 1,120 & 195 & 551.5 & 227\\
pard300-7 & 8 & 4.08 & 4,395 & 337 & 3.3 & 1,996 & 250 & 13.68 & 1,182 & 168 & $>600$& 210\\
pard300-7 & 10 & 5.08 & 5,494 & 334 & 1.77 & 1,716 & 107 & 11.13 & 464 & 26 & 199.0 & 63\\
pard300-7 & 12 & 0.71 & 1,278 & 37 & 0.56 & 455 & 18 & 11.94 & 373 & 14 & 577.8 & 208\\
pard300-7 & nc & 1.43 & 2,940 & 37 & 0.59 & 615 & 13 & 11.53 & 374 & 12 & $>600$ & 200\\\midrule
pard300-8 & 6 & 5.28 & 3,174 & 589 & 6.05 & 3,052 & 549 & 13.53 & 2,232 & 282 & 331.9 & 113\\
pard300-8 & 8 & 11.74 & 7,725 & 1,034 & 15.13 & 7,912 & 983 & 20.73 & 5,125 & 658 & 523.2 & 185\\
pard300-8 & 10 & 23.65 & 18,550 & 1,174 & 12.02 & 8,273 & 602 & 18.47 & 3,646 & 354 & 368.3 & 130\\
pard300-8 & 12 & 7.02 & 8,034 & 331 & 4.33 & 2,786 & 171 & 14.83 & 1,447 & 104 & 234.7 & 89\\
pard300-8 & nc & 8.88 & 9,420 & 333 & 7.09 & 6,558 & 287 & 19.36 & 4,719 & 239 & $>600$ & 220\\\midrule
pard300-9 & 6 & 12.26 & 13,400 & 1,033 & 15.08 & 8,177 & 931 & 20.94 & 4,948 & 620 & 538.5 & 195\\
pard300-9 & 8 & 97.90 & 31,670 & 2,356 & 76.99 & 30,940 & 2,192 & 77.93 & 24,080 & 1,823 & $>600$ & 207\\
pard300-9 & 10 & 215.2 & 97,800 & 2,741 & 118.1 & 64,980 & 1,907 & 66.37 & 36,790 & 1,113 & $>600$ & 201\\
pard300-9 & 12 & 11.08 & 11,750 & 268 & 9.68 & 8,423 & 196 & 18.18 & 3,710 & 117 & $>600$& 269\\
pard300-9 & nc & 24.51 & 25,400 & 284 & 10.65 & 8,604 & 225 & 31.34 & 10,640 & 196 & $>600$ & 240\\\midrule
pard300-10 & 6 & 5.13 & 3,884 & 583 & 7.5 & 3,589 & 503 & 15.05 & 2,227 & 234 & 262.5 & 93\\
pard300-10 & 8 & 9.54 & 6,685 & 803 & 11.44 & 5,257 & 687 & 17.12 & 3,190 & 300 & 289.0 & 107\\
pard300-10 & 10 & 6.02 & 3,417 & 486 & 4.96 & 2,097 & 380 & 14.25 & 1,215 & 229 & 259.5 & 99\\
pard300-10 & 12 & 13.33 & 9,969 & 388 & 5.35 & 3,815 & 207 & 14.56 & 1,690 & 84 & $>600$ & 195\\
pard300-10 & nc & 26.77 & 16,430 & 410 & 15.2 & 8,326 & 336 & 35.89 & 9,676 & 319 & $>600$ & 175\\
\bottomrule
\end{tabular}
\label{tab:comparisonwithmisocp-pard300plus}
\end{table}

\begin{table}[h]
\centering\footnotesize
\caption{Performance of the outer-approximation method vs. CPLEX's MISOCO method on the $400^+$ instances generated by \citet{frangioni2006perspective}, with a time budget of $600$s per approach, with $\kappa=0$, $\gamma=\frac{1000}{n}$. We run all approaches on one thread. Note that ``nc'' refers to an instance without an explicit cardinality constraint. }
\begin{tabular}{@{}l l r r r r r r r r r r r r@{}} \toprule
Problem & $k$ & \multicolumn{3}{c@{\hspace{0mm}}}{Algorithm \ref{alg:refinedCuttingPlaneMethod}} & \multicolumn{3}{c@{\hspace{0mm}}}{Algorithm \ref{alg:refinedCuttingPlaneMethod} $+$ in-out} & \multicolumn{3}{c@{\hspace{0mm}}}{Algorithm \ref{alg:refinedCuttingPlaneMethod} $+$ in-out $+$ 50} &  \multicolumn{2}{c@{\hspace{0mm}}}{CPLEX MISOCO} \\
\cmidrule(l){3-5} \cmidrule(l){6-8} \cmidrule(l){9-11} \cmidrule(l){12-13}  & & Time & Nodes & Cuts & Time & Nodes & Cuts & Time & Nodes & Cuts & Time & Nodes \\\midrule
pard400-1 & 6 & 64.88 & 18,730 & 2,963 & 69.91 & 17,430 & 2,670 & 86.00 & 17,793 & 2,320 & $>600$ & 95\\
pard400-1 & 8 & 364.5 & 54,420 & 5,400 & 283.6 & 66,130 & 4,734 & 232.7 & 41,000 & 4,086 & $>600$ & 94\\
pard400-1 & 10 & 36.33 & 24,850 & 980 & 13.25 & 8,578 & 554 & 24.45 & 3,740 & 318 & $>600$ & 97\\
pard400-1 & 12 & 14.19 & 11,480 & 328 & 9.40 & 6,732 & 214 & 24.06 & 4,030 & 144 & $>600$ & 100\\
pard400-1 & nc & 49.56 & 35,030 & 336 & 17.41 & 11,580 & 261 & 44.27 & 16,000 & 238 & $>600$ & 74\\\midrule
pard400-2 & 6 & 0.18 & 71 & 21 & 0.12 & 24 & 8 & 3.31 & 18 & 6 & 160.9 & 17\\
pard400-2 & 8 & 0.31 & 227 & 24 & 0.13 & 12 & 8 & 2.17 & 14 & 8 & 350.9 & 53\\
pard400-2 & 10 & 0.14 & 87 & 10 & 0.05 & 0 & 2 & 0.05 & 0 & 2 & 178.2 & 23\\
pard400-2 & 12 & 0.12 & 54 & 6 & 0.24 & 10 & 3 & 1.13 & 5 & 3 & $>600$ & 69\\
pard400-2 & nc & 0.12 & 52 & 5 & 0.25 & 9 & 2 & 2.57 & 12 & 2 & $>600$ & 70\\\midrule
pard400-3 & 6 & 1.38 & 1,335 & 149 & 1.70 & 895 & 121 & 16.37 & 534 & 84 & $>600$ & 100\\
pard400-3 & 8 & 3.22 & 2,458 & 271 & 3.24 & 1,862 & 206 & 18.40 & 1,353 & 133 & $>600$ & 80\\
pard400-3 & 10 & 8.81 & 9,924 & 347 & 3.22 & 2,500 & 146 & 19.15 & 1,351 & 55 & $>600$ & 82\\
pard400-3 & 12 & 0.45 & 838 & 12 & 0.26 & 102 & 2 & 13.71 & 59 & 2 & 582.0 & 92\\
pard400-3 & nc & 0.56 & 1,259 & 10 & 0.22 & 130 & 2 & 15.34 & 74 & 2 & $>600$ & 100\\\midrule
pard400-4 & 6 & 53.04 & 18,490 & 1,712 & 54.56 & 14,360 & 1,677 & 57.55 & 10,890 & 1,237 & $>600$ & 99\\
pard400-4 & 8 & 183.9 & 47,460 & 3,660 & 179.0 & 42,140 & 3,522 & 166.1 & 37,070 & 2,923 & $>600$ & 90\\
pard400-4 & 10 & 259.1 & 76,390 & 2,153 & 516.5 & 81,900 & 5,782 & 439.3 & 96,750 & 3,614 & $>600$ & 90\\
pard400-4 & 12 & 1.88 & 2,428 & 98 & 0.64 & 311 & 21 & 15.56 & 220 & 12 & 407.6 & 67\\
pard400-4 & nc & 3.76 & 4,738 & 105 & 2.14 & 1,795 & 66 & 18.33 & 1,008 & 53 & $>600$ & 90\\\midrule
pard400-5 & 6 & 11.07 & 5,100 & 658 & 11.41 & 4,793 & 586 & 23.78 & 3,363 & 363 & $>600$ & 94\\
pard400-5 & 8 & 17.42 & 12,060 & 939 & 20.97 & 9,459 & 811 & 35.33 & 6,300 & 507 & $>600$& 94\\
pard400-5 & 10 & 213.7 & 74,590 & 2,312 & 175.9 & 73,740 & 1,933 & 100.6 & 42,380 & 1,184 & $>600$ & 89\\
pard400-5 & 12 & 9.29 & 9,720 & 272 & 4.13 & 2,538 & 105 & 19.9 & 765 & 59 & 485.4 & 95\\
pard400-5 & nc & 17.16 & 15,750 & 306 & 19.33 & 12,320 & 235 & 38.46 & 9,137 & 244 & $>600$ & 70\\\midrule
pard400-6 & 6 & 0.27 & 116 & 30 & 0.25 & 32 & 9 & 6.34 & 37 & 6 & 356.2 & 61\\
pard400-6 & 8 & 0.17 & 92 & 15 & 0.06 & 0 & 2 & 0.08 & 0 & 2 & 208.4 & 33\\
pard400-6 & 10 & 0.42 & 317 & 36 & 0.18 & 44 & 9 & 4.81 & 18 & 2 & 200.9 & 33\\
pard400-6 & 12 & 4.8 & 7,491 & 76 & 0.91 & 545 & 22 & 19.21 & 324 & 15 & $>600$ & 97\\
pard400-6 & nc & 5.4 & 8,154 & 82 & 1.36 & 654 & 17 & 19.44 & 372 & 15 & $>600$ & 83\\\midrule
pard400-7 & 6 & 48.05 & 16,100 & 2,514 & 90.72 & 12,480 & 2,376 & 122.5 & 14,130 & 2,233 & $>600$ & 98\\
pard400-7 & 8 & 114.0 & 39,650 & 3,412 & 178.3 & 30,110 & 3,224 & 194.2 & 27,800 & 2,819 & $>600$ & 86\\
pard400-7 & 10 & 31.51 & 21,060 & 1,304 & 35.49 & 19,670 & 985 & 32.94 & 8,230 & 497 & $>600$ & 86\\
pard400-7 & 12 & 1.38 & 1,567 & 70 & 0.42 & 164 & 13 & 12.63 & 60 & 5 & 169.5 & 25\\
pard400-7 & nc & 1.77 & 2,063 & 83 & 2.09 & 1,410 & 64 & 17.73 & 893 & 42 & $>600$ & 61\\\midrule
pard400-8 & 6 & 118.1 & 27,150 & 3,187 & 165.0 & 28,200 & 3,025 & 185.1 & 24,550 & 2,753 & $>600$ & 98\\
pard400-8 & 8 & 342.8 & 91,060 & 5,120 & 335.7 & 62,250 & 5,377 & 356.2 & 58,570 & 4,889 & $>600$ & 97\\
pard400-8 & 10 & 229.3 & 113,200 & 2,704 & 105.9 & 60,640 & 1,546 & 86.51 & 36,100 & 1,100 & $>600$ & 91\\
pard400-8 & 12 & 3.14 & 3,999 & 100 & 1.31 & 948 & 31 & 19.45 & 248 & 16 & 375.6 & 55\\
pard400-8 & nc & 3.14 & 3,717 & 92 & 4.26 & 3,786 & 78 & 22.05 & 1,974 & 63 & $>600$ & 57\\\midrule
pard400-9 & 6 & 77.79 & 22,580 & 2,345 & 103.5 & 20,500 & 2,242 & 107.7 & 17,480 & 1,788 & $>600$ & 88\\
pard400-9 & 8 & 466.0 & 60,570 & 4,064 & 227.1 & 63,950 & 3,900 & 217.7 & 54,390 & 3,298 & $>600$ & 89\\
pard400-9 & 10 & 409.0 & 126,200 & 3,610 & 16.71 & 8,440 & 448 & 34.38 & 6,406 & 315 & $>600$ & 77\\
pard400-9 & 12 & 0.69 & 747 & 52 & 0.72 & 238 & 33 & 14.43 & 212 & 20 & $>600$ & 96\\
pard400-9 & nc & 0.61 & 629 & 43 & 0.77 & 458 & 33 & 14.5 & 379 & 22 & $>600$ & 100\\\midrule
pard400-10 & 6 & 170.0 & 23,610 & 3,587 & 168.1 & 22,860 & 3,473 & 227.1 & 21,270 & 3,172 & $>600$ & 90\\
pard400-10 & 8 & 245.2 & 45,870 & 5,375 & 380.3 & 53,370 & 5,307 & 410.4 & 53,740 & 4,974 & $>600$ & 92\\
pard400-10 & 10 & 391.6 & 108,400 & 3,236 & 177.5 & 67,620 & 2,292 & 72.6 & 26,350 & 1,162 & $>600$ & 80\\
pard400-10 & 12 & 3.79 & 4,910 & 152 & 1.01 & 557 & 42 & 16.65 & 351 & 22 & 360.0 & 57\\
pard400-10 & nc & 4.70 & 4,035 & 143 & 4.08 & 3253 & 130 & 20.43 & 2,239 & 113 & $>600$ & 37\\
\bottomrule
\end{tabular}
\label{tab:comparisonwithmisocp-pard400plus}
\end{table}

We now present the instance-wise runtimes (in seconds) for the smallest instances generated by \citet{frangioni2006perspective}, without any diagonal matrix extraction, and $k \in \{6,8,10,12,n\}$. Table \ref{tab:comparisonwithmisocp-noddterm} demonstrates that not using any diagonal matrix extraction technique substantially slows our approach.

\begin{table}[h]
\centering\footnotesize
\caption{Performance of the outer-approximation method on the $200^+$ instances generated by \citet{frangioni2006perspective}, with a time budget of $600$s per approach, $\kappa=0$, $\gamma=\frac{1000}{n}$, and no diagonal matrix extraction. We run all approaches on one thread. Note that ``nc'' refers to an instance without an explicit cardinality constraint.}
\begin{tabular}{@{}l l r r r r r r r r r r r r@{}} \toprule
Problem & $k$ & \multicolumn{3}{c@{\hspace{0mm}}}{Algorithm \ref{alg:refinedCuttingPlaneMethod}} & \multicolumn{3}{c@{\hspace{0mm}}}{Algorithm \ref{alg:refinedCuttingPlaneMethod}$+$in-out} & \multicolumn{3}{c@{\hspace{0mm}}}{Algorithm \ref{alg:refinedCuttingPlaneMethod} $+$ in-out $+$ 50} & \multicolumn{3}{c@{\hspace{0mm}}}{CPLEX MISOCO}\\
\cmidrule(l){3-5} \cmidrule(l){6-8} \cmidrule(l){9-11} \cmidrule(l){12-13}  & & Time & Nodes & Cuts & Time & Nodes & Cuts & Time & Nodes & Cuts & Time & Nodes\\\midrule
pard200-1 & 6 & $>600$ & 143,800 & 10,650 & $>600$ & 217,200 & 9,671 & $>600$ & 120,900 & 7,530 & $>600$ & 700\\
pard200-1 & 8 & $>600$ & 94,900 & 11,830 & 400.6 & 118,300 & 6,461 & $>600$ & 106,900 & 7,822 & $>600$ & 600\\
pard200-1 & 10 & $>600$ & 173,000 & 10,020 & $>600$ & 71,050 & 10,340 & $>600$ & 63,470 & 5,422 & $>600$ & 700\\
pard200-1 & 12 & $>600$ & 223,700 & 6,852 & $>600$ & 68,600 & 10,500 & $>600$ & 42,450 & 6,632 & $>600$ & 686\\
pard200-1 & nc & $>600$ & 286,600 & 8,236 & $>600$ & 59,070 & 12,560 & $>600$ & 53,580 & 8,221 & $>600$ & 800\\\midrule
pard200-2 & 6 & $>600$ & 243,600 & 10,930 & 239.9 & 47,180 & 8,305 & $>600$ & 90,460 & 8,538 & $>600$ & 500\\
pard200-2 & 8 & $>600$ & 328,100 & 9,365 & $>600$ & 234,100 & 8,547 & $>600$ & 207,800 & 2,011 & $>600$ & 242\\
pard200-2 & 10 & $>600$ & 468,400 & 9,134 & 0.2 & 47 & 12 & 5.32 & 28 & 16 & $>600$ & 700\\
pard200-2 & 12 & $>600$ & 273,300 & 7,406 & 0.47 & 3 & 18 & 1.12 & 3 & 18 & $>600$ & 700\\
pard200-2 & nc & $>600$ & 505,000 & 8,497 & 0.21 & 65 & 12 & 9.95 & 69 & 12 & $>600$ & 600\\\midrule
pard200-3 & 6 & $>600$ & 112,800 & 12,020 & 1.36 & 236 & 32 & 43.69 & 515 & 32 & $>600$ & 505\\
pard200-3 & 8 & $>600$ & 235,000 & 10,150 & 1.24 & 385 & 24 & 59.49 & 750 & 24 & $>600$ & 500\\
pard200-3 & 10 & $>600$ & 245,500 & 6,940 & 3.56 & 6,937 & 12 & 136.5 & 44,310 & 18 & $>600$ & 510\\
pard200-3 & 12 & $>600$ & 292,360 & 7,016 & 57.13 & 63,850 & 944 & $>600$ & 139,800 & 935 & $>600$ & 346\\
pard200-3 & nc & $>600$ & 334,300 & 8,382 & 42.01 & 77,870 & 2,384 & 285.73 & 120,900 & 3,660 & $>600$ & 600\\\midrule
pard200-4 & 6 & $>600$ & 86,780 & 13,660 & 0.3 & 0 & 14 & 0.49 & 0 & 14 & $>600$ & 500\\
pard200-4 & 8 & $>600$ & 272,600 & 10,700 & 0.4 & 642 & 20 & 50.41 & 620 & 20 & $>600$ & 561\\
pard200-4 & 10 & $>600$ & 289,900 & 8,193 & 0.61 & 516 & 30 & 46.85 & 879 & 32 & $>600$ & 498\\
pard200-4 & 12 & $>600$ & 303,600 & 6,790 & 0.9 & 228 & 16 & 51.31 & 216 & 16 & $>600$ & 294\\
pard200-4 & nc & $>600$ & 366,100 & 7,999 & 0.38 & 86 & 22 & 13.87 & 80 & 22 & $>600$ & 587\\\midrule
pard200-5 & 6 & $>600$ & 112,000 & 9,616 & $>600$ & 141,200 & 9,208 & $>600$ & 127,100 & 7,060 & $>600$ & 700\\
pard200-5 & 8 & $>600$ & 132,600 & 11,370 & 135.0 & 52,970 & 4,089 & $>600$ & 93,390 & 7,365 & $>600$ & 600\\
pard200-5 & 10 & $>600$ & 183,100 & 9,788 & $>600$ & 51,300 & 10,010 & $>600$ & 68,410 & 6,456 & $>600$ & 700\\
pard200-5 & 12 & $>600$ & 203,500 & 5,519 & $>600$ & 53,730 & 9,813 & $>600$ & 39,300 & 5,365 & $>600$ & 600\\
pard200-5 & nc & $>600$ & 315,400 & 9,270 & $>600$ & 59,070 & 12,680 & $>600$ & 86,240 & 7,025 & $>600$ & 800\\\midrule
pard200-6 & 6 & $>600$ & 116,700 & 11,260 & $>600$ & 181,500 & 8,698 & $>600$ & 131,200 & 6,758 & $>600$ & 691\\
pard200-6 & 8 & $>600$ & 161,100 & 10,700 & $>600$ & 202,060 & 9,698 & $>600$ & 105,400 & 8,449 & $>600$ & 700\\
pard200-6 & 10 & $>600$ & 141,200 & 11,370 & $>600$ & 80,800 & 11,170 & $>600$ & 56,700 & 7,828 & $>600$ & 700\\
pard200-6 & 12 & $>600$ & 236,600 & 7,586 & $>600$ & 52,560 & 11,000 & $>600$ & 41,630 & 5,415 & $>600$ & 400\\
pard200-6 & nc & $>600$ & 338,100 & 9,120 & $>600$ & 71,790 & 13,850 & $>600$ & 50,980 & 8,922 & $>600$ & 800\\\midrule
pard200-7 & 6 & $>600$ & 106,000 & 8,527 & 3.66 & 4,736 & 392 & 72.66 & 6,080 & 391 & $>600$ & 500\\
pard200-7 & 8 & $>600$ & 149,800 & 12,640 & $>600$ & 193,900 & 10,840 & $>600$ & 222,800 & 3,503 & $>600$ & 500\\
pard200-7 & 10 & $>600$ & 214,200 & 9,858 & $>600$ & 139,000 & 8,465 & $>600$ & 154,300 & 5,859 & $>600$ & 485\\
pard200-7 & 12 & $>600$ & 313,000 & 7,902 & $>600$ & 206,900 & 9,568 & $>600$ & 215,100 & 3,773 & $>600$ & 500\\
pard200-7 & nc & $>600$ & 220,800 & 6,976 & $>600$ & 194,500 & 8,507 & $>600$ & 209,600 & 4,971 & $>600$ & 700\\\midrule
pard200-8 & 6 & $>600$ & 167,600 & 11,200 & $>600$ & 217,700 & 9,307 & $>600$ & 173,200 & 7,609 & $>600$ & 650\\
pard200-8 & 8 & $>600$ & 183,200 & 11,140 & $>600$ & 213,900 & 10,100 & $>600$ & 103,100 & 5,954 & $>600$ & 700\\
pard200-8 & 10 & $>600$ & 182,400 & 10,190 & $>600$ & 115,600 & 9,215 & $>600$ & 35,700 & 4,307 & $>600$ & 360\\
pard200-8 & 12 & $>600$ & 217,000 & 6,149 & $>600$ & 52,420 & 10,550 & $>600$ & 52,300 & 9,066 & $>600$ & 700\\
pard200-8 & nc & $>600$ & 241,600 & 6,751 & $>600$ & 56,070 & 11,820 & $>600$ & 50,590 & 7,702 & $>600$ & 471\\\midrule
pard200-9 & 6 & $>600$ & 126,500 & 11,140 & $>600$ & 128,500 & 9,437 & $>600$ & 106,200 & 7,029 & $>600$ & 600\\
pard200-9 & 8 & $>600$ & 130,300 & 10,580 & $>600$ & 135,700 & 8,579 & $>600$ & 88,600 & 6,653 & $>600$ & 700\\
pard200-9 & 10 & $>600$ & 186,500 & 11,160 & $>600$ & 63,560 & 10,570 & $>600$ & 57,800 & 5,646 & $>600$ & 600\\
pard200-9 & 12 & $>600$ & 266,100 & 7,433 & $>600$ & 61,000 & 10,700 & $>600$ & 42,550 & 6,381 & $>600$ & 689\\
pard200-9 & nc & $>600$ & 256,100 & 7,535 & $>600$ & 57,280 & 10,540 & $>600$ & 38,560 & 6,376 & $>600$ & 800\\\midrule
pard200-10 & 6 & $>600$ & 179,600 & 12,700 & $>600$ & 245,900 & 7,980 & $>600$ & 155,100 & 7,300 & $>600$ & 600\\
pard200-10 & 8 & $>600$ & 123,600 & 10,060 & $>600$ & 198,370 & 8,670 & $>600$ & 36,770 & 2,414 & $>600$ & 680\\
pard200-10 & 10 & $>600$ & 164,800 & 8,273 & $>600$ & 46,100 & 9,375 & $>600$ & 29,000 & 4,824 & $>600$ & 393\\
pard200-10 & 12 & $>600$ & 193,300 & 5,776 & $>600$ & 52,820 & 9,616 & $>600$ & 37,900 & 5,834 & $>600$ & 364\\
pard200-10 & nc & $>600$ & 193,700 & 5,523 & $>600$ & 40,770 & 8,887 & $>600$ & 34,300 & 6,096 & $>600$ & 432\\
\bottomrule
\end{tabular}
\label{tab:comparisonwithmisocp-noddterm}
\end{table}

{\color{black}
Next, we present the aggregate runtimes (in seconds) for all instances generated by \citet{frangioni2006perspective}, when we run \verb|CPLEX|'s MISOCP solver after first supplying the cuts generated by the \verb|in-out| method. Note that, to allow \verb|CPLEX| to benefit from the \verb|in-out| cuts, we introduce an auxiliary variable $\tau$, change the objective to minimizing $\tau$, impose the constraint $\tau\geq \frac{1}{2}\bm{x}^\top \bm{\Sigma}\bm{x}+\frac{1}{2\gamma}\bm{e}^\top \bm{\theta}-\kappa \bm{mu}^\top \bm{x}$ to model the MISOCO objective, and impose the cuts from the in-out method using the epigraph variable $\tau$.

\begin{table}[h]
\centering\footnotesize
\caption{\color{black}Average runtime in seconds per approach with $\kappa=0$, $\gamma=\frac{1000}{n}$ for the problems generated by \citet{frangioni2006perspective}. We impose a time limit of $600$s and run all approaches on one thread. If a solver fails to converge, we report the number of explored nodes at the time limit, use $600$s in lieu of the solve time, and report the number of failed instances (out of $10$) next to the solve time in brackets. {\color{black}Note that the minimum investment constraints impose an implicit cardinality constraint with $k \approx {\color{black}13}$.}}
\begin{tabular}{@{}l l r r@{}} \toprule
Problem & $k$ &  \multicolumn{2}{c@{\hspace{0mm}}}{CPLEX MISOCO in-out} \\
\cmidrule(l){3-4} & & Time & Nodes \\\midrule
200+ & 6 & 255.5 (1)& 92.3\\
200+ & 8 & 272.2 (1)& 98.4\\
200+ & 10 & 375.2 (3)& 116.4\\
200+ & 12 & 447.8 (6) & 216\\
200+ & 200 & $>600$ (10) & 332.1\\\midrule
300+ & 6 & 573.4 (9)& 92\\
300+ & 8 & 578.4 (9)& 89.2\\
300+ & 10 & $>600$ (10)& 92.1\\
300+ & 12 & 575.65 (9) & 117.3\\
300+ & 200 & $>600$ (10) & 108.5\\\midrule
400+ & 6 & 569.4 (9)& 43.3\\
400+ & 8 & 563.3 (9)& 45.1\\
400+ & 10 & 562.4 (9)& 45.7\\
400+ & 12 & 593.8 (9) & 58.3\\
400+ & 200 & $>600$ (10) & 41.4\\
\bottomrule
\end{tabular}
\label{tab:comparisonwithminreturn_fullpard_cplexinout}
\end{table}

}

Finally, we present the instance-wise runtimes (in seconds) for the smallest instances generated by \citet{frangioni2006perspective}, with the diagonal matrix extraction technique proposed by \citet{zheng2014improving}, and $k \in \{10,n\}$ (we restrict the values $k$ can take to use the diagonal matrices pre-computed by \citet{frangioni2017improving}). Table \ref{tab:comparisonwithmisocp-pardusingzheng} demonstrates that using the diagonal matrix extraction technique proposed by \citet{zheng2014improving} substantially slows our approach; the results for $n \in \{300, 400\}$ are similar. Indeed, this technique is only faster for the pard$200$-$1$ problem with $k=10$, and is slower in the other $95\%$ of instances (sometimes substantially so).

\begin{table}[h]
\centering\footnotesize
\caption{Performance of the outer-approximation method on the $200^+$ instances generated by \citet{frangioni2006perspective}, with a time budget of $600$s per approach, $\kappa=0$, $\gamma=\frac{1000}{n}$, and the diagonal matrix extraction technique proposed by \citet{zheng2014improving}. We run all approaches on one thread. Note that ``nc'' refers to an instance without an explicit cardinality constraint.}
\begin{tabular}{@{}l l r r r r r r r r r r@{}} \toprule
Problem & $k$ & \multicolumn{3}{c@{\hspace{0mm}}}{Algorithm \ref{alg:refinedCuttingPlaneMethod}} & \multicolumn{3}{c@{\hspace{0mm}}}{Algorithm \ref{alg:refinedCuttingPlaneMethod} $+$ in-out} & \multicolumn{3}{c@{\hspace{0mm}}}{Algorithm \ref{alg:refinedCuttingPlaneMethod} $+$ in-out $+$ 50} \\
\cmidrule(l){3-5} \cmidrule(l){6-8} \cmidrule(l){9-11}  & & Time & Nodes & Cuts & Time & Nodes & Cuts & Time & Nodes & Cuts \\\midrule
pard200-1 & 10 & 0.74 & 130 & 30 & 0.03 & 0 & 4 & 0.03 & 0 & 4\\
pard200-1 & nc & $>600$ & 887,400 & 1,386 & $>600$  & 539,900 & 1,070 & $>600$ & 369,500 & 675\\\midrule
pard200-2 & 10 & 234.3 & 239,500 & 875 & 57.14 & 45,230 & 193 & 78.88 & 39,000 & 196\\
pard200-2 & nc & $>600$ & 995,900 & 823 & $>600$  & 157,100 & 135 & $>600$  & 226,100 & 66\\\midrule
pard200-3 & 10 & 245.1 & 207,200 & 1,195 & 71.95 & 55,050 & 365 & 76.64 & 30,910 & 249\\
pard200-3 & nc & $>600$ & 903,500 & 888 & $>600$ & 357,200 & 246 & $>600$  & 268,400 & 259\\\midrule
pard200-4 & 10 & $>600$  & 442,500 & 1,967 & 344.7 & 223,700 & 1,053 & 228.9 & 135,500 & 760\\
pard200-4 & nc & 535.1 & 913,500 & 1,092 & $>600$  & 297,600 & 94 & 529.8 & 212,600 & 94\\\midrule
pard200-5 & 10 & $>600$  & 439,900 & 4,965 & 48.87 & 70,200 & 206 & 69.71 & 52,300 & 204\\
pard200-5 & nc & $>600$  & 1,340,000 & 1,314 & $>600$  & 531,400 & 1,408 & $>600$  & 573,200 & 1,358\\\midrule
pard200-6 & 10 & $>600$  & 311,800 & 4,922 & 6.54 & 6,382 & 116 & 36.29 & 12,370 & 107\\
pard200-6 & nc & $>600$ & 1,280,000 & 1,016 & $>600$  & 479,900 & 789 & $>600$  & 557,600 & 580\\\midrule
pard200-7 & 10 & 549.8 & 389,000 & 2,542 & 515.6 & 228,100 & 1,336 & 292.4 & 105,900 & 743\\
pard200-7 & nc & $>600$ & 1,245,000 & 522 & $>600$  & 496,200 & 183 & $>600$  & 502,600 & 119\\\midrule
pard200-8 & 10 & $>600$  & 399,200 & 3,419 & 2.32 & 1,638 & 46 & 20.19 & 1,716 & 45\\
pard200-8 & nc & $>600$ & 1,337,000 & 862 & $>600$  & 674,300 & 552 & $>600$  & 507,900 & 420\\\midrule
pard200-9 & 10 & 589.7 & 576,100 & 1,756 & 6.31 & 8,746 & 122 & 26.53 & 8,290 & 143\\
pard200-9 & nc & $>600$  & 1,264,000 & 1,941 & $>600$  & 703,900 & 1,977 & $>600$  & 498,000 & 1,970\\\midrule
pard200-10 & 10 & $>600$  & 416,200 & 3,798 & 288.4 & 160,300 & 1,070 & 422.0 & 192,800 & 1,002\\
pard200-10 & nc & $>600$  & 974,400 & 1,584 & $>600$  & 498,100 & 1,513 & $>600$  & 313,400 & 1482\\
\bottomrule
\end{tabular}
\label{tab:comparisonwithmisocp-pardusingzheng}
\end{table}
\newpage
\FloatBarrier
\section{Additional Pseudocode}\label{sec:auxpsueodocode}

In this appendix, we provide auxiliary pseudocode pertaining to the experiments run in Section \ref{sec:compexperiments}. Specifically, we provide pseudocode pertaining to our implementation of the \verb|in-out| method of \cite{ben2007acceleration}, which we have applied before running Algorithm \ref{alg:refinedCuttingPlaneMethod} in some problems in Section \ref{sec:compexperiments}.

\begin{algorithm*}
\caption{The in-out method of \cite{ben2007acceleration}, as applied at the root node.}
\label{alg:refinedCuttingPlaneMethodinout}
\begin{algorithmic}
\REQUIRE Optimal solution to Problem \eqref{mainproblemmisocp_dualversion} $\bm{z}^\star$, objective value $\theta_{\text{soco}}$
\STATE $\epsilon \gets 10^{-10}, \lambda \gets 0.1, \delta \gets 2 \epsilon$
\STATE $t \leftarrow 1 $
\REPEAT
\STATE Compute $\bm{z}_0, \theta_0$ solution of
\begin{align*}
\min_{\bm{z} \in \mathrm{Conv}\left(\mathcal{Z}_k^n\right), \theta} \: \theta \quad \mbox{ s.t. } \quad \theta \geq f(\bm{z}_i) + g_{\bm{z}_i}^\top (\bm{z}-\bm{z}_i) \quad \forall i \in [t].
\end{align*}
\IF{$z_0$ has not improved for $5$ consecutive iterations}
\STATE Set $\lambda=1$
\IF{$z_0$ has not improved for $10$ consecutive iterations}
\STATE Set $\delta=0$
\ENDIF
\ENDIF
\STATE Set $z_{t+1} \gets \lambda z_0+(1-\lambda)z_{\text{soco}}+\delta \bm{e}$.
\STATE Round $\bm{z}_{t+1}$ coordinate-wise so that $\bm{z}_{t+1} \in [0, 1]^n$.
\STATE Compute $f(\bm{z}_{t+1})$ and $g_{\bm{z}_{t+1}} \in \partial f (\bm{z}_{t+1})$.
\STATE Apply cut $\theta \geq f(\bm{z}_{t+1}) + g_{\bm{z}_{t+1}}^\top (\bm{z}-\bm{z}_{t+1})$ at root node of integer model.
\STATE $t \leftarrow t+1 $
\UNTIL{$ f(\bm{z}_0)-\theta_0\leq \varepsilon$ or $t>200$}
\RETURN $\bm{z}_t$
\end{algorithmic}
\end{algorithm*}

\end{APPENDICES}
\end{document}